\documentclass[twoside,a4paper]{amsart}

\usepackage[colorlinks=true, breaklinks=true, urlcolor=webbrown, linkcolor=RoyalBlue, citecolor=webgreen, backref=page]{hyperref}
\usepackage[nochapters,beramono,eulermath,pdfspacing,dottedtoc,subfig]{classicthesis}
\usepackage{epsfig, graphicx, subfigure}
\usepackage{verbatim, setspace}
\usepackage{amsmath, amssymb}

\rohead{\mbox{\color{halfgray} \rightmark\hfil \hspace{0.5em} \rlap{\small \vline \kern1em\color{black}\thepage}}}
\lehead{\mbox{\llap{\small\thepage\kern1em\color{halfgray} \vline}\color{halfgray}\hspace{0.5em}\rightmark\hfil}} 
\numberwithin{equation}{section}

\newcommand{\R}		{\mathbb{R}}
\newcommand{\C}		{\mathbb{C}}
\newcommand{\N}		{\mathbb{N}}
\newcommand{\Z}		{\mathbb{Z}}

\newcommand{\dist}{\mathsf{dist}}

\newcommand{\supp}{\mathsf{supp}}
\newcommand{\diag}{\mathsf{diag}}
\newcommand{\re}{\mathsf{Re}}
\newcommand{\im}{\mathsf{Im}}
\newcommand{\sgn}{\mathsf{sgn}}
\newcommand{\n}	{\textnormal{\bf n}}
\renewcommand{\arg}{\mathsf{arg}}
\renewcommand{\det}{\mathsf{det}}
\renewcommand{\deg}{\mathsf{deg}}
\newcommand{\RS}{\mathfrak{R}}
\newcommand{\mdp}{\mathrm{mod~periods~}}
\newcommand{\z}	{{\boldsymbol z}}
\newcommand{\w}	{{\boldsymbol w}}
\newcommand{\tr}	{{\boldsymbol t}}
\newcommand{\e}	{{\boldsymbol e}}
\newcommand{\ualpha}{{\boldsymbol\alpha}}
\newcommand{\ubeta}{{\boldsymbol\beta}}
\newcommand{\ugamma}{{\boldsymbol\gamma}}

\newcommand{\rhy}   {\textnormal{RHP}-${\boldsymbol Y}$}
\newcommand{\rhs}   {\textnormal{RHP}-${\boldsymbol S}$}
\newcommand{\rhz}   {\textnormal{RHP}-${\boldsymbol Z}$}
\newcommand{\rhn}   {\textnormal{RHP}-${\boldsymbol N}$}
\newcommand{\rhp}   {\textnormal{RHP}-$\pmb{P}$}
\newcommand{\rhr}   {\textnormal{RHP}-${\boldsymbol R}$}
\newcommand{\rhpsi}   {\textnormal{RHP}-${\boldsymbol \Psi}$}
\newcommand{\rhphi}   {\textnormal{RHP}-${\boldsymbol \Phi}$}
\newcommand{\rhup}   {\textnormal{RHP}-${\boldsymbol \Upsilon}$}

\newtheorem{theorem}{Theorem}
\newtheorem{proposition}[theorem]{Proposition}

\newtheorem{condition}[theorem]{Condition}

\begin{document}

\title[Hermite-Pad\'e approximants for a pair of functions]{Hermite-Pad\'e approximants for a pair of Cauchy transforms with overlapping symmetric supports}

\author[A.I. Aptekarev]{Alexander I. Aptekarev}
\address{Keldysh Institute of Applied Mathematics, Russian Academy of Science, Moscow, Russian Federation}
\email{aptekaa@keldysh.ru}

\thanks{The research of the first author was carried out with support from a grant of the Russian Science Foundation (project RScF-14-21-00025).}

\author[W. Van Assche]{Walter Van Assche}
\address{Department of Mathematics, KU Leuven, Celestijnenlaan 200B box 2400, BE-3001 Leuven, Belgium}
\email{walter@wis.kuleuven.be}

\thanks{The second author was supported by KU Leuven research grant OT/12/073, FWO research project G.0934.13 and the Belgian Interuniversity Attraction Poles Programme P7/18.}

\author[M. Yattselev]{Maxim L. Yattselev}
\address{Department of Mathematical Sciences, Indiana University-Purdue University Indianapolis, 402~North Blackford Street, Indianapolis, IN 46202, USA}
\email{maxyatts@math.iupui.edu}

\begin{abstract}
Hermite-Pad\'e approximants of type II are vectors of rational functions with common denominator that interpolate a given vector of power series at infinity with maximal order. We are interested in the situation when the approximated vector is given by a pair of Cauchy transforms of smooth complex measures supported on the real line. The convergence properties of the approximants are rather well understood when the supports consist of two disjoint intervals (Angelesco systems) or two intervals that coincide under the condition that the ratio of the measures is a restriction of the Cauchy transform of a third measure (Nikishin systems). In this work we consider the case where the supports form two overlapping intervals (in a symmetric way) and the ratio of the measures extends to a holomorphic function in a region that depends on the size of the overlap. We derive Szeg\H{o}-type formulae for the asymptotics of the approximants, identify the convergence and divergence domains (the divergence domains appear for Angelesco systems but are not present for Nikishin systems), and show the presence of overinterpolation (a feature peculiar for Nikishin systems but not for Angelesco systems). Our analysis is based on a Riemann-Hilbert problem for multiple orthogonal polynomials (the common denominator).
\end{abstract}



\maketitle

\setcounter{tocdepth}{3}
\tableofcontents

\section{Introduction}

For $p\in\N$, let $\vec f = \big(f_1,\ldots,f_p\big)$ be a vector of holomorphic germs at infinity. Given a multi-index $\vec n=(n_1,\ldots,n_p)\in\N^p$, the \emph{type II Hermite-Pad\'e approximant} to $\vec f$ corresponding to $\vec n$,
\begin{equation}
\label{pi_vecn}
\vec \pi_{\vec n} = \left(\pi_{\vec n}^{(1)},\ldots,\pi_{\vec n}^{(p)}\right), \quad \pi_{\vec n}^{(i)} := \frac{P_{\vec n}^{(i)}}{Q_{\vec n}},
\end{equation}
is a vector of rational functions with common denominator $Q_{\vec n}$ satisfying
\begin{equation}
\label{R_vecn}
\left\{
\begin{array}{l}
\deg(Q_{\vec n}) \leq |\:\vec n\:|=n_1+\cdots+n_p \medskip \\
R_{\vec n}^{(i)}(z):=\left(Q_{\vec n}f_i - P_{\vec n}^{(i)}\right)(z) = \mathcal{O}\big(z^{-n_i-1}\big) \quad \mbox{as} \quad z\to\infty
\end{array}
\right.
\end{equation}
for each $i\in\{1,\ldots,p\}$. We shall not deal with type I Hermite-Pad\'e approximants in this work and therefore henceforth we will drop the ``type II'' modifier. Such approximants were introduced by Hermite \cite{Herm73} for the vector of exponentials $(1,e^z,\ldots,e^{(p-1)z})$, with the interpolation taking place at the origin rather than at infinity, as a tool in proving the transcendence of $e$. Later, his student Pad\'e systematically studied the scalar case $p=1$ \cite{Pade92} and such approximants are now called \emph{Pad\'e approximants}.

From our perspective interpolating at infinity is more convenient than interpolating at the origin in the following sense. Any holomorphic function can be written as a Cauchy integral of its boundary values
on any curve encircling a domain of analyticity. For a holomorphic function on a domain we will use the terminology \emph{trace} to mean the boundary values of the function on the boundary of the domain.
When the function is holomorphic at infinity, such an integral representation can in some important cases be deformed analytically into an integral over a ``one dimensional'' set. A particular fruitful example of this principle is when we deal with Markov functions.
Let $\mu_i$ be a positive measure and take
\begin{equation}
\label{markov}
f_i(z) = \int\frac{\mathrm{d}\mu_i(x)}{x-z},  \qquad \supp(\mu_i)\subseteq[a_i,b_i]\subset\R.
\end{equation}
The vector $\vec f$ is called an \emph{Angelesco system} if $[a_j,b_j]\cap[a_i,b_i]=\varnothing$ for $j\neq i$ (such systems were initially considered by Angelesco \cite{Ang19} and later rediscovered by Nikishin \cite{Nik79}). Angelesco has shown that $Q_{\vec n}$ has exactly $n_i$ zeros on $[a_i,b_i]$. This means that an Angelesco system is an example of a \emph{perfect} system, i.e., a system for which every multi-index is \emph{normal}, that is, $\deg(Q_{\vec n})=|\:\vec n\:|$. As far as the asymptotic behavior of $\pi^{(j)}_{\vec n}$ is concerned, convergence properties are not as straightforward as one could hope. Given an Angelesco system and a sequence of multi-indices such that $n_i/|\:\vec n\:|\to c_i>0$, $c_1+\cdots+c_p=1$, Gonchar and Rakhmanov \cite{GRakh81} have shown that for each $j$ the complement of $\cup_{i=1}^p[a_i,b_i]$ is separated by a system of analytic arcs into two domains, say $D_j^+$, containing the point at infinity, and $D_j^-$, possible empty, such that $\pi^{(j)}_{\vec n}$ converges to $f_j$ in $D_j^+$ and diverges to infinity in $D_j^-$. Moreover, the polynomials $Q_{\vec n}$ can have an asymptotic zero distribution (in the sense of weak$^*$ convergence) on a strict subset of $\cup_{i=1}^p\supp(\mu_i)$ (\emph{pushing effect}). The pushing effect always implies existence of a divergence region but the reverse implication is not true. More detailed (strong) asymptotics for Hermite-Pad\'e approximants to Angelesco systems when $p=2$ and the weights $\mathrm{d}\mu_i/\mathrm{d}x$ satisfy the so-called \emph{Szeg\H{o} condition}, was obtained by the first author in \cite{Ap88}.

Another class of Markov functions for which positive convergence results were obtained is now known as \emph{Nikishin systems} \cite{Nik80}. The functions $f_i$ in \eqref{markov} form a Nikishin system if they all are supported on the same interval $[a,b]$ and the Radon-Nikodym derivatives $\mathrm{d}\mu_j/\mathrm{d}\mu_1$, $j\in\{2,\ldots,p\}$, form a Nikishin system of order $p-1$ with respect to some interval $[c,d]$ such that $[c,d]\cap[a,b]=\varnothing$. Nikishin himself \cite{Nik80} has shown that such systems are perfect when $p=2$ and the Hermite-Pad\'e approximants converge uniformly outside of the interval $[a,b]$ in this case. This puts Nikishin systems more in line with the Pad\'e case $p=1$ (Markov theorem \cite{Mar95}) as neither the pushing effect nor the possibility of non-empty divergence regions appears for them. However, Nikishin systems do possess one new phenomenon, namely overinterpolation. It turns out that $R_{\vec n}^{(2)}$ has zeros on $[c,d]$ that are dense on this interval. It took 30 years to prove that Nikishin systems are perfect for any $p$ \cite{FPrLL11}. In \cite{FPrLL11}, Fidalgo Prieto and L\'opez Lagomasino also proved uniform convergence for multi-indices close to the diagonal. Strong asymptotics in the case of diagonal multi-indices and Szeg\H{o} weights was derived by the first author in \cite{Ap99}.

It is interesting to observe that the first result on strong asymptotics of Hermite-Pad\'e approximants was obtained by Kalaygin \cite{Kal79} for the case of two touching symmetric intervals (the limiting case of an Angelesco system).

As often happens in mathematics, the treatment of Angelesco and Nikishin systems can be unified under the umbrella of \emph{generalized Nikishin systems} (GN-systems) as introduced in \cite{GRakhS97}, where Gonchar, Rakhmanov, and Sorokin defined  a system of Markov functions with the help of a rooted tree graph and considered the question of uniform convergence of Hermite-Pad\'e approximants to such a system. In such a set-up, an Angelesco system corresponds to a tree where the root is connected by $p$ edges to $p$ leaves and a Nikishin system corresponds to a tree in which every node except for the final leaf has exactly one child. Strong asymptotics of Hermite-Pad\'e approximants to GN-systems of Markov functions generated by more general (than rooted tree) graphs (admitting cycles) was derived by the first author and Lysov in \cite{ApLy10}. An example of GN-system from \cite{ApLy10} is a pair of two Markov functions in \eqref{markov} where the support of one of them is strictly included in the support of the other, i.e., $ \supp(\mu_2)=[a_2,b_2] \subset \supp(\mu_1)=[a_1,b_1] $ and the Radon-Nikodym derivative $\mathrm{d}\mu_2/\mathrm{d}\mu_1$ along $[ a_1,b_1] $ is a Markov function with support $[ a_3,b_3] $ where $[a_3,b_3]\cap[a_1,b_1]=\varnothing$. Weak asymptotics of Hermite-Pad\'e approximants to this example of Markov functions was derived by Rakhmanov in \cite{Rakh11}.

We emphasize that in all the results we listed above the \textit{geometry} of the problem is \textit{real}, i.e., the supports of the limiting distributions of the poles of the approximants and the overinterpolation points belong to $\R$.

In this paper we consider Hermite-Pad\'e approximants to a pair of Cauchy transforms (Markov functions) of generally speaking complex measures with overlapping supports and aim at strong asymptotics. The set up does not fall into the framework of GN-systems, even when the measures are positive, as their supports overlap (for GN-systems the supports are either disjoint or one coincides with or contains the other). It turns out that both phenomena, the pushing effect and overinterpolation, appear in this case. Hermite-Pad\'e approximants to a pair of Markov functions with overlapping supports were first considered by Stahl with the goal of proving weak asymptotics \cite{ St87, uSt2, AptSt92}. He had the important insight that the geometry of this problem is complex, i.e., the overinterpolation points are distributed on analytic arcs in the complex plane (later, a similar effect was observed in \cite{AptBlK05}). This discovery was very unusual at the time because the input geometry (i.e., the supports of the measures generating the Markov functions) is completely real. Unfortunately, Stahl's results have never been published. This work was strongly motivated by the desire to provide a detailed proof of his findings (in an even more delicate setting of strong asymptotics).

Opting here for complex measures is natural from the point of view of complex analysis. However, many techniques, like those in \cite{Ap88,Ap99}, do not apply as they use positivity in an essential way. An approach that does not rely on positivity was outlined by Nuttall in his seminal paper \cite{Nut84}. There Nuttall conjectured that the main term of the asymptotics of Hermite-Pad\'e approximants is a function solving of a certain \emph{explicit} boundary value problem on some \emph{unknown} Riemann surface. He identified this surface only in a handful of special cases. Elaborating on Nuttall's approach, the first two authors and Kuijlaars \cite{ApKVA08} pinpointed the algebraic equation which defines the appropriate Riemann surface in the case of two Cauchy transforms of complex measures supported on two arcs joining pairs of  branch points in the complex plane (the simplest example is a complexified Angelesco system) and derived formulae of strong asymptotics in the case when the Riemann surface has genus zero. Below, we build upon the ideas developed in \cite{ApKVA08} and extend the results of \cite{ApKVA08} to the cases when the appropriated Riemann surface has positive genus (elliptic and ultra- elliptic case).

In Section~\ref{sec:RS} we identify the Riemann surface in Nuttall's program  by an algebraic equation, discuss its realization as a ramified cover of $\overline\C$, and  construct a certain function on this surface whose level lines will geometrically describe convergence and divergence domains of the approximants. In Section~\ref{sec:NS}, we construct the Nuttall-Szeg\H{o} functions that will provide the leading term of the asymptotics of the Hermite-Pad\'e approximants. Finally, in Section~\ref{sec:HP} we state the main result of this work. The remaining part of the paper is devoted to the proofs of all the stated results.

\section{Riemann Surface}
\label{sec:RS}

Let $a\in(0,1)$ be given. Our goal is to investigate Hermite-Pad\'e approximants to a pair of Markov-type functions generated by measures with supports $[-1, a]$ and $[-a, 1]$. To this end we consider the algebraic equation
\begin{equation}
\label{h}
A(z)h^3 - 3B_2(z)h - 2B_1(z) = 0,
\end{equation}
where the polynomials $A(z)$, $B_2(z)$, and $B_1(z)$ are defined by
\[
\left\{
\begin{array}{lcl}
A(z) &:=& (z^2-1)(z^2-a^2), \medskip \\
B_2(z) &:=& z^2-p^2, \medskip \\
B_1(z) &:=& z,
\end{array}
\right.
\]
for some parameter $p>0$. Denote by $h_k$, $k\in\{0,1,2\}$, the three distinct branches of the algebraic function $h$ determined by \eqref{h}. Naturally, these branches satisfy
\begin{equation}
\label{branches}
\left\{
\begin{array}{l}
h_0+h_1+h_2\equiv0, \medskip \\
h_0h_1+h_0h_2+h_1h_2 = -3B_2/A, \medskip \\
h_0h_1h_2  = 2B_1/A.
\end{array}
\right.
\end{equation}
Hence, it necessarily holds that $h_k(\infty)=0$ and therefore $h_k(z)=c_k/z+\cdots$ as $z\to\infty$. It readily follows from the above equations that the constants $c_k$ are the solutions of $0=c^3-3c-2=(c+1)^2(c-2)$. Thus, we put
\begin{equation}
\label{hatinfinity}
\left\{
\begin{array}{lcr}
h_0(z) &=& \displaystyle \frac2z + \cdots \smallskip\\
h_i(z) &=& \displaystyle -\frac1z + \cdots
\end{array}
\right. \quad \mbox{as} \quad z\to\infty
\end{equation}
for $i\in\{1,2\}$.  It can easily be checked that all three solutions of \eqref{h} are real for positive large $x$. Hence, we can label the branches so that for all $x>0$ large enough
\begin{equation}
\label{label}
h_0(x) > h_1(x) > h_2(x).
\end{equation}

Denote by $\RS$ the Riemann surface of $h$. It is a three-sheeted ramified cover of $\overline\C$. We shall denote by $\RS^{(k)}$, $k\in\{0,1,2\}$, the sheet on which $h$ coincides with $h_k$ (a particular realization of this surface is specified in Theorem~\ref{thm:RS} below). We denote by $\z$ a generic point on $\RS$ with \emph{natural projection} $\pi(\z)=z\in\overline\C$. If we want to specify the sheet, we write $z^{(k)}$ for $z^{(k)}\in\RS^{(k)}$ so that $\pi(z^{(k)})=z$. Thus, $h_k(z)=h(z^{(k)})$. Consider the differential
\begin{equation}
\label{GreenDiff}
\mathrm{d}\mathcal N(\z) := h(\z)\mathrm{d}z,
\end{equation}
where $\z$ is a generic point on $\RS$. The choice of the parameter $p$ is driven by the following condition:
\begin{equation}
\label{condition}
N(\z) := \re\left(\int^\z\mathrm{d}\mathcal N\right) \quad \text{is a well defined (\emph{single-valued}) harmonic function on} \; \RS.
\end{equation}
When it exists, $N(\z)$ is defined up to an additive constant. If we denote by $N_k$ the restriction of $N$ to $\RS^{(k)}$, then it is easy to see that $N_1+N_2+N_3$ is a well defined harmonic function in $\overline\C$ and therefore it is constant. Thus, we normalize $N$ so that
\begin{equation}
\label{2.6_1}
N_1(z)+N_2(z)+N_3(z) \equiv 0, \quad z\in\overline\C.
\end{equation}
Notice also that \eqref{condition} is equivalent to requiring that all the periods of the differential $\mathrm{d}\mathcal N$ are purely imaginary.

\begin{figure}[!ht]
\centering
\subfigure[Case I]{\includegraphics[scale=.5]{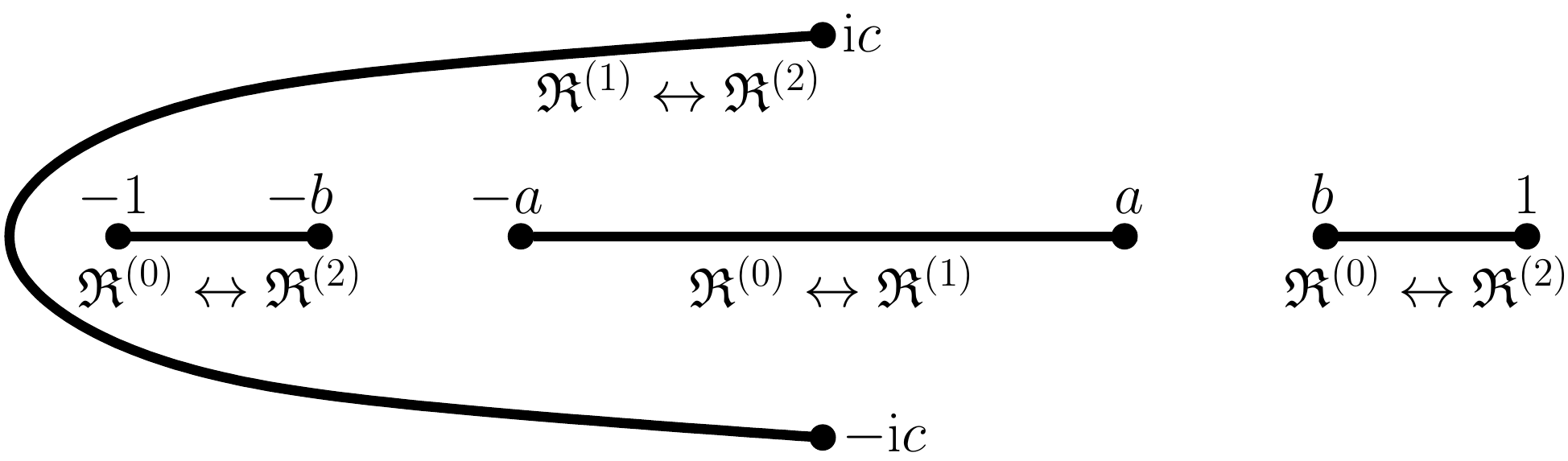}}
\subfigure[Case II]{\includegraphics[scale=.5]{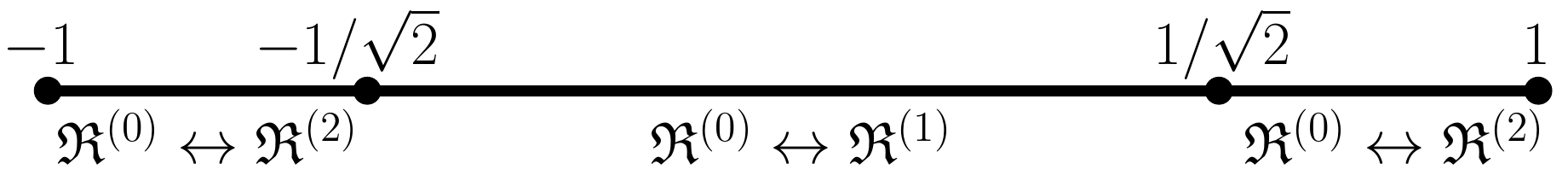}}
\subfigure[Case III]{\includegraphics[scale=.5]{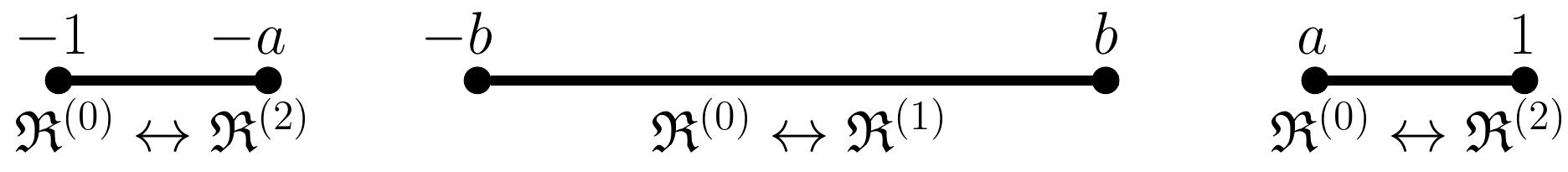}}
\caption{\small Ramification points of $\RS$ and the cuts (black curves) along which the sheets $\RS^{(0)},\RS^{(1)},\RS^{(2)}$ are glued to each other.}
\label{RS}
\end{figure}

\begin{theorem}
\label{thm:RS}
Consider the algebraic equation \eqref{h} with $a\in(0,1)$.
\begin{itemize}
\item[(I)] If $a\in\big(0,1/\sqrt 2\big)$, then there exists $p\in\big(a,\sqrt{(1+a^2)/3}\big)$ such that condition \eqref{condition} is fulfilled. In this case $\RS$ has eight ramification points whose projections are $\{\pm1,\pm a\}$ and $\{\pm b,\pm\mathrm{i}c\}$ for some uniquely determined $b\in(a,p)$ and $c>0$. Moreover, the surface can be realized as on Figure~\ref{RS}(a);
\item[(II)] If  $a=1/\sqrt 2$, then condition \eqref{condition} is fulfilled for $p=1/\sqrt 2$. In this case $\RS$ has four ramification points whose projections are $\big\{\pm1,\pm1/\sqrt2\big\}$ and it can be realized as on Figure~\ref{RS}(b);
\item[(III)] If $a\in\big(1/\sqrt 2,1\big)$, then condition \eqref{condition} is fulfilled for
    for $p=\sqrt{(1+a^2)/3}$. In this case $\RS$ has six ramification points whose projections are $\{\pm1,\pm a\}$ and $\{\pm b\}$ for some uniquely determined $b\in(p,a)$. Moreover, the surface can be realized as on Figure~\ref{RS}(c).
\end{itemize}
\end{theorem}

The points $b$ and $c$ in Case I and the point $b$ in Case III can be explicitly computed as they are solutions of a certain explicit quadratic or linear (in $z^2$) equation whose parameters depend on $a$ and $p$.

It follows from Theorem~\ref{thm:RS} that $\RS$ has genus $g=2$ when $a\in\big(0,1/\sqrt2\big)$ and genus $g=1$ otherwise. Moreover, in Cases I and III, all the ramification points have order 2 while in Case II the points $\pm1$ have order 2 while $\pm1/\sqrt 2$ have order 3.

In the light of Theorem~\ref{thm:RS}, it will be convenient to fix the following notation:
\begin{equation}
\label{chains}
\left\{
\begin{array}{lll}
\Delta_0 &:=& \pi\big(\text{cycle that separates } \RS^{(1)} \text{ and } \RS^{(2)} \big), \medskip \\
\Delta_1 &:=& \pi\big(\text{cycle that separates } \RS^{(0)} \text{ and } \RS^{(1)} \big), \medskip \\
\Delta_2 &:=& \Delta_{21} \cup \Delta_{22}, \medskip \\
\Delta_{21} &:=& \pi\big(\text{the left cycle of the chain that separates } \RS^{(0)} \text{ and } \RS^{(2)} \big), \medskip \\
\Delta_{22} &:=& \pi\big(\text{the right cycle of the chain that separates } \RS^{(0)} \text{ and } \RS^{(2)} \big).
\end{array}
\right.
\end{equation}
Clearly, $\Delta_0$ is defined only in Case I.

\begin{figure}[!ht]
\centering
\subfigure[Case I]{\includegraphics[scale=.5]{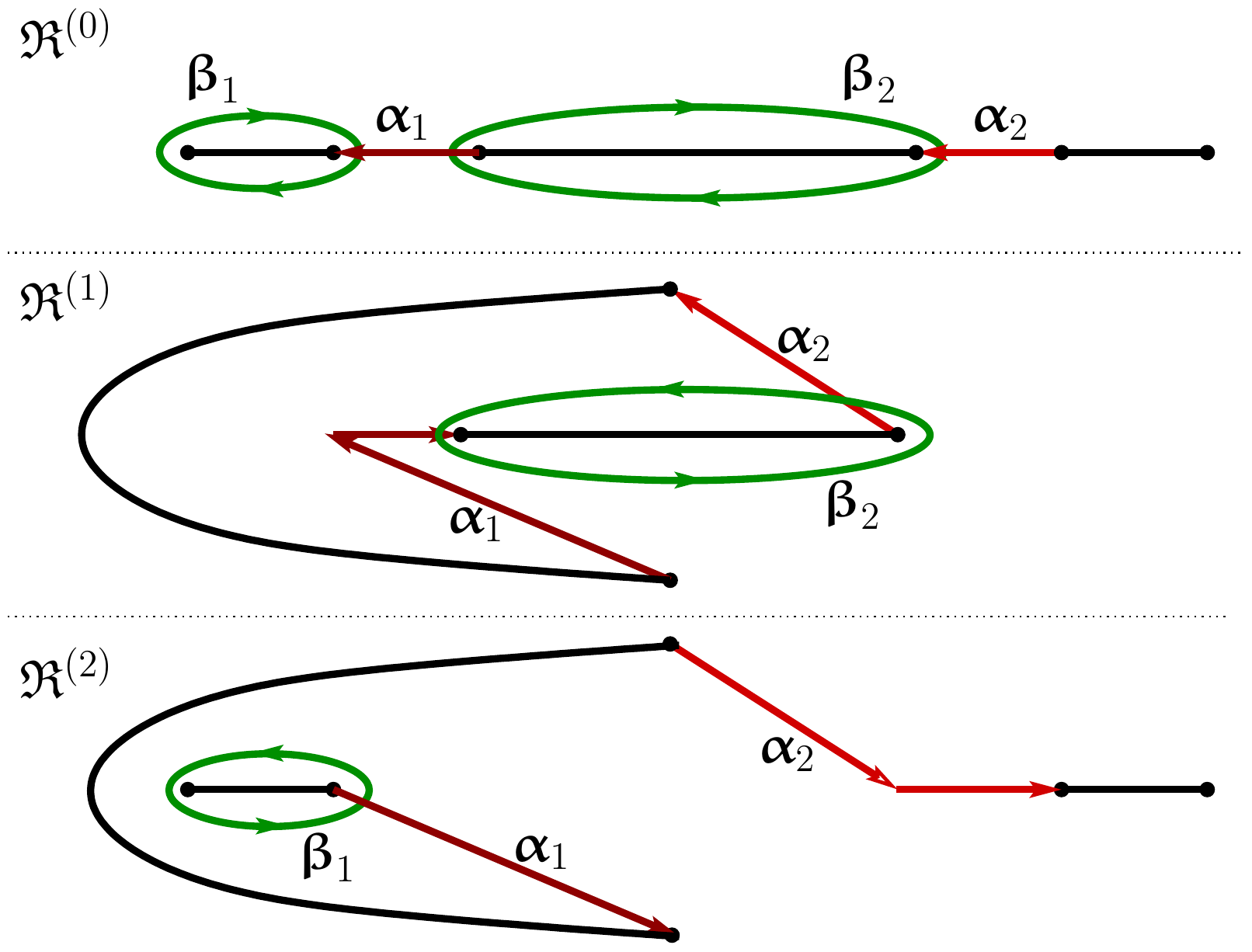}}
\subfigure[Cases II and III (there is no gap between the cuts on $\RS^{(0)}$ in Case II, but this does not affect the choice of the cycles)]{\includegraphics[scale=.5]{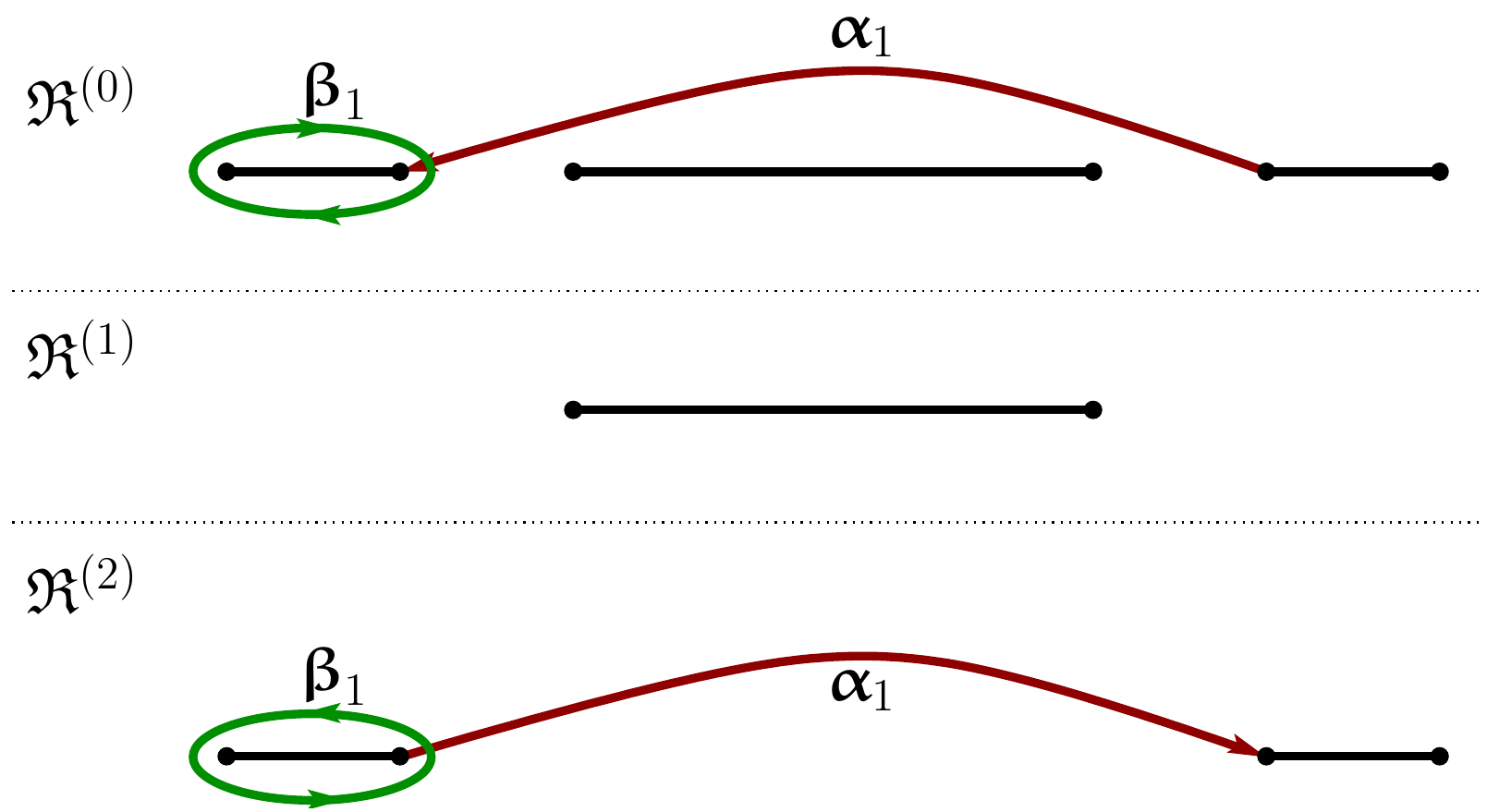}}
\caption{\small Homology basis for $\RS$.}
\label{fig:HB}
\end{figure}

Let $\mathrm{d}\mathcal N$ be defined by \eqref{GreenDiff}.
The function
\begin{equation}
\label{Phi}
\Phi(\z) := \exp\left\{\int^{\z}\mathrm{d}\mathcal N\right\}
\end{equation}
is analytic, except for a double pole at $\infty^{(0)}$, and multiplicatively multi-valued on $\RS$. Moreover, it is single-valued in $\RS_{\ualpha,\ubeta}:=\RS\setminus\bigcup_{i=1}^g(\ualpha_i\cup\ubeta_i)$, where $\{\ualpha_i,\ubeta_i\}_{i=1}^g$ is a homology basis on $\RS$.   Later on, see Figure~\ref{fig:HB}, we shall specify the basis in more detail, but right now it is sufficient to assume that each cycle  $\ugamma\in\{\ualpha_i,\ubeta_i\}_{i=1}^g$ possesses a \emph{projective involution}: $|\ugamma\cap\pi^{-1}(z)|=2$ for any $z\in\pi(\ugamma)$ which is not a branch point of $\RS$ (the involution is then defined by mapping a point on $\ugamma$ to the other one with the same projection). We normalize $\Phi$ so that
\begin{equation}
\label{normal}
\Phi^{(0)}\Phi^{(1)}\Phi^{(2)} \equiv 1 \quad \text{in} \quad \overline\C,
\end{equation}
where $\Phi^{(k)}$ is the pullback to $\C$ of the restriction of $\Phi$ to $\RS^{(k)}$. Let us show that such a normalization is indeed possible. Since $\RS_{\ualpha,\ubeta}$ is simply connected and $\mathrm{d}\mathcal N$ has integer residues, the restriction of $\Phi$ to $\RS_{\ualpha,\ubeta}$ is single-valued. It satisfies
\begin{equation}
\label{Phi-jumps}
\Phi^+ = \Phi^-\left\{
\begin{array}{ll}
\exp\big\{2\pi \mathrm{i}\omega_i\big\} & \mbox{on} \quad \ualpha_i, \medskip \\
\exp\big\{2\pi \mathrm{i}\tau_i\big\} & \mbox{on} \quad \ubeta_i,
\end{array}
\right.
\end{equation}
$1\leq i\leq g$, where the constants $\omega_i$ and $\tau_i$ are real (this is guaranteed by \eqref{condition}) and given by
\begin{equation}
\label{GreenPeriods}
\omega_i :=  -\frac{1}{2\pi \mathrm{i}}\oint_{\ubeta_i}\mathrm{d}\mathcal N, \qquad \tau_i:=\frac{1}{2\pi \mathrm{i}}\oint_{\ualpha_i}\mathrm{d}\mathcal N.
\end{equation}
At times, it will be convenient to use vector notation $\vec\tau=(\tau_1,\tau_2)^\mathsf{T}$ and $\vec\omega=(\omega_1,\omega_2)^\mathsf{T}$ if $g=2$ and $\vec\tau=(\tau_1)$ and $\vec\omega=(\omega_1)$ if $g=1$.  Since $\Phi$ has a double pole at $\infty^{(0)}$ and simple  zeros at $\infty^{(1)},\infty^{(2)}$, we can write
\begin{equation}
\label{Phi_jatInfty}
\left\{
\begin{array}{lcl}
\Phi^{(0)}(z) &=& \displaystyle C_0z^2 + \cdots \medskip\\
\Phi^{(i)}(z) &=& \displaystyle C_iz^{-1} + \cdots
\end{array}
\right. \quad \mbox{as} \quad z\to\infty.
\end{equation}
It follows from \eqref{2.6_1} that $\log|\Phi^{(0)}\Phi^{(1)}\Phi^{(2)}|\equiv0$ on $\overline\C$. Hence, if we choose the normalization $C_0C_1C_2=1$, then \eqref{normal} is fulfilled due to our choice of the homology basis: since the cycles possess projective involutions, exactly two pullbacks of $\Phi$ have jumps at each point belonging to the projection of a homology cycle, moreover, the jumps are reciprocal as the projected parts of the cycle coincide as sets but have opposite orientations.

Due to the symmetries of the Riemann surface, the vectors $\vec\omega$ and $\vec\tau$ take special forms. To be more specific, let us fix a homology basis. We choose $\pi(\ubeta_1)=\Delta_{21}$ and, in Case I, $\pi(\ubeta_2)=\Delta_1$, see \eqref{chains}, while the $\ualpha$-cycles are as on Figure~\ref{fig:HB}(a) in Case I and the $\ualpha_1$-cycle should be chosen as on Figure~\ref{fig:HB}(b) in Cases II and III.

\begin{proposition}
\label{prop:periods}
Let $\{\ualpha_i,\ubeta_i\}_{i=1}^g$ be the homology basis which we just fixed. In Case I, one has
\begin{equation}
\label{not-rational1}
\vec\omega = \big(\omega,2(1-\omega)\big)^\mathsf{T} \quad \text{and} \quad \vec\tau = \big(\tau,-\tau\big)^\mathsf{T}
\end{equation}
for some real constants $\tau=\tau(a)$ and $\omega=\omega(a)\in(1/2,1)$. In Cases II and III, one has
\begin{equation}
\label{not-rational2}
\vec\omega = \vec\tau = (1/2).
\end{equation}
\end{proposition}

Recall that $|\Phi(\z)|=\exp\{N(\z)\}$ is a single-valued function on $\RS$, see \eqref{2.6_1}. The asymptotics of the Hermite-Pad\'e approximants will depend on the relative sizes of the different branches of $|\Phi|$ (in other words, on the size of the branches of $N$). To this end, define $\Omega_{ijk}$ to be open subset of $\overline\C$ such that
\[
\Omega_{ijk} := \left\{z:\big|\Phi^{(i)}(z)\big|>\big|\Phi^{(j)}(z)\big|>\big|\Phi^{(k)}(z)\big|\right\}.
\]
We also define the closed set $\Gamma:=\Gamma_{01}\cup\Gamma_{02}\cup\Gamma_{12}$ by
\[
\Gamma_{ij} := \left\{z:\big|\Phi^{(i)}(z)\big|=\big|\Phi^{(j)}(z)\big|,~i\neq j\right\}.
\]
Clearly, $\overline\C\setminus\Gamma=\bigcup_{i\neq j\neq k\neq i}\Omega_{ijk}$. Then the following theorem holds.

\begin{figure}[!ht]
\centering
\subfigure[Case I]{\includegraphics[scale=.5]{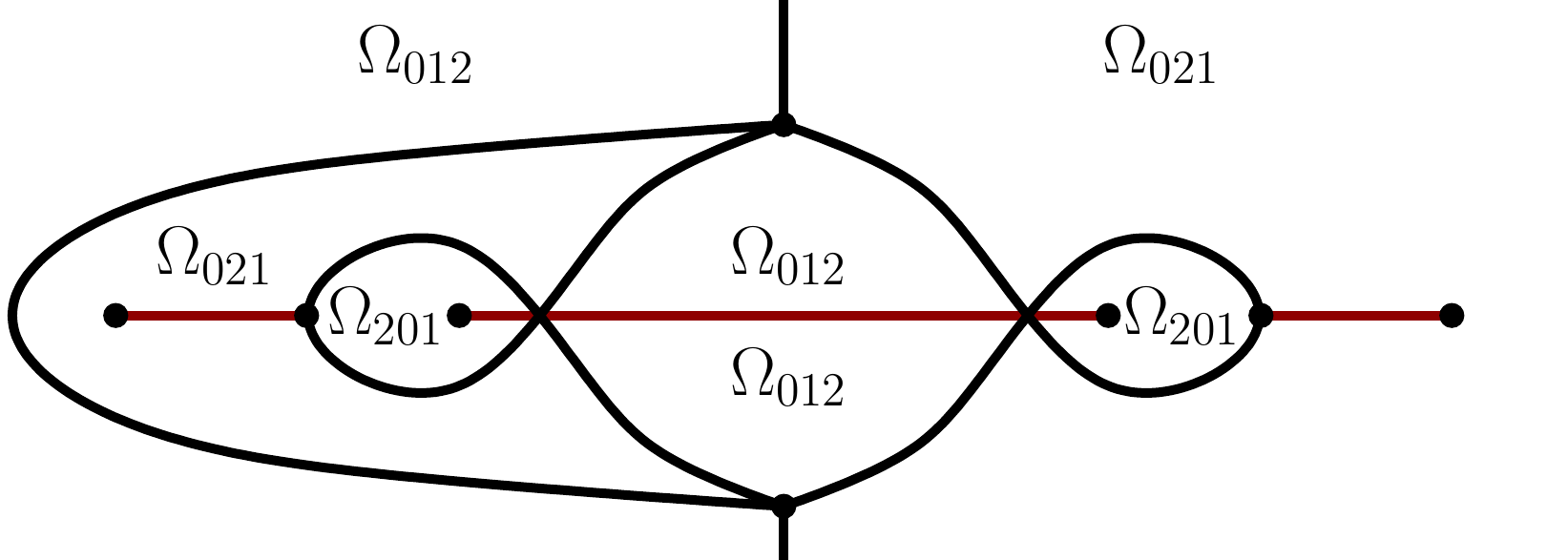}}
\subfigure[Case II]{\includegraphics[scale=.5]{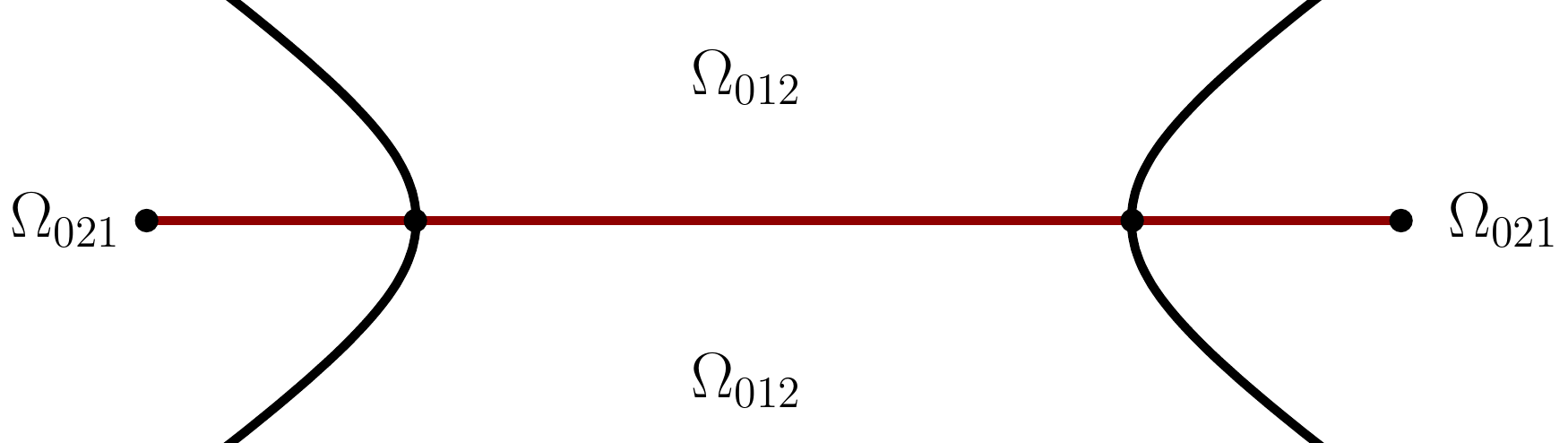}}\\
\subfigure[Case IIIa]{ \includegraphics[scale=.5]{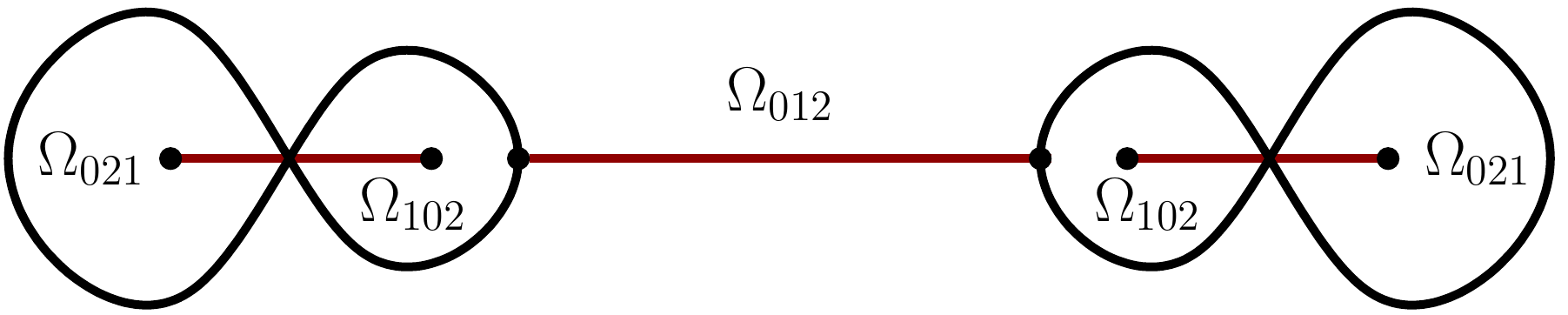}} \\
\subfigure[Case IIIb]{\includegraphics[scale=.5]{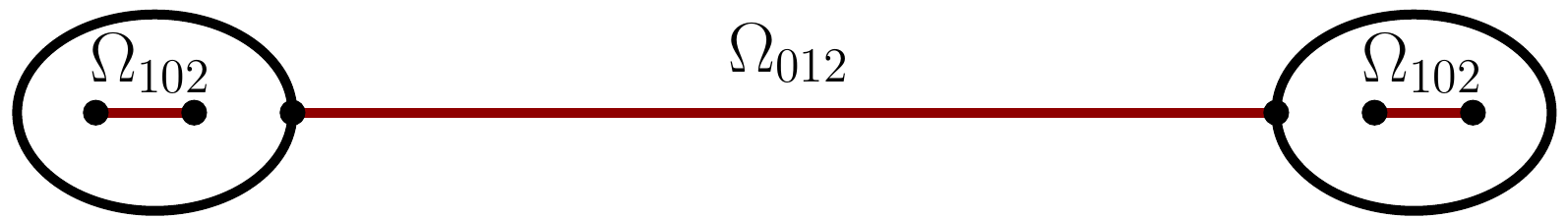}}
\caption{\small Domains $\Omega_{ijk}$.}
\label{Oms}
\end{figure}

\begin{theorem}
\label{thm:Omegas}
In Case I, the regions $\Omega_{ijk}$ are distributed as on Figure~\ref{Oms}(a); in Case II, the domains are distributed as on Figure~\ref{Oms}(b); in Case III, the domains are distributed as on Figure~\ref{Oms}(c,d).
\end{theorem}

We prove Theorems~\ref{thm:RS} and~\ref{thm:Omegas} as well as Proposition~\ref{prop:periods} in Section~\ref{sec:Geom}.

\section{Nuttall-Szeg\H{o} Functions}
\label{sec:NS}

 To define the Nuttall-Szeg\H{o} functions, we first need to formulate a certain \emph{Jacobi inversion problem}. To this end, denote by
\[
\vec\Omega(\z)=\big(\Omega_1(\z),\ldots,\Omega_g(\z)\big)^\mathsf{T}
\]
the vector of  $g$ normalized  holomorphic integrals on $\RS$. That is, $\vec\Omega(\z)$ is a vector of analytic and additively multi-valued functions on $\RS$ that are single-valued on $\RS_{\ualpha,\ubeta}$ and satisfy
\begin{equation}
\label{Omega-jump}
\vec\Omega^{+}-\vec\Omega^- = \left\{
\begin{array}{ll}
-\mathbf B\vec e_i & \mbox{on} \quad \ualpha_i,\\
\vec e_i & \mbox{on} \quad \ubeta_i,
\end{array}
\right. \quad \mathbf{B} := \left[\oint_{\ubeta_j}\mathrm{d}\Omega_i\right]_{j,i=1}^g,
\end{equation}
where $\vec e_i$ is the $i$-th vector of the standard basis in $\R^g$. It is known that $\mathbf{B}$ is a symmetric matrix with positive definite imaginary part when $g=2$ and is a complex number with $\mathrm{Im}(\mathbf{B})>0$ when $g=1$. The Jacobi inversion problem we need consists of finding an integral divisor\footnote{Recall that an \emph{integral divisor} of order $d$ is a formal expression $\mathcal{D}=\sum_in_i\z_i$, where $n_i\in\N$ and $\sum_in_i=d$. A \emph{principal divisor} is an expression $\sum_{i=1}^kn_i\tr_i-\sum_{i=1}^jm_i\w_i$ such that there exists a rational function on $\RS$ with a zero of multiplicity $n_i$ at $\tr_i$ for each $i\in\{1,\ldots,k\}$, a pole of order $m_i$ at $\w_i$ for each $i\in\{1,\ldots,j\}$, and is otherwise non-vanishing and finite (in particular, one has that $\sum_{i=1}^kn_i=\sum_{i=1}^jm_i$).} of order $g$, say $\mathcal{D}_n$, such that
\begin{equation}
\label{one-jip}
\vec\Omega\big(\mathcal{D}_n\big) \equiv \vec\Omega\big(g\infty^{(2)}\big) + \vec c_\rho + n\big(\vec\omega+\mathbf{B}\vec\tau\big) \quad \big(\mdp \mathrm{d}\vec\Omega\big),
\end{equation}
where $\vec c_\rho$ is a vector of constants that depends only on the approximated functions and is defined further below in \eqref{Srho}, $\vec\tau$ and $\vec\omega$ are the vectors defined after \eqref{GreenPeriods}, and the equivalence  of two vectors $\vec c,\vec e\in\C^g$ is defined by
\begin{equation}
\label{equivalence}
\vec c\equiv \vec e \quad \left(\mdp \mathrm{d}\vec\Omega\right) \quad \Leftrightarrow \quad \vec c - \vec e=\vec j+\mathbf{B}\vec m, \qquad \vec j,\vec m\in\Z^g.
\end{equation}
It is known that the Jacobi inversion problem has a unique solution when $g=1$. Hence, $\mathcal{D}_n$ is well defined for all $n$ in this case. Moreover, the following proposition holds.
\begin{proposition}
\label{prop:genus1} Assuming $g=1$, let $\mathcal{D}_n$ be the unique solution of \eqref{one-jip}. Then $\mathcal{D}_{2m} = \mathcal{D}_0$ and $\mathcal{D}_{2m+1} = \mathcal{D}_1$ for all $m\geq0$. Moreover,  either $\mathcal{D}_0\neq\infty^{(0)}$ or $\mathcal{D}_1\neq\infty^{(0)}$ and therefore $\N_*:=\big\{n\in\N:~\mathcal{D}_n\neq\infty^{(0)}\big\}$ is infinite.
\end{proposition}

The unique solvability of a Jacobi inversion problem is no longer guaranteed when $g=2$. However, it is known that solutions are either unique or given by any \emph{special divisor}\footnote{On genus 2 surfaces, a \emph{special divisor} is an integral divisor of order 2 such that there exists a rational function on $\RS$ with simple poles at the elements of the divisor (a double pole if they coincide) and otherwise regular.}. Notice that integral divisors of order 2 can be considered as elements of $\RS^2/\Sigma_2$, where $\Sigma_2$ is the symmetric group of two elements. This is a compact topological space. Hence, we can talk about limit points of sequences of integral divisors.
\begin{proposition}
\label{prop:genus2}
 Let $g=2$, then \eqref{one-jip} has a unique solution for at least one of the indices $n-1,n$ for any $n\in\N$. Moreover, there always exists an infinite subsequence $\N_*$ such that no limit point of $\big\{\mathcal{D}_n\big\}_{n\in\N_*}$ is a special divisor or is of the form $\infty^{(0)}+\w$ for some $\w\in\RS$.
\end{proposition}

The following theorem is the main result of this section.

\begin{theorem}
\label{thm:N-S}
Let $\rho_1$ and $\rho_2$ be functions holomorphic and non-vanishing in a neighborhood of $[-1,1]$. In Case I, assume in addition that the ratio $\rho_2/\rho_1$ holomorphically extends to a non-vanishing function in a neighborhood of $\Delta_0$. Denote by $\boldsymbol\Delta$ the chain that separates $\RS^{(0)}$, $\RS^{(1)}$, and $\RS^{(2)}$ ($\pi(\boldsymbol\Delta)$ is a union of sets defined in \eqref{chains}). For each $n$ for which \eqref{one-jip} is uniquely solvable, there exists a function, say $\Psi_n$, meromorphic in $\RS\setminus\boldsymbol\Delta$ such that the zero/pole multi-set of $\Psi_n$ is given by
\begin{equation}
\label{Psi-Div}
\mathcal{D}_n + (n+1)\big(\infty^{(1)}+\infty^{(2)}\big)-2n\infty^{(0)}.
\end{equation}
Its traces are bounded except at the branch points of $\RS$ where $\Psi_n$ behaves like
\begin{equation}
\label{around-branching}
\big|\Psi_n(\z)|\sim |z-e|^{-1/4} \quad \text{as} \quad \z\to \e,
\end{equation}
for every branch point $\e$ of $\RS$ except when $e=\pm1/\sqrt2$ (Case II only) where the exponent $-1/4$ should be replaced by $-1/3$. Furthermore
\begin{equation}
\label{BVP}
\left\{
\begin{array}{lllll}
\left(\Psi_n^{(1)}\right)^\pm &=& \pm\left(\Psi_n^{(0)}\right)^\mp\rho_1 &\text{on} &  \Delta_1^\circ, \medskip \\
\left(\Psi_n^{(2)}\right)^\pm &=& \mp\left(\Psi_n^{(0)}\right)^\mp\rho_2 &\text{on} &  \Delta_{21}^\circ, \medskip \\
\left(\Psi_n^{(2)}\right)^\pm &=& \pm\left(\Psi_n^{(0)}\right)^\mp\rho_2 &\text{on} &  \Delta_{22}^\circ, \medskip \\
\left(\Psi_n^{(2)}\right)^\pm &=& \pm\left(\Psi_n^{(1)}\right)^\mp(\rho_2/\rho_1) & \text{on} & \Delta_0^\circ,
\end{array}
\right.
\end{equation}
where $\Delta^\circ$ is the interior of the arc $\Delta$. Moreover, if $\Psi$ is a function meromorphic in $\RS\setminus\boldsymbol\Delta$ satisfying \eqref{BVP}, \eqref{around-branching}, and \eqref{Psi-Div} with $\mathcal{D}_n$ replaced by any other integral divisor $\mathcal{D}$, then $\Psi$ is a constant multiple of $\Psi_n$.
\end{theorem}

It might seem, especially after looking at \eqref{one-jip}, that the construction of the Nuttall-Szeg\H{o} functions $\Psi_n$, particularly, the divisors $\mathcal{D}_n$, depends on the choice of the homology basis. However, this is precisely the point of the uniqueness part of Theorem~\ref{thm:N-S} that the functions $\Psi_n$ are independent of the intermediate steps of their construction.

The Nuttall-Szeg\H{o} functions $\Psi_n$ are certain multiples of the functions $\Phi^n$. In fact, their ratios form a normal family on each $\RS_\varepsilon$ obtained from $\RS$ by excising  circular neighborhoods of radius $\varepsilon$ around the branch points.

\begin{theorem}
\label{thm:N-Sasymp}
Let $\N_*\subseteq\N$ be a subsequence as in Propositions~\ref{prop:genus1} or~\ref{prop:genus2}. For each $\varepsilon>0$, there exists a constant $C_\varepsilon(\N_*)>1$ such that
\begin{equation}
\label{Psi-Est}
\left\{
\begin{array}{lll}
\left|\Psi_n\right| \leq C_\varepsilon(\N_*)\left|\Phi^n\right| & \text{in} & \RS_\varepsilon, \medskip \\
\left|\Psi_n\right| \geq C_\varepsilon(\N_*)^{-1}\left|\Phi^n\right| & \text{in} & \RS^{(0)}\cap\pi^{-1}\big\{|z|>1/\varepsilon\big\}.
\end{array}
\right.
\end{equation}
\end{theorem}

We prove Propositions~\ref{prop:genus1} and~\ref{prop:genus2} as well as Theorems~\ref{thm:N-S} and~\ref{thm:N-Sasymp} in Section~\ref{sec:N-S}.

\section{Asymptotics of Hermite-Pad\'e Approximants}
\label{sec:HP}

Below, we consider vector functions $\vec f:=\big(f_1,f_2\big)$ of the form
\begin{equation}
\label{appr-fun}
f_j(z) := \frac{1}{2\pi\mathrm{i}}\int_{F_j}\frac{\rho_j(x)}{x-z}\,\mathrm{d}x, \quad j\in\{1,2\},
\end{equation}
where $F_1=[-1,a]$ and $F_2=[-a,1]$, $a\in(0,1)$, and $\rho_j$ are holomorphic and non-vanishing in a neighborhood of $[-1,1]$. Additionally, we impose the following condition on the functions $\rho_j$:
\begin{condition}
\label{cond}
The ratio $\rho_2/\rho_1$ extends from $(-a,a)$ to a holomorphic and non-vanishing function
\begin{itemize}
\item in a domain that contains in its interior the closure of all the bounded components of the regions $\Omega_{ijk}$ in Case I, see Figure~\ref{Oms}(a);
\item in a domain whose complement is compact and belongs to the right-hand component of $\Omega_{021}$ in Cases II and IIIa, see Figures~\ref{Oms}(b,c);
\item in the extended complex plane, i.e., the ratio is a non-zero constant, in Case IIIb.
\end{itemize}
\end{condition}

Condition~\ref{cond} needs to be regarded in the following context.  As conjectured in \cite{Ap08b} and supported by the results for Nikishin systems, given $\vec f$ of the form \eqref{appr-fun}, the appropriate Riemann surface $\RS$ must depend on the analytic continuation of $\rho_2/\rho_1$ from $F_1\cap F_2$. In this work, on the other hand, we fixed the surface by considering \eqref{h}, which necessitates a condition on the continuation of $\rho_2/\rho_1$ of the above type.

In what follows, we assume that $a\in(0,1)$ is fixed and the vector $\vec f=(f_1,f_2)$ is given by  \eqref{appr-fun}, where $\rho_1$ and $\rho_2$ satisfy Condition~\ref{cond}; $\Psi_n$ are the Nuttall-Szeg\H{o} functions of Theorem~\ref{thm:N-S} corresponding to $a$ and the functions $\rho_1$ and $\rho_2$, constructed on the Riemann surface $\RS$ of $h$ realized as in Theorem~\ref{RS}; the constants $C_n$ are defined by
 \begin{equation}
 \label{Cn}
C_n^{-1}:=\lim_{z\to\infty} \Psi_n^{(0)}(z)z^{-2n};
 \end{equation}
$Q_{\vec n}$ are the monic denominators of the Hermite-Pad\'e approximants $\vec \pi_{\vec n}$ to $\vec f$, which are also multiple orthogonal polynomials with respect to the weights $\rho_1$ and $\rho_2$, see \eqref{ortho}), and $R_{\vec n}^{(i)}$, $i\in\{1,2\}$, are the linearized errors of approximation corresponding to the multi-indices $\vec n=(n,n)$, $n\in \N$, see \eqref{pi_vecn}--\eqref{R_vecn}.

\begin{figure}[!ht]
\centering
\subfigure[Case I]{\includegraphics[scale=.5]{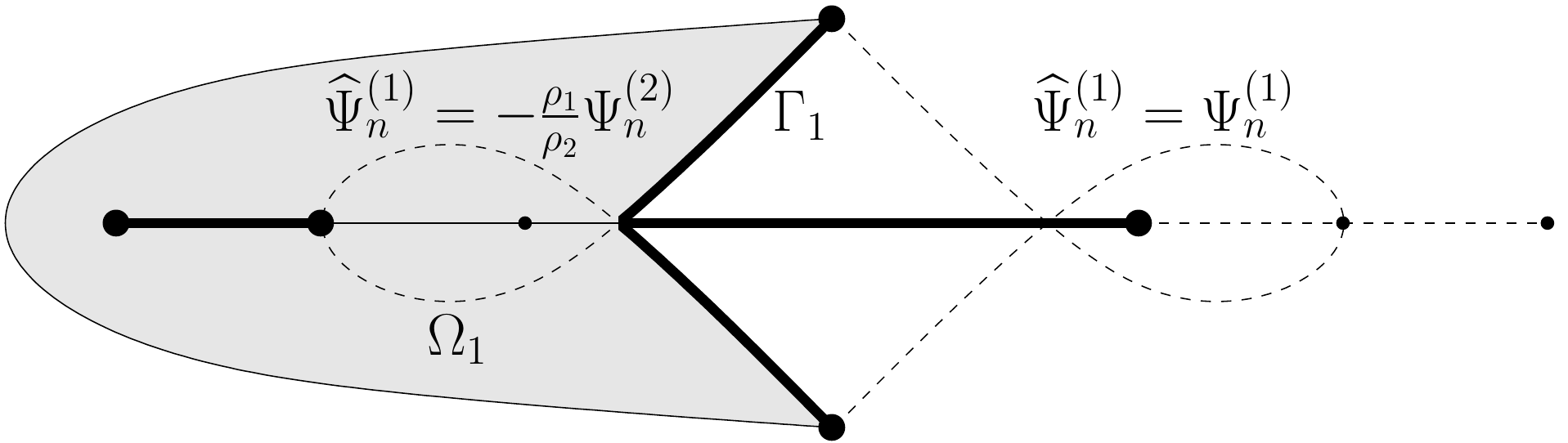}}
\subfigure[Case I]{\includegraphics[scale=.5]{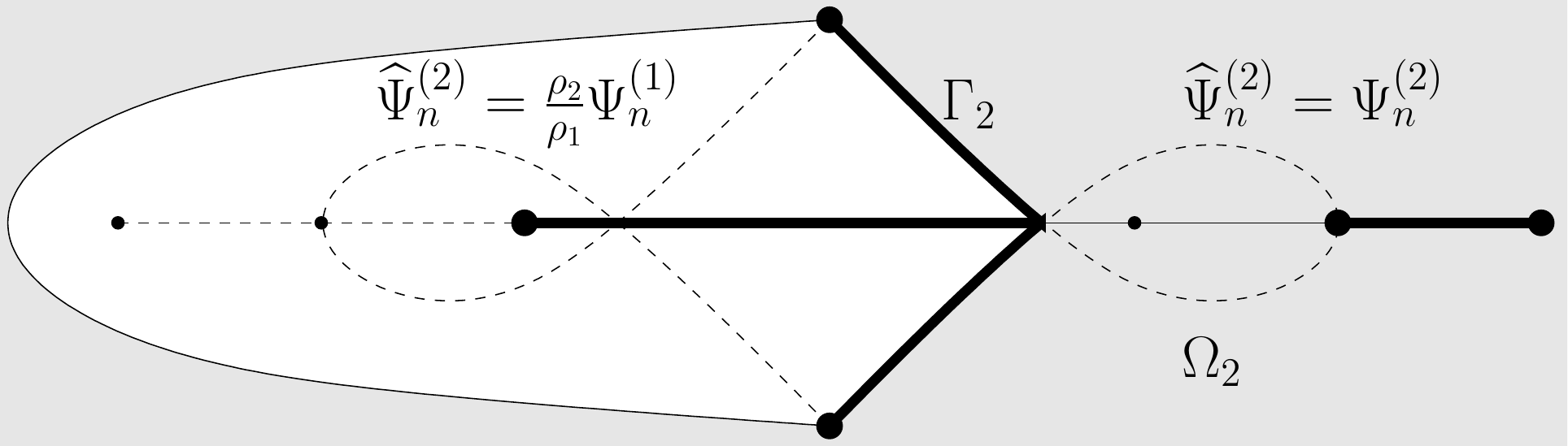}}
\subfigure[Case II]{\includegraphics[scale=.5]{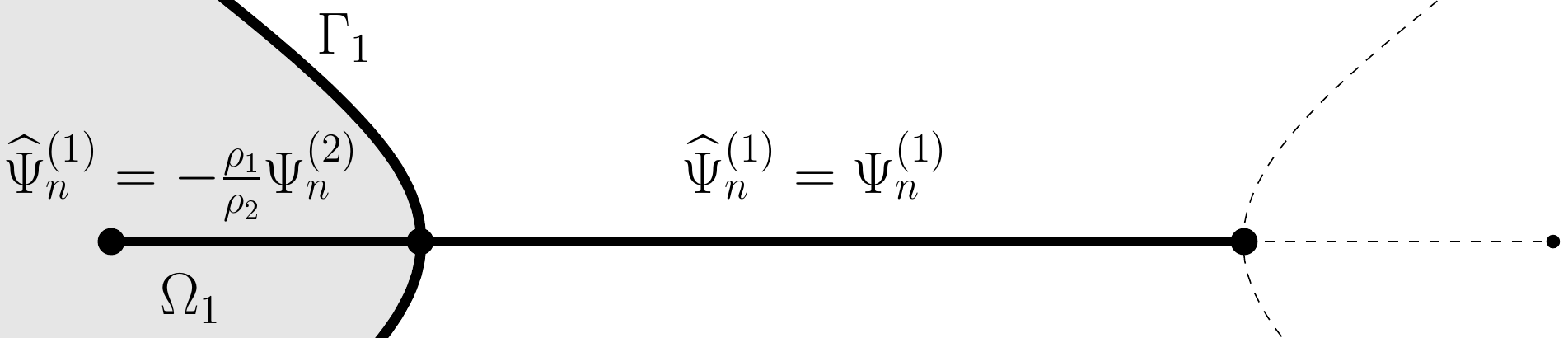}}
\subfigure[Case II]{\includegraphics[scale=.5]{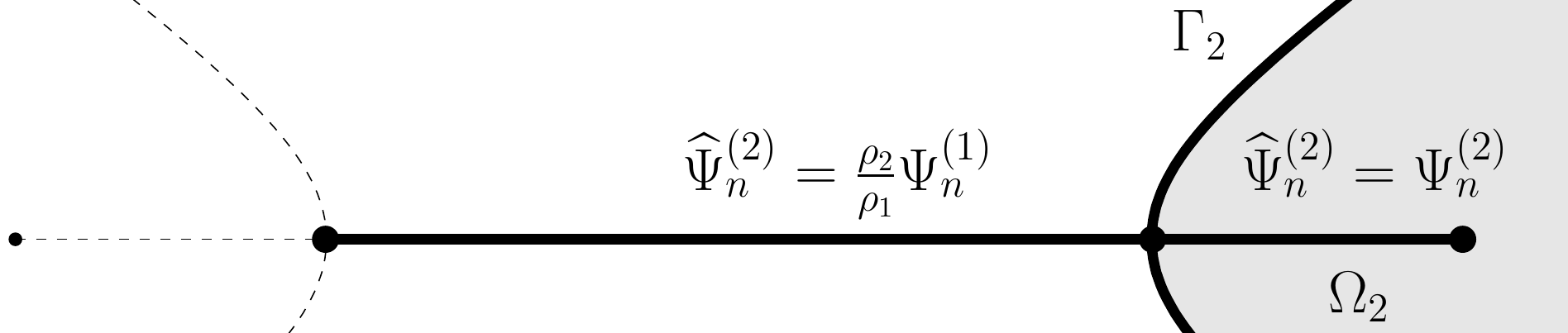}}
\subfigure[Case IIIa]{\includegraphics[scale=.5]{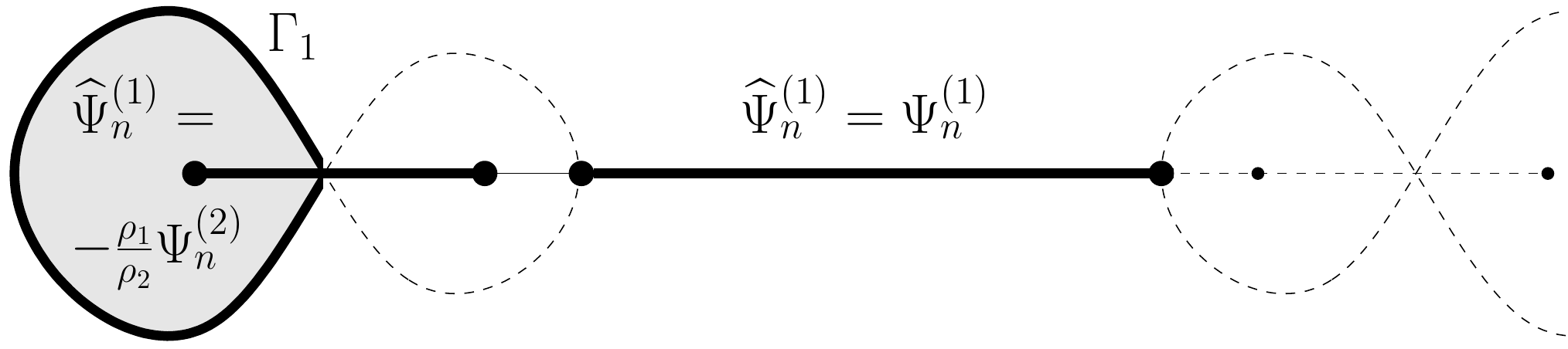}}
\subfigure[Case IIIa]{\includegraphics[scale=.5]{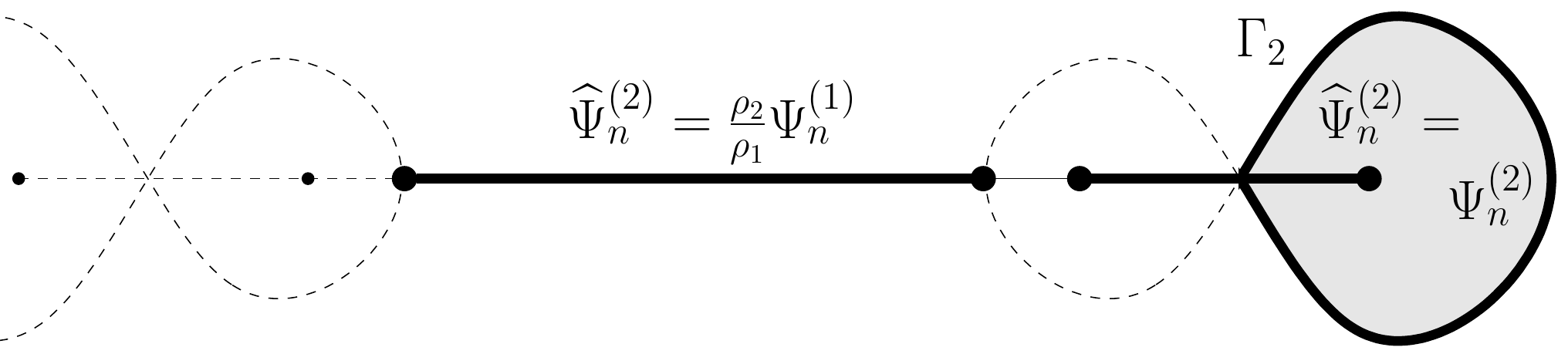}}
\caption{\small The sets $\Gamma_i$ (non-horizontal bold lines), the domains $\Omega_i$ (shaded regions), and the boundary of the domain of holomorphy of $\widehat\Psi_n^{(1)}$ (bold lines).}
\label{fig:results}
\end{figure}

To describe the asymptotics of $R_{\vec n}^{(i)}$, it will be convenient to put
\[
\left\{
\begin{array}{lll}
\Gamma_1 &:=& \Gamma_{12}^-\setminus\Delta_0, \medskip \\
\Gamma_2 &:=& \Gamma_{12}^+,
\end{array}
\right. \quad \text{and} \quad \left\{
\begin{array}{lll}
\Omega_1 &:=& \big(\overline\Omega_{021}^-\cup\overline\Omega_{201}^-\big)\setminus\big(\Gamma_{12}\cup\Delta_1\cup\Delta_2\big), \medskip \\
\Omega_2 &:=& \big(\overline\Omega_{021}^+\cup\overline\Omega_{201}^+\cup\overline\Omega\big)\setminus\big(\Gamma_{12}^+\cup\Delta_1\cup\Delta_2\big),
\end{array}
\right.
\]
where $\Delta_1$ and $\Delta_2$ were introduced in \eqref{chains}, $\Gamma_{12}^\pm:=\Gamma_{12}\cap\{\pm\re(z)>0\}$, $\Omega_{ijk}^\pm:=\Omega_{ijk}\cap\{\pm\re(z)>0\}$, and $\Omega$ is present only in Case I and is equal to the unbounded component of $\Omega_{012}^-$, see Figures~\ref{fig:results} and~\ref{Oms}. Observe that $\Omega_1=\Omega_2=\varnothing$ in Case IIIb. Define
\begin{equation}
\label{PsiHat}
\widehat\Psi_n^{(1)} :=
\left\{
\begin{array}{rl}
-\frac{\rho_1}{\rho_2}\Psi_n^{(2)} & \text{in } \Omega_1, \medskip \\
\Psi_n^{(1)} & \text{otherwise},
\end{array}
\right. \quad \text{and} \quad \widehat\Psi_n^{(2)} :=
\left\{
\begin{array}{rl}
\Psi_n^{(2)} & \text{in } \Omega_2, \medskip \\
\frac{\rho_2}{\rho_1}\Psi_n^{(1)} & \text{otherwise}.
\end{array}
\right.
\end{equation}
In Case I, it follows from \eqref{BVP} that $\widehat\Psi_n^{(i)}$ is the analytic continuation of $\Psi_n^{(i)}$ across $\Delta_0$ until $\Gamma_i$, which always exists by the analyticity of $\rho_i$'s. In Cases II and III it still holds by \eqref{BVP} that $\widehat\Psi_n^{(i)}$ is an analytic continuation of $\Psi_n^{(i)}$. However, this time one needs to continue $\Psi_n^{(1)}$ through $\Delta_1^\pm$ into $\RS^{(0)}$ and then through $\Delta_{21}^\mp$ into $\RS^{(2)}$ while $\Psi_n^{(2)}$ needs to be continued through $\Delta_{22}^\pm$ into $\RS^{(0)}$ and then through $\Delta_1^\mp$ into $\RS^{(1)}$.

For any $\delta>0$, we further define
\[
\mathcal N_\delta := \left\{z:~|\im(z)|<\delta,~\re(z)\in\Delta_1\cup\Delta_2,~\dist\left(\re(z),\{\text{endpoints of }\Delta_1\text{ and }\Delta_2\}\right)>\delta\right\},
\]
and
\[
\mathcal N_\delta^{(i)} := \left\{
\begin{array}{ll}
\left\{z=t+x:~x\in(-\delta,\delta),~t\in \Gamma_i,~\delta<|\im(t)|<c-\delta\right\} & \text{in Case I}, \medskip \\
\left\{z=t+x:~x\in(-\delta,\delta),~t\in \Gamma_i,~\delta<|\im(t)|<1/\delta\right\} & \text{in Case II}, \medskip \\
\left\{z:~\dist(z,\Gamma_i)<\delta,~|\im(z)|>\delta \text{ when } |\re(z)|<1\right\} & \text{in Case III}.
\end{array}
\right.
\]
Finally, let us introduce the following notation. Given a sequence of functions $F_n$, a sequence of finite multi-sets $X_n$, and positive numbers $\epsilon_n$, we write
\[
F_n=\mathcal{O}\left(\epsilon_n;X_n\right) \quad \Leftrightarrow \quad F_n(z)\prod_{x\in X_n}\frac{|z-x|}{\sqrt{(1+|z|^2)(1+|x|^2)}} =\mathcal{O}(\epsilon_n).
\]
Then the following theorem holds.

\begin{theorem}
\label{thm:main}
Let $\N_*$ be a subsequence from either Proposition~\ref{prop:genus1} or Proposition~\ref{prop:genus2} (depending on whether $a\geq 1/\sqrt2$ or $a<1/\sqrt2$), and $i\in\{1,2\}$. Then
\begin{equation}
\label{main-A}
\left\{
\begin{array}{lll}
Q_{\vec n} &=& C_n\Psi_n^{(0)}\left(1 + \mathcal{O}\big(\epsilon_n;X_n\big) \right), \medskip \\
R_{\vec n}^{(i)} &=& C_n\widehat\Psi_n^{(i)}\left(1 + \mathcal{O}\big(\epsilon_n;X_n^{(i)}\big)\right),
\end{array}
\right. \quad n\in \N_*,
\end{equation}
locally uniformly in $\overline\C\setminus(\Delta_1\cup\Delta_2)$ and $\overline\C\setminus\big(F_i\cup \Gamma_i\big)$, respectively, where $X_n$ is the multi-set of zeros of $\Psi_n^{(0)}$ in $\overline\C\setminus(\Delta_1\cup\Delta_2)$, $X_n^{(i)}$ is the multi-set of zeros of $\widehat\Psi_n^{(i)}$ in $\overline\C\setminus\big(F_i\cup \Gamma_i\big)$, and\footnote{In Cases I and III one has $\epsilon_n=n^{-1}$ and $\epsilon_n=n^{-1/6}$ in Case II.} $\epsilon_n\to 0$. Moreover, there exists $\delta_0>0$ such that for all $\delta\in(0,\delta_0)$ one has
\begin{equation}
\label{main-B}
\left\{
\begin{array}{lll}
Q_{\vec n} &=& C_n\Psi_{n+}^{(0)}\left(1 + \mathcal{O}\big(\epsilon_n;X_{n+}\big) \right) + C_n\Psi_{n-}^{(0)}\left(1 + \mathcal{O}\big(\epsilon_n;X_{n-}\big) \right), \medskip \\
R_{\vec n}^{(i)} &=& C_n\widehat\Psi_{n+}^{(i)}\left(1 + \mathcal{O}\left(\epsilon_n;X_{n+}^{(i)}\right)\right) + C_n\widehat\Psi_{n-}^{(i)}\left(1 + \mathcal{O}\left(\epsilon_n;X_{n-}^{(i)}\right)\right),
\end{array}
\right. \quad n\in \N_*,
\end{equation}
locally uniformly in $\mathcal N_\delta$ and $\mathcal N_\delta^{(i)}$, respectively, where $\Psi_{n\pm}^{(0)}$ are functions holomorphic in $\mathcal N_\delta$ that coincide with $\big(\Psi_n^{(0)})^\pm$ on $\Delta_1^\circ\cup\Delta_2^\circ$  and $X_{n\pm}$ are their multi-sets of zeros in $\mathcal N_\delta$, and $\widehat\Psi_{n\pm}^{(i)}$ are functions holomorphic in $\mathcal{N}_\delta^{(i)}\setminus F_i$ that coincide with $\big(\widehat\Psi_n^{(i)}\big)^\pm$ on $\Gamma_i^\circ\setminus F_i$ and $X_{n\pm}^{(i)}$ are their multi-sets of zeros in $\mathcal{N}_\delta^{(i)}$.
\end{theorem}

It follows from Theorems~\ref{thm:Omegas} and~\ref{thm:N-Sasymp} that in Cases I and III we observe the pushing effect, i.e, $a\neq b$, and the presence of divergence regions (both phenomena are observed in Angelesco systems). Indeed, according to \eqref{R_vecn} and \eqref{main-A}, one has locally uniformly in $\overline\C \setminus\big([-1,1]\cup \Gamma_i\big)$ that
\[
f_i - \pi_{\vec n}^{(i)} = \frac{\widehat\Psi_n^{(i)}}{\Psi_n^{(0)}} \frac{1 + \mathcal{O}\big(\epsilon_n;X_n^{(i)}\big)}{1 + \mathcal{O}\big(\epsilon_n; X_n\big)}.
\]
Hence, the error of approximation $f_i - \pi_{\vec n}^{(i)}$ is geometrically small in $\{|z|\geq1/\varepsilon\}$ for all $\varepsilon$ small enough. On the other hand, in Case I, the error $f_1 - \pi_{\vec n}^{(1)}$ is at least geometrically big on compact subsets of $\Omega_{201}^-$, see Figure~\ref{Oms}(a), (it could be infinite if the elements of $\mathcal{D}_n$ belong to $\RS^{(0)}$ and project onto this component) and the error $f_2 - \pi_{\vec n}^{(2)}$ is at least geometrically big on compact subsets of $\Omega_{201}^+$. In Case III, both components of the approximant diverge in $\Omega_{102}$ and in Case II there are no divergence domains.

As to the zeroes of the functions, it can be deduced from Theorem~\ref{thm:Omegas} and \eqref{main-B} that $Q_{\vec n}$ must vanish in $\mathcal N_\delta$ and $R_{\vec n}^{(i)}$ must vanish in $\mathcal N_\delta^{(i)}$, which is precisely the phenomenon of overinterpolation first observed in Nikishin systems.

\section{Multiple Orthogonal Polynomials}  \label{MOP}

The basis of our approach to asymptotics of Hermite-Pad\'e approximants lies in their connection with multiple orthogonal polynomials. It is quite simple to verify that if the functions $f_i$ are of the form \eqref{appr-fun}, then \eqref{R_vecn} is fulfilled if and only if
\begin{equation}
\label{ortho}
\int_{F_i}Q_{\vec n}(x)x^k\rho_i(x)\,\mathrm{d}x =0, \quad k\in\{0,\ldots,n_i-1\}.
\end{equation}
Moreover, the linearized error functions $R_{\vec n}^{(i)}$ admit the following integral representation:
\begin{equation}
\label{secondkind}
R_{\vec n}^{(i)}(z) = \frac{1}{2\pi\mathrm{i}}\int_{F_i}\frac{Q_{\vec n}(x)\rho_i(x)}{x-z}\,\mathrm{d}x.
\end{equation}

The analysis of the system \eqref{ortho}--\eqref{secondkind} then proceeds via its reformulation as a matrix Riemann-Hilbert problem. This fundamental fact in the theory of orthogonal polynomials was first revealed by Fokas, Its, and Kitaev \cite{FIK91,FIK92} and the extension to multiple orthogonal polynomials was given in \cite{GerKVA01}. Set $\vec n_1:=(n-1,n)$ and $\vec n_2:=(n,n-1)$, and assume that the index $n$ is such that
\begin{equation}
\label{normality}
m_n^{(i)}R_{\vec n_i}^{(i)} = z^{-n} + \cdots
\end{equation}
for some constants $m_n^{(i)}$. Under condition \eqref{normality}, the matrix
\begin{equation}
\label{Y}
{\boldsymbol Y}:=
\left(
\begin{array}{ccc}
Q_{\vec n} & R_{\vec n}^{(1)} & R_{\vec n}^{(2)} \smallskip \\
m_n^{(1)}Q_{\vec n_1} & m_n^{(1)}R_{\vec n_1}^{(1)} & m_n^{(1)}R_{\vec n_1}^{(2)}  \smallskip \\
m_n^{(2)}Q_{\vec n_2} & m_n^{(2)}R_{\vec n_2}^{(1)} & m_n^{(2)}R_{\vec n_2}^{(2)}
\end{array}
\right)
\end{equation}
solves the following Riemann-Hilbert Problem (\rhy):
\begin{itemize}
\label{rhy}
\item[(a)] ${\boldsymbol Y}$ is analytic in $\C\setminus[-1,1]$ and
\[
\lim_{z\to\infty} {\boldsymbol Y}(z)\ \diag\left(z^{-2n},z^n,z^n\right) = \boldsymbol I,
\]
where $\diag(\cdot,\cdot,\cdot)$ is the diagonal matrix and ${\boldsymbol I}$ is the identity matrix;
\item[(b)] ${\boldsymbol Y}$ has continuous traces on $(-1,1)\setminus\{\pm a\}$ that satisfy $\displaystyle {\boldsymbol Y}_+ = {\boldsymbol Y}_- \boldsymbol J(\chi_1\rho_1 , \chi_2\rho_2)$, where
\begin{equation}
\label{J}
\boldsymbol J(x,y) =
\left(
\begin{array}{ccc}
1 & x & y \smallskip \\
0 & 1 & 0  \smallskip \\
0 & 0 & 1
\end{array}
\right)
\end{equation}
and $\chi_i$ is the indicator function of $F_i$;
\item[(c)] $\boldsymbol Y(z)=\boldsymbol{\mathcal{O}}\left(\log|z-e|\right)$ as $[-1,1]\not\ni z\to e$ for $e\in\{\pm1,\pm a\}$\footnote{In fact, in each case the entries of at least two columns remain bounded. However, the above simplification does not affect the forthcoming analysis.}.
\end{itemize}
Vice versa, if \hyperref[rhy]{\rhy} is solvable, then the solution is necessarily of the form \eqref{Y} and \eqref{normality} holds. To prove Theorem~\ref{thm:main}, we then follow the framework of the non-linear steepest descent method for matrix Riemann-Hilbert problems, first introduced in the $2\times2$ case by Deift and Zhou \cite{DZ93}. The proof of Theorem~\ref{thm:main} is carried out in  Sections~\ref{sec:CaseI}--\ref{sec:CaseIII}.

\section{Geometry}
\label{sec:Geom}

This section is devoted to proving Theorems~\ref{thm:RS} and~\ref{thm:Omegas} together with Proposition~\ref{prop:periods}. In Section~\ref{ssec:Real-RS} we establish that $\RS$ can indeed be realized as in Figure~\ref{RS}. In Section~\ref{ssec:param} we justify the the choice of the parameter $p$, thus, finishing the proof of Theorem~\ref{thm:RS}, while simultaneously proving Proposition~\ref{prop:periods}. Finally, we prove Theorem~ \ref{thm:Omegas} in Section~\ref{ssec:domains}.

\subsection{Realization of $\RS$}
\label{ssec:Real-RS}

It can be readily computed that the discriminant of \eqref{h} is equal to
\begin{equation}
\label{discriminant}
\begin{array}{lll}
D(z) &=& 108A(z)\big[B_2^3(z)-A(z)B_1^2(z)\big] \medskip \\
&=& 108 A(z)\big[(1+a^2-3p^2)z^4+(3p^4-a^2)z^2-p^6\big].
\end{array}
\end{equation}

Assume first that we are in Case I, i.e.,
\begin{equation}
\label{6.1_I}
a\in\left(0,\frac1{\sqrt2}\right) \quad \text{and} \quad p\in\left(a,\sqrt{\frac{1+a^2}3}\right).
\end{equation}
Since the polynomial $B_2^3-AB_1^2$, which is symmetric and of degree 4, is negative at the origin and has positive leading coefficient, it follows that it has four zeros, which we denote by $\pm b$ and $\pm\mathrm{i}c$, where $b,c>0$. Furthermore, $b\in(a,p)$ since $B_2^3-AB_1^2$ is positive at $p$ and negative at $a$. Observe that if a point is a branch point of $h$ of order 3 (all three branches coincide) and $h$ is finite at this point, then necessarily all three branches are equal to zero there. Hence, $\pm b,\pm\mathrm{i}c$ are branch points of order 2. Furthermore, the first equation in \eqref{branches} implies that neither of the points $\pm1,\pm a$ can be a pole of one the branches while the second equation implies that the branches cannot have a cubic root singularity there. Hence, the points $\pm1,\pm a$ are branch points of order 2, two branches are infinite at them and one is finite.

It follows from the above discussion that we can analytically continue the branches $h_k$ so that the inequalities in \eqref{label} hold for $x>1$. At 1 two branches blow up and, of course, all three add up to zero. This is possible only if $h_0$ tends to $\infty$, $h_2$ tends to $-\infty$, and $h_1$ remains bounded. Thus, 1 is a branch point joining $\RS^{(0)}$ and $\RS^{(2)}$. It can easily be seen from \eqref{h} that all the branches must satisfy
\begin{equation}
\label{conj-sym}
h_k(z)=\overline{h_k(\overline z)}.
\end{equation}
Therefore, all the branch cuts must be conjugate-symmetric. Thus, the branch cut starting at 1 must end at $b$. That is, $h_{0\pm}=h_{2\mp}$ on $(b,1)$, which immediately implies
\begin{equation}
\label{likearoot}
(h_0-h_2)_+ + (h_0-h_2)_-  \equiv 0
\end{equation}
on $(b,1)$. Moreover, \eqref{conj-sym} implies that the traces above are purely imaginary. As $h_0,h_2$ are unbounded near 1 and bounded near $b$, one has
\begin{equation}
\label{dif1}
(h_0-h_2)^2(x) = \frac{x-b}{x-1}f_1^2(x),
\end{equation}
where $f_1$ is holomorphic, non-vanishing, and real on $(a,\infty)$.  As $f_1$ must be positive for $x>1$ and is non-vanishing,  it is, in fact, positive for $x>a$. Therefore, $(h_0-h_2)(x) > 0$, $x\in(a,b)\cup(1,\infty)$. Moreover, since $h_1$ is real and non-vanishing for $x>a$ and is negative for $x$ large enough by \eqref{hatinfinity}, $h_1$ is negative for all $x>a$. Thus,
\begin{equation}
\label{arrange1}
h_0(x) > h_2(x) > h_1(x), \quad x\in(a,b).
\end{equation}
As in the case of 1, two branches are infinite at $a$ and one is finite. The inequalities in \eqref{arrange1} imply that the unbounded branches are $h_0$ and $h_1$ and therefore $a$ is a branch point between $\RS^{(0)}$ and $\RS^{(1)}$. The branch cut is $(-a,a)$ since it must be along the real axis. That is, $h_{0\pm}=h_{1\mp}$ on $(-a,a)$ and
\begin{equation}
\label{likearoot1}
(h_0-h_1)_+ + (h_0-h_1)_-  \equiv 0
\end{equation}
on $(-a,a)$ where the traces are purely imaginary. As $h_0,h_1$ are unbounded near $\pm a$, one has
\begin{equation}
\label{dif2}
(h_0-h_1)^2(x) = \frac{f_a^2(x)}{x^2-a^2},
\end{equation}
where $f_a$ is holomorphic, non-vanishing, and real on $(-b,b)$. The function $f_a$ is positive on $(a,b)$ according to \eqref{arrange1} and therefore is positive on $(-b,b)$. Hence, $(h_0-h_1)(x) < 0$ for  $x\in(-b,-a)$. Since both branches blow up at $-a$, one has
\begin{equation}
\label{arrange2}
h_1(x) > h_2(x) > h_0(x), \quad x\in(-b,-a).
\end{equation}
Furthermore, it is quite simple to deduce from \eqref{h} that the branches $h_k$ must satisfy $h_k(-x)=-h_{j_k}(x)$ for all $k\in\{0,1,2\}$, where $j_k\in\{0,1,2\}$ depends on $k$. According to \eqref{arrange1} and \eqref{arrange2}, the considered continuations of $h_k$ satisfy $j_k=k$ for $|x|\in(a,b)$. Therefore, $-b$ is a branch point for $h_0$ and $h_2$ as well as $-1$. In particular, \eqref{likearoot} remains valid for $x\in(-1,-b)$. Hence,
\begin{equation}
\label{dif3}
(h_0-h_2)^2(x) = \frac{x+b}{x+1}f_{-1}^2(x),
\end{equation}
where $f_{-1}$ is holomorphic, non-vanishing, and real on $(d,-a)$ for some $d<-1$. From \eqref{arrange2}, we know that $f_{-1}$ is negative on $(-b,-a)$ and therefore it is negative on $(d,-a)$. Thus, $(h_0-h_2)(x) < 0$ for $x\in(d,-1)$. This and the blowing up of $h_0$ and $h_2$ at $-1$ imply that
\begin{equation}
\label{arrange3}
h_2(x) > h_1(x) > h_0(x), \quad x\in(d,-1).
\end{equation}
On the other hand, write
\begin{equation}
\label{foromega1}
h_i(z) = -\frac1z + \frac{\alpha_i}{z^2}+\cdots
\end{equation}
for $i\in\{1,2\}$. Plugging this expansion into \eqref{h} and considering the $1/z$ term on the left-hand side, we get that
\begin{equation}
\label{foromega2}
\alpha_i^2 = \frac{1+a^2}3 - p^2>0.
\end{equation}
This means that the inequality $h_1(x)>h_2(x)$ holds for all $|x|$ large enough. That is,
\[
h_1(x) > h_2(x) > h_0(x), \quad x\in(-\infty,d^\prime),
\]
where necessarily $d^\prime<d$. As there are no branch points between $d^\prime$ and $d$, there should be a branch cut passing between them and this cut should necessarily be between $\RS^{(1)}$ and $\RS^{(2)}$. In other words, $\pm\mathrm{i}c$ are branch points of $h_1$ and $h_2$. This finishes the proof of the claim that $\RS$ can be realized as in Figure~\ref{RS}(a) in Case I.

Assume now that we are in Case II, i.e.,
\[
a=p=1/\sqrt2.
\]
Then the discriminant of \eqref{h} is equal to $D(z) = 27(z^2-1)(z^2-1/2)^2$. As before, $\pm1$ are branch points of order 2. Furthermore, we know from the third identity in \eqref{branches} that some of the branches are unbounded near $\pm1/\sqrt{2}$. However, if the branching were of order 2, the left-hand side of the second relation in \eqref{branches} would be unbounded near $\pm1/\sqrt{2}$, but it is bounded there. Hence, $\pm1/\sqrt{2}$ are branch points of order 3.

Plugging a power expansion for $h_i$, $i\in\{1,2\}$, into \eqref{h}, one can compute that
\begin{equation}
\label{case2hi}
h_i(z) = -\frac1z + \frac{\beta_i}{z^3} + \cdots,
\end{equation}
where the $\beta_i$'s are solutions of $\beta^2+\beta+1/6=0$. In particular, $(h_1-h_2)(z)z^3=1/\sqrt3+\mathcal{O}(1/z)$ for all $z$ large. Then, by repeating the initial steps of the analysis for Case I, we see that \eqref{label} extends for all $x>1$ and that
\begin{equation}
\label{arrange6}
h_2(x) > h_1(x) > h_0(x), \quad x\in(-\infty,-1).
\end{equation}
Therefore, $(-1,-1/\sqrt2)$ and $(1/\sqrt2,1)$ are the branch cuts for $h_0$ and $h_2$. Around the point $1/\sqrt2$ it can be directly verified that \eqref{h} is solved by the following Puiseux series:
\begin{equation}
\label{hxi}
h(z;\xi) = \frac1{\sqrt{1-z^2}}\left[-\xi\left(z-\frac1{\sqrt2}\right)^{-1/3}H^{1/3}(z) + \xi^2 \left(z-\frac1{\sqrt2}\right)^{1/3}H^{-1/3}(z) \right],
\end{equation}
where $\xi$ is any solution of $\xi^3=1$, $H(z) := \big(1+2z\sqrt{1-z^2}\big)/\big(2z+\sqrt2\big)$, and all the roots are principal. Since $h_1$ is negative and holomorphic for $x>1/\sqrt2$, one finds that $h_1(z) = h(z;1)$ locally around $1/\sqrt2$. On the other hand,
\[
(h_0-h_2)^2(x) = \frac{f^2(x)}{x^2-1},
\]
where $f(x)$ is non-vanishing and holomorphic for $|x|>1/\sqrt 2$. Since $f(x)>0$ for $x>1$, we get that $f(x)>0$ for $x>1/\sqrt2$ and
\begin{equation}
\label{dif4}
h_{0\pm}(x) = h_{2\mp}(x) = -\frac12h_1(x) \mp \mathrm{i}\frac{f(x)}{\sqrt{1-x^2}}, \quad x\in\big(1/\sqrt 2,1\big).
\end{equation}
Thus, $h_{0+}$ has values in the fourth quadrant and $h_{0-}$ has values in the first quadrant. As the first summand in \eqref{hxi} is dominant around $1/\sqrt2$, we can conclude that
\begin{equation}
\label{Htoh}
h_0(z) = \left\{
\begin{array}{ll}
h\big(z;e^{2\pi\mathrm{i}/3}\big), & \im(z)>0, \smallskip \\
h\big(z;e^{4\pi\mathrm{i}/3}\big), & \im(z)<0,
\end{array}\right.
\quad \text{and} \quad  h_2(z) = \left\{
\begin{array}{ll}
h\big(z;e^{4\pi\mathrm{i}/3}\big), & \im(z)>0, \smallskip \\
h\big(z;e^{2\pi\mathrm{i}/3}\big), & \im(z)<0,
\end{array}\right.
\end{equation}
locally around $1/\sqrt2$. From this it is easy to see that $h_2$ is holomorphic across the interval $\big(-1/\sqrt2,1/\sqrt2\big)$ and this interval is the branch cut for $h_1$ and $h_0$.

Finally, assume that we are in Case III, that is,
\[
a\in\left(\frac1{\sqrt2},1\right) \quad \text{and} \quad p=\sqrt{\frac{1+a^2}3}<a.
\]
In this case $B_2^3-AB_1^2$ is a polynomial of degree 2 which has two roots $\pm b$ satisfying $b\in(p,a)$. Exactly as in Case I, we see that all the branch points, namely $\{\pm1,\pm a,\pm b\}$, are of order 2. Furthermore, the same reasoning as in Case I gives that \eqref{likearoot} holds on $(a,1)$. Since $h_0$ and $h_2$ are unbounded at both $1$ and $a$, we get that
\begin{equation}
\label{dif5}
(h_0-h_2)^2(x) = \frac{f_1^2(x)}{(x-1)(x-a)},
\end{equation}
where $f_1$ is holomorphic, non-vanishing, and real on $(b,\infty)$. Since $f_1$ is positive for $x>1$ and hence for $x>b$, we can conclude that $(h_0-h_2)(x)<0$ for $x\in(b,a)$. The blowing up of $h_0$ and $h_2$ at $a$ implies that
\begin{equation}
\label{arrange4}
h_2(x) > h_1(x) > h_0(x), \quad x\in(b,a).
\end{equation}
It further follows from the third equation in \eqref{branches} that two branches of $h$ are negative and one branch is positive on $(b,a)$. The inequalities in \eqref{arrange4} show that the negative branches are $h_0$ and $h_1$. This, in turn, implies that $b$ is a branch point of $h_0$ and $h_1$ and so is $-b$. Hence, \eqref{likearoot1} holds on $(-b,b)$. As all the branches are bounded at $\pm b$, we have that
\begin{equation}
\label{dif7}
(h_0-h_1)^2(x) = (x^2-b^2)f_b^2(x)
\end{equation}
for $x\in(-a,a)$, where $f_b$ is holomorphic, non-vanishing, and is negative on $(b,a)$. Therefore, it is negative on the whole interval $(-a,a)$ and we get that $(h_0-h_1)(x) > 0$ for $x\in(-a,-b)$. As before, the third equation in \eqref{branches} implies that two branches of $h$ are positive and one branch is negative on $(-a,-b)$. Moreover, it also implies that $h_2(0)=0$. As $h_2(x)>0$ for $x>0$ and it has no other zeros, it must be the negative branch. Thus,
\begin{equation}
\label{arrange5}
h_0(x) > h_1(x) > h_2(x), \quad x\in(-a,-b).
\end{equation}
Now, as the branches that meet at $-a$ are unbounded, they must be $h_0$ and $h_2$ by \eqref{arrange5}. Hence, \eqref{likearoot} holds on $(-1,-a)$ and we have that
\begin{equation}
\label{dif6}
(h_0-h_2)^2(x) = \frac{f_{-1}^2(x)}{(x+1)(x+a)},
\end{equation}
where $f_{-1}$ is holomorphic, non-vanishing, and positive on $(-a,-b)$. As before, this means that $f_{-1}$ is positive for all $x<-b$ and \eqref{arrange6} holds in this case as well. This finishes the proof of the claim that $\RS$ can be realized as in Figure~\ref{RS}(c) in Case III.

\subsection{Choice of the Parameter}
\label{ssec:param}

We start with Cases II and III as they are much simpler.  Here we show that for the choice of the parameter $p$ as in Theorem~\ref{thm:RS} the condition \eqref{condition} is fulfilled. That is, the period of the Nuttall differential $\mathrm{d}\mathcal N$ over any given chain on $\RS$ is purely imaginary. The latter simplifies to showing that the periods of $\mathrm{d}\mathcal N$ on the cycles of a homology basis are purely imaginary. In fact, \eqref{not-rational2} implies that both periods of $\mathrm{d}\mathcal N$ are equal to $\pi\mathrm{i}$.

In proving \eqref{not-rational2} we shall rely on the following observation: one has that $h_k(-z)=-h_{j_k}(z)$ for some $j_k\in\{0,1,2\}$, which can be deduced immediately from \eqref{h}. In fact, we see from \eqref{label} and \eqref{arrange6} that $h_k(-z)=-h_k(z)$ for all $k\in\{0,1,2\}$. This implies that
\[
\omega_1 = -\frac1{2\pi\mathrm{i}}\int_{\Delta_{21}}(h_{0+}-h_{0-})(x)\,\mathrm{d}x = -\frac1{2\pi\mathrm{i}}\int_{\Delta_{22}}(h_{0+}-h_{0-})(x)\,\mathrm{d}x.
\]
The sum of the last two integrals is equal to 1 by Cauchy's residue theorem applied to $h_2$ (recall that $h_{0\pm}=h_{2\mp}$ on $\Delta$ and $h_2(z)=-1/z+\mathcal{O}(1/z^2)$). This gives the desired conclusion about $\omega_1$. Furthermore, from the choice of the $\ualpha$-cycle, see Figure~\ref{fig:HB}(b), one has for Case III that
\[
\tau_1 = \frac1{2\pi\mathrm{i}}\left( \int_a^bh_0(x)\,\mathrm{d}x + \int_b^{-b}h_{0+}(x)\,\mathrm{d}x + \int_{-b}^{-a}h_0(x)\,\mathrm{d}x + \int_{-a}^ah_2(x)\,\mathrm{d}x \right).
\]
As $h_0$ is an odd function, the sum of the first and the third integrals is equal to zero. Since $h_2$ is odd as well, the fourth integral is zero too. Using the symmetry considerations once more, we can get that
\[
\int_b^{-b}h_{0+}(x)\,\mathrm{d}x = \int_{-b}^bh_{0-}(x)\,\mathrm{d}x.
\]
Then by applying Cauchy's residue theorem to $h_1$ (recall that $h_1(z)=-1/z+\mathcal{O}(1/z^2)$), we see that
\begin{eqnarray*}
\tau_1 &=& \frac1{4\pi\mathrm{i}} \left( \int_b^{-b}h_{0+}(x)\,\mathrm{d}x + \int_{-b}^bh_{0-}(x)\,\mathrm{d}x \right) \\
 &=& \frac1{4\pi\mathrm{i}} \left( \int_b^{-b}h_{1-}(x)\,\mathrm{d}x + \int_{-b}^bh_{1+}(x)\,\mathrm{d}x \right) = 1/2.
\end{eqnarray*}
In Case II, we have that
\[
\tau_1 = \frac1{2\pi\mathrm{i}}\left( \int_{1/\sqrt2}^{-1/\sqrt2}h_{0+}(x)\,\mathrm{d}x + \int_{-1/\sqrt2}^{1/\sqrt2}h_2(x)\,\mathrm{d}x \right),
\]
and the conclusion $\tau_1=1/2$ follows from an analogous symmetry argument. This finishes the proof for Cases II and III of Theorem~\ref{thm:RS}.

In Case I we again start with \eqref{not-rational1}. The symmetry of the surface implies in this case that $h_0$ is an odd function, i.e., $h_0(-z)=-h_0(z)$, the functions $h_i$, $i\in\{1,2\}$, are odd within the bounded domain delimited by $\Delta_0$ and $-\Delta_0:=\{z:-z\in\Delta_0\}$, and $h_1(-z)=-h_2(z)$ for $z$ within the unbounded domain delimited by $\Delta_0$ and $-\Delta_0$.  Then it follows from \eqref{hatinfinity} and the Cauchy residue theorem that
\begin{eqnarray*}
2 &=& -\frac1{2\pi\mathrm{i}}\int_{\Delta_2\cup\Delta_1}(h_{0+}-h_{0-})(x)\,\mathrm{d}x  \\
&=& -\frac1{2\pi\mathrm{i}}\int_{\Delta_2}(h_{0+}-h_{2+})(x)\,\mathrm{d}x - \frac1{2\pi\mathrm{i}}\int_{\Delta_1}(h_{0+}-h_{1+})(x)\,\mathrm{d}x  \\
&=& \frac1{2\pi}\int_{\Delta_{21}}\big(-f_{-1}(x)\big)\sqrt{\frac{-b-x}{x+1}}\,\mathrm{d}x + \frac1{2\pi}\int_{\Delta_1}\frac{f_a(x)}{\sqrt{z^2-x^2}}\,\mathrm{d}x + \frac1{2\pi}\int_{\Delta_{22}} f_1(x)\sqrt{\frac{x-b}{1-x}}\,\mathrm{d}x ,
\end{eqnarray*}
where we used the notation from \eqref{dif1}, \eqref{dif2}, and \eqref{dif3}. The analysis right after \eqref{dif1}, \eqref{dif2}, and \eqref{dif3} implies that each of the last three integrals is positive. From the definition of our homology basis, it can easily be seen that the first integral is equal to $\omega_1$, the second one is equal to $\omega_2$, and the third is equal to the first as $h_0$ is an odd function. That is, $\omega_i>0$ and $2\omega_1+\omega_2=2$. The latter clearly implies the first relation in \eqref{not-rational1} and the fact that $\omega=\omega_1\in(0,1)$. Applying now the Cauchy residue theorem to $h_2$, we get that
\begin{eqnarray*}
1 & = & \frac1{2\pi\mathrm{i}} \int_{\Delta_2}(h_{2+}-h_{2-})(x)\,\mathrm{d}x + \frac1{2\pi\mathrm{i}} \int_{\Delta_0}(h_{2+}-h_{2-})(x)\,\mathrm{d}x  \\
& = & 2\omega_1 + \frac1{2\pi\mathrm{i}} \int_{\Delta_0}(h_{2+}-h_{1+})(x)\,\mathrm{d}x,
\end{eqnarray*}
where $\Delta_0$ is oriented from $-\mathrm{i}c$ to $\mathrm{i}c$. Hence, we can conclude that $\omega_1\in(1/2,1)$, if we show that the last integral is negative. Since $h_{2+}-h_{1+}$ is a trace of a holomorphic function, we can deform the path of integration to get
\begin{eqnarray}
\frac1{2\pi\mathrm{i}} \int_{\Delta_0}(h_{2+}-h_{1+})(x)\,\mathrm{d}x &=& -\frac1{2\pi\mathrm{i}}\left(\int_{-\mathrm{i}\infty}^{-\mathrm{i}c}+\int_{\mathrm{i}c}^{\mathrm{i}\infty}\right)(h_2(t)-h_1(t))\,\mathrm{d}t \nonumber \\
& = & - \frac1{2\pi}\left(\int_{-\infty}^{-c}+\int_c^{\infty}\right)(h_2(\mathrm{i}x)-h_1(\mathrm{i}x))\,\mathrm{d}x \nonumber \\
\label{Delta0cycle}
& = & - \frac1\pi \int_c^{\infty} (h_2(\mathrm{i}x)-h_1(\mathrm{i}x))\,\mathrm{d}x,
\end{eqnarray}
where we used the symmetry $h_2(-z)=-h_1(z)$ to get the last equality. Notice also that the branches satisfy $h_k(z)=\overline{h_k(\bar z)}$ as follows directly from \eqref{h}. Hence,
\[
h_2(\mathrm{i}x)-h_1(\mathrm{i}x) = h_2(-\mathrm{i}x)-h_1(-\mathrm{i}x) = \overline{h_2\big(-\overline{\mathrm{i}x}\big)-h_1\big(-\overline{\mathrm{i}x}\big)} =  \overline{h_2(\mathrm{i}x)-h_1(\mathrm{i}x)}
\]
for $x>c$. Hence, this difference is real there. Moreover, $h_2(\mathrm{i}x)-h_1(\mathrm{i}x)>0$ for all $x$ large as follows from \eqref{foromega1} and \eqref{foromega2}. Since the difference $h_2-h_1$ can be equal to zero only at the branch points of $\RS$, $h_2(\mathrm{i}x)-h_1(\mathrm{i}x)>0$ for all $x>c$, which shows that the integral in \eqref{Delta0cycle} is negative as desired.

To prove the second relation in \eqref{not-rational1}, observe that our choice of the homology basis can be made so that the $\ualpha$-cycles are contained within the bounded domain delimited by $\Delta_0$ and $-\Delta_0$, see Figure~\ref{fig:HB}(a). Moreover, they can freely be deformed within the domain of holomorphy of $\mathrm{d}\mathcal N$. Thus, it follows from \eqref{GreenDiff} and \eqref{GreenPeriods} that
\begin{eqnarray*}
2\pi\mathrm{i}\tau_1 & = &  \int_{-a}^{-b}h_0(x)\,\mathrm{d}x + \int_{-b}^0h_2(x)\,\mathrm{d}x + \int_0^{-\mathrm{i}c}h_2(x)\,\mathrm{d}x + \int_{-\mathrm{i}c}^0 h_1(x)\,\mathrm{d}x + \int_0^{-a}h_{0-}(x)\,\mathrm{d}x  \\
& = & \int_a^bh_0(x)\,\mathrm{d}x + \int_b^0h_2(x)\,\mathrm{d}x + \int_0^{\mathrm{i}c}h_2(x)\,\mathrm{d}x + \int_{\mathrm{i}c}^0 h_1(x)\,\mathrm{d}x + \int_0^ah_{0+}(x)\,\mathrm{d}x  \\
& = & -2\pi\mathrm{i}\tau_2,
\end{eqnarray*}
where we used the fact $h_k(-z)=-h_k(z)$ with the bounded domain delimited by $\Delta_0$ and $-\Delta_0$, which, in particular, implies that $h_{1-}(-x)=h_{1+}(x)$ for $x\in(-a,a)$. This finishes the proof of \eqref{not-rational1}. It only remains to show that there exists a choice of the parameter $p$ in \eqref{h} so that $\tau$ from \eqref{not-rational1} is real.

Showing that there is a choice of the parameter $p\in(a,\sqrt{(1+a^2)/3})$ such that $\tau$ in \eqref{not-rational1} is real, is equivalent to proving that $I=0$, where
\begin{equation}
\label{A.11-1}
 I:=\re\left(\oint_{\ualpha_2} h(\z)\,\mathrm dz\right).
\end{equation}
In order to prove \eqref{A.11-1} we consider the limiting values of $I=I(p)$ as $p\to\sqrt{(1+a^2)/3}$ and $p\to a$. If they have opposite signs, then clearly such a choice of $p$ is indeed possible.

We start with the case $p\to\sqrt{(1+a^2)/3}$. Define
\[
\left\{
\begin{array}{lll}
N_k(z;p) &:=& \re\left(\int_b^zh_k(x)\,\mathrm dx\right), \quad k\in\{0,2\}, \medskip \\
N_1(z;p) &:=& \re\left(\int_a^zh_1(x)\,\mathrm dx\right) + N_0(a;p),
\end{array}
\right. \quad z\in\big\{\re(z),\im(z)\geq 0\big\}.
\]
Then it obviously holds that $I(p) = N_1(\mathrm ic;p) - N_2(\mathrm ic;p)$. In Section~\ref{ssec:domains} further below we shall argue that
\begin{equation}
\label{sHWRIG}
N(\z;p^*) := \re\left(\int^\z h(\boldsymbol x;p^*)\,\mathrm dx\right), \quad p^*:= \sqrt{(1+a^2)/3},
\end{equation}
is a well defined harmonic function on the Riemann surface of $h(\cdot;p^*)$ and that
\begin{equation}
\label{sdvgads}
N\big(\infty^{(1)};p^*\big) - N\big(\infty^{(2)};p^*\big)<0.
\end{equation}
As $I(p)$ depends continuously on the parameter $p$, we can conclude that $I(p^*)<0$.

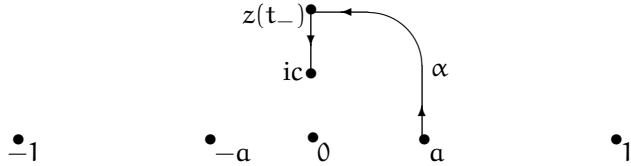
\begin{figure}[ht!]
\parbox{5in}{\centering
\unitlength 0.8pt
\linethickness{0.5pt}
\begin{picture}(350,125)
\put(40,67){$\bullet$}\put(38,60){$-1$}
\put(130,67){$\bullet$}\put(135,60){$-a$} 
\put(178,68){$\bullet$}\put(183,60){${\small 0}$}
\put(230,67){$\bullet$}\put(235,60){$a$}
\put(320,67){$\bullet$}\put(325,60){$1$}
\put(180,103){\line(0,1){30}}
\put(236,100){$\alpha$}
\put(180,70){\oval(104,120)[tr]}
\put(232,85){\vector(0,1){}}
\put(195,130){\vector(-1,0){}}
\put(177,128){$\bullet$}\put(148,125){$z(t_-)$}
\put(177,98){$\bullet$}\put(167,98){$\mathrm i c$}
\put(180,115){\vector(0,-1){}}
\end{picture}}
\vspace{-1.8cm}
\caption{Contour $\alpha$ for the limiting case $p=a<1/\sqrt2$.}
\label{Fig.gamma3}
\end{figure}

Let now $p=a$. In this case equation \eqref{h} becomes
\begin{equation}\label{A.15}
(z^{2}-1)(z^{2}-a^{2})\,h^{3}-3(z^{2}-a^{2})\,h-2z=0\;,\quad a^{2}<1/2.
\end{equation}
The branch points with projections $a$ and $b$ of the curve \eqref{A.15} merge together into a triple branch point with projection $a$ (this can be observed directly from \eqref{discriminant}). To verify that $I(a)>0$, we deform the cycle $\ualpha_2=\ualpha_2(a)$ into a cycle $\ualpha$ which is involution-symmetric and whose projection from $\RS^{(1)}$, say $\alpha$, is as on Figure~\ref{Fig.gamma3}.
That is, $\alpha$ emanates from $a$ into the first quadrant along some special arc $z(t)$, $t\in[0,t_{-}]$, and then proceeds along the imaginary axis from $z(t_-)$ down to $\mathrm ic$.

As in \eqref{Delta0cycle} one has $\overline{h_i(\bar z;a)}=h_i(z;a)$ and $h_1(-z)=-h_2(z)$ for all $z$ large (and hence for all $z\in(\mathrm ic,\mathrm i\infty)$), which leads to the same conclusion that $(h_1-h_2)(z)$ is real for  $z\in(\mathrm ic,\mathrm i\infty)$ and therefore the real part of the integral of this difference on $[\mathrm ic,\mathrm i\infty)$ is equal to zero. Hence, we get that
\[
I(a) = \re\left(\int_\alpha(h_1(z;a)-h_2(z;a))\,\mathrm dz\right).
\]
Therefore, in order to complete the proof of Theorem~\ref{thm:RS}, it remains to show that there exists an arc $z(t)$, such that
\begin{equation}
\label{A.16-1}
\left\{
\begin{array}{l}
z(0)=a, \quad \re(z(t_{-}))=0, \quad \im(z(t_{-}))>c, \medskip \\
\re\big((h_{1}-h_{2})(z(t))z^\prime(t)\big)>0, \quad t\in(0,t_{-}).
\end{array}
\right.
\end{equation}
To construct this arc we use a parametrization of the algebraic curve \eqref{A.15}. This parametrization was suggested in \cite{uApToulYa} and it has the form
\begin{equation}
\label{A.17}
\left\{\begin{array}{l}
h=-\mathrm i\left(\xi t+\displaystyle\frac{1}{\xi t}\right)C(t) \medskip \\
z=\displaystyle\frac{t^{3}}{1-t^{6}}\big(\tilde{a}_{-}(t)+\mathrm i\tilde{a}_{+}(t)\big)
\end{array}\right.,
\quad \xi^{3}=1, \quad t\in\mathbb{C},
\end{equation}
where
\begin{equation}
\label{A.18}
\tilde{a}_{\pm}(t):=\sqrt{\pm(1-a^{2})\mp a^{2}t^{\pm 6}}\;,\quad C(t):=\frac{\tilde{a}_{+}(t)+\mathrm it^{6}\tilde{a}_{-}(t)}{(1-a^{2})\,(1+t^{6})}.
\end{equation}
Since the expression for $h$ in \eqref{A.17} has the Cardano form
\[
h=\xi\,A+\xi^{-1}B,\quad A=-\mathrm itC,\quad B=-\frac{\mathrm i}{t}\,C,
\]
the equivalence of \eqref{A.17} and \eqref{A.15} can be written (using the Vieta relations) as
\[
\left\{\begin{array}{l}
-AB=C^{2}=-\displaystyle\frac{1}{z^{2}-1} \medskip \\
A^{3}+B^{3}=\mathrm iC^{3}\,(t^{3}+t^{-3})=\displaystyle\frac{2z}{(z^{2}-1)\,(z^{2}-a^{2})}.
\end{array}\right.
\]
With some work one can verify that plugging in the expression for $C(t)$ from \eqref{A.18} and the expression for $z(t)$ from \eqref{A.17} into the above equations produces an identity, which proves that \eqref{A.17} is indeed a parametrization of \eqref{A.15}, (for details see \cite{uApToulYa}). Observe that
\begin{equation}
\label{A.19}
\tilde a_\pm(t)\geq0, \quad t\in[0,t_{-}], \quad t_{-}:=\left(\frac{a^{2}}{1-a^{2}}\right)^{1/6},
\end{equation}
and that the parametrization \eqref{A.17} defines a Jordan arc with end points
\[
z(0)=a \quad \text{and} \quad z(t_{-})=\frac{\mathrm ia}{\sqrt{1-2a^{2}}}, \quad \text{while} \quad \mathrm ic=\frac{\mathrm ia^{2}}{\sqrt{1-2a^{2}}},
\]
As $a>a^2$, the first line of \eqref{A.16-1} is satisfied. The local analysis around the point $a$ shows that $h_1$ is given by \eqref{A.17} with  $\xi=(\mathrm i\sqrt{3}-1)/2$ and for $h_2$ one needs to select $\xi=(\mathrm -i\sqrt{3}-1)/2$. Then, after a tedious computation, we get
\[
(h_{1}-h_{2})(z(t))z^\prime(t)= \frac{3\sqrt{3}\,t}{t^{4}+t^{2}+1}\left(\frac{1-2a^2}{\tilde a_+\tilde a_-}-\mathrm i\right).
\]
Hence, the second line of \eqref{A.16-1} follows from \eqref{A.19} as $a^{2}<1/2$. This finishes the proof of Theorem~\ref{thm:RS} granted we can prove that $N(\z;p^*)$ in \eqref{sHWRIG} is a well defined harmonic function, which we do at the end of the upcoming subsection.

\subsection{The Regions $\Omega_{ijk}$}
\label{ssec:domains}

We start with two general observations that are consequences of the single-valuedness of $N(\z)$ on $\RS$. First, the regions $\Omega_{ijk}$ could be equivalently defined by
\[
\Omega_{ijk} := \big\{z:N_{j}(z) > N_{i}(z) > N_{k}(z)\big\},
\]
and this definition does not depend on the initial point of integration chosen in \eqref{condition} as changing the initial point results in adding the same constant to all $N_k$ simultaneously. Secondly, let $H$ be the analytic continuation of $h_1-h_2$ from the point at infinity. Then $H$ is an algebraic function and $\re(H)$ is a well defined harmonic function on the Riemann surface of $H$. Hence, the set $\Gamma$ is a subset of a projection of the zero level line of $\re(H)$ to $\overline\C$. As such it cannot be dense in an open set.

\begin{figure}[!ht]
\centering
\includegraphics[scale=.6]{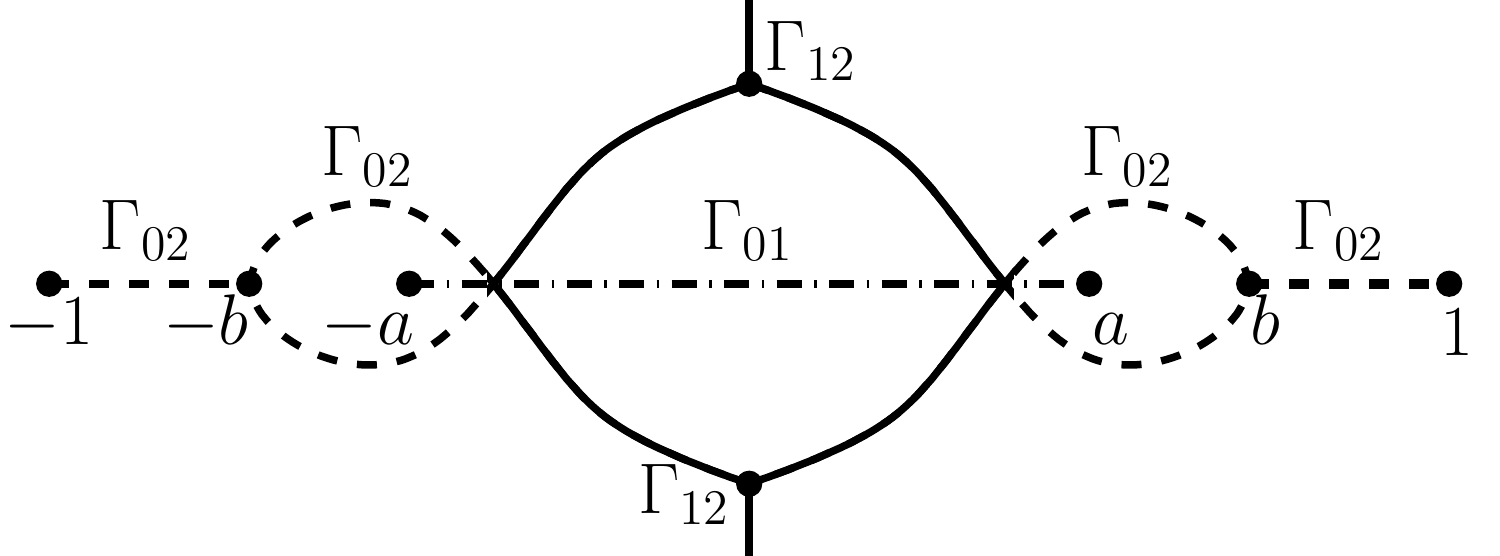}
\caption{\small Case I: $\Gamma_{12}$ - solid lines, $\Gamma_{02}$ - dashed lines, and $\Gamma_{01}$ - dotted dashed line.}
\label{Fig.3A}
\end{figure}

We start with Case I. It will be convenient to consider a slightly different realization of $\RS$, namely, we shall suppose that
\begin{equation}
\label{A.7}
\Delta_0 = [-\mathrm i\infty,-\mathrm ic]\cup[\mathrm ic,\mathrm i\infty].
\end{equation}
Our first goal is to show that $\Gamma$ in this case has the form as shown in Figure~\ref{Fig.3A}.
It follows from our first observation that lines $\Gamma_{ij}$ emanating from a projection of a branch point $\e$ of $\RS$ can be described by
\[
0 = \re\left(\int_e^z(h_i-h_j)(x)\,\mathrm dx\right),
\]
for $z\in\Gamma_{ij}$ locally around $e$. That is, $\Gamma_{ij}$ is a trajectory of a quadratic differential $(h_i-h_j)^2(x)\,\mathrm dx^2$. The local behavior of the trajectories is well known, \cite{Strebel}. This implies that exactly one line of the set $\Gamma$ emanates from the points $\pm 1$ and $\pm a$ (``hard edges") and exactly three lines (with angle $\pi/3$ between them) emanate  from the points $\pm b,\,\pm \mathrm ic$ (``soft edges").

In order to proceed, recall the inequalities \eqref{label} and \eqref{arrange1} as well as the decompositions \eqref{dif1} and \eqref{dif2}. Furthermore, it follows from \eqref{conj-sym} and the fact that $h_0(-z)=-h_0(z)$ that $h_0$ is purely imaginary on $\mathrm i\R$. Analogously, we can conclude that $h_i$, $i\in\{1,2\}$, are purely imaginary on $(-\mathrm ic,0)\cup(0,\mathrm ic)$.
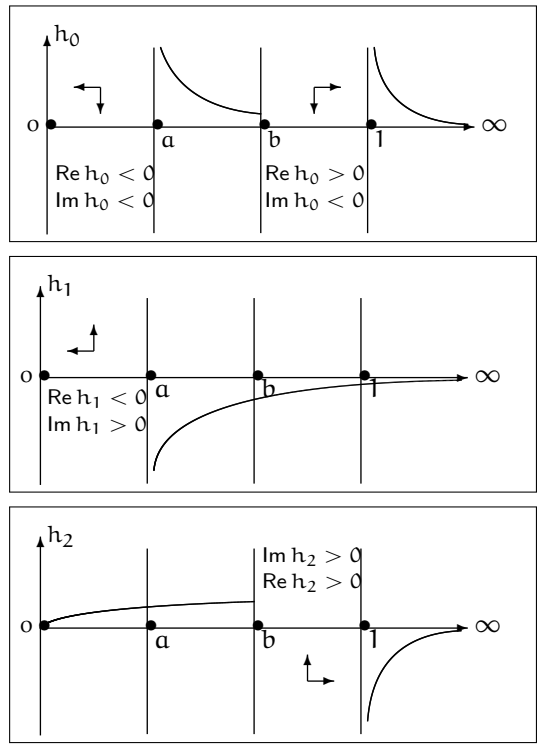
\begin{figure}[!ht]
\parbox{5.5in}{\centering
\framebox{
\unitlength 0.5pt
\linethickness{0.5pt}
\begin{picture}(370,165)(-10,0)
\put(0,80){\vector(1,0){320}}
\put(-15,77){{\small 0}}
\put(325,77){$\infty$}
\put(0,-1){\vector(0,1){150}}
\put(5,145){{\small $ h_{0}$}}
\put(80,0){\line(0,1){140}} \put(160,0){\line(0,1){140}}
\put(240,0){\line(0,1){140}} 
\put(85,65){{\small $a$}}\put(165,65){{\small $b$}} \put(245,65){{\small $1$}}
\put(-2,77){$\bullet$} \put(78,77){$\bullet$}
\put(158,77){$\bullet$} \put(238,77){$\bullet$}
\qbezier(85,140)(100,95)(160,90) \qbezier(245,140)(250,85)(315,82)
\put(40,110){\vector(-1,0){20}} \put(40,110){\vector(0,-1){20}}
\put(200,110){\vector(1,0){20}} \put(200,110){\vector(0,-1){20}}
\put(6,40){{\footnotesize $\re\,h_{0}<0$}}
\put(6,20){{\footnotesize $\im\,h_{0}<0$}}
\put(163,40){{\footnotesize $\re\,h_{0}>0$}} \put(163,20){{\footnotesize $\im\,h_{0}<0$}}
\end{picture}}
\\ \medskip
\framebox{
\unitlength 0.5pt
\linethickness{0.5pt}
\begin{picture}(370,165)(-5,0)
\put(0,80){\vector(1,0){320}}
\put(-15,77){{\small 0}}
\put(325,77){$\infty$}
\put(0,-1){\vector(0,1){150}}
\put(5,145){{\small $ h_{1}$}}
\put(80,0){\line(0,1){140}}
\put(160,0){\line(0,1){140}}
\put(240,0){\line(0,1){140}}
\put(85,65){$a$}
\put(165,65){$b$}
\put(245,65){$1$}
\put(-2,77){$\bullet$}
\put(78,77){$\bullet$}
\put(158,77){$\bullet$}
\put(238,77){$\bullet$}
\qbezier(85,10)(90,75)(315,78)
\put(40,100){\vector(-1,0){20}}
\put(40,100){\vector(0,1){20}}
\put(5,60){{\footnotesize $\re\,h_{1}<0$}}
\put(5,40){{\footnotesize $\im\,h_{1}>0$}}
\end{picture}}\\ \medskip
\framebox{
\unitlength 0.5pt
\linethickness{0.5pt}
\begin{picture}(370,165)(-5,0)
\put(0,80){\vector(1,0){320}}
\put(-15,77){{\small 0}}
\put(325,77){$\infty$}
\put(0,-1){\vector(0,1){150}}
\put(5,145){{\small $h_{2}$}}
\put(80,0){\line(0,1){140}}
\put(160,0){\line(0,1){140}}
\put(240,0){\line(0,1){140}}
\put(85,65){$a$}
\put(165,65){$b$}
\put(245,65){$1$}
\put(-2,77){$\bullet$}
\put(78,77){$\bullet$}
\put(158,77){$\bullet$}
\put(238,77){$\bullet$}
\qbezier(0,80)(10,95)(160,100)
\qbezier(245,10)(250,75)(315,78)
\put(200,40){\vector(1,0){20}}
\put(200,40){\vector(0,1){20}}
\put(166,110){{\footnotesize $\re\,h_{2}>0$}}
\put(166,130){{\footnotesize $\im\,h_{2}>0$}}
\end{picture}}
\caption{\small Behavior of the branches $h$ along the upper bank of $\R_+$.}
\label{Fig.3B-3C}}
\end{figure}

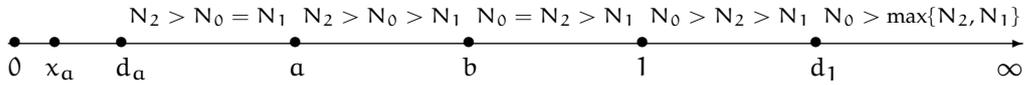
\begin{figure}[ht!]
\centering
\unitlength 1pt
\linethickness{0.5pt}
\begin{picture}(450,20)(0,0)
\put(0,7){\vector(1,0){380}}
\put(0,-5){$0$}
\put(370,-5){$\infty$}
\put(0,5){$\bullet$}
\put(15,5){$\bullet$}
\put(40,5){$\bullet$}
\put(105,5){$\bullet$}
\put(170,5){$\bullet$}
\put(235,5){$\bullet$}
\put(300,5){$\bullet$}
\put(14,-5){$x_{a}$}
\put(40,-5){$d_{a}$}
\put(105,-5){$a$}
\put(170,-5){$b$}
\put(235,-5){$1$}
\put(300,-5){$d_{1}$}
\put(45,15){{\footnotesize$N_{2}>N_{0}=N_{1}$}}
\put(110,15){{\footnotesize$N_{2}>N_{0}>N_{1}$}}
\put(175,15){{\footnotesize$N_{0}=N_{2}>N_{1}$}}
\put(240,15){{\footnotesize$N_{0}>N_{2}>N_{1}$}}
\put(305,15){{\footnotesize$N_{0}>\max\{N_{2},N_{1}\}$}}
\end{picture}
\caption{\small The result of the qualitative analysis of the branches of $N$ along $\mathbb{R}_+$.}
\label{Fig.3D}
\end{figure}

It follows immediately from \eqref{dif1} and \eqref{arrange1}, by choosing the initial point of integration to be $a$ in \eqref{condition}, that
\begin{equation}
\label{A.0}
N_{0}(x)=N_{1}(x),\quad x\in\Delta_1, \quad \mbox{and} \quad N_{0}(x)>N_1(x),\quad x\in(a,b).
\end{equation}
Similarly, by choosing the initial point of integration in \eqref{condition} to be $b$, we deduce from \eqref{dif2} and \eqref{arrange1} that
\begin{equation}
\label{A.3}
N_{0}(x)=N_2(x)>N_1(x),\quad x\in\Delta_{22}, \quad \mbox{and} \quad N_2(x)>N_0(x),\quad x\in(a,b).
\end{equation}
Hence, we get from the symmetry of the surface that
\begin{equation}
\label{A.1}
\Delta_1\subset\Gamma_{01}, \quad \Delta_2\subset\Gamma_{02}, \quad (-b,-a)\cup(a,b)\in\Omega_{201}.
\end{equation}
Moreover we have
\begin{equation}\label{A.2}
N_{2}(a)>N_{0}(a)=N_{1}(a)\,,\quad N_{2}(b)=N_{0}(b)>N_{1}(b)\,.
\end{equation}
To conclude our qualitative analysis based on  Figure~\ref{Fig.3B-3C}, we continue the integration beyond the points $1$ (to the right) and $a$ (to the left). We get that there exists $d_1 >1$ and $d_a>0$ such that
\begin{equation}
\label{A.4}
      (1,d_{1})\in\Omega_{021} \quad \mbox{and} \quad N_{0}(x)=N_{1}(x)>N_{2}(x)\,,\quad x\in(d_{a},a].
\end{equation}
Indeed, if we start integrating from $1$ in the positive direction, then $N_{0}$ increases and $N_{2}$ decreases. Thus, \eqref{A.4} follows from
 \eqref{A.3} and \eqref{A.1} by continuity. We summarize the order of the branches of $N$ along $\R_{+}$ in Figure~\ref{Fig.3D}.

We note that $N_2$ decreases if its argument moves from $d_{a}$ to the left while $N_{0}$ and $N_{1}$ increase, see Figure~\ref{Fig.3B-3C}. Thus, it is possible that there exists $x_{a}\in(0,d_a)$ such that
\begin{equation}
\label{A.5}
\quad N_{0}(x_{a})=N_{1}(x_{a})=N_{2}(x_{a}).
\end{equation}
Moreover, if such a point exists, it is unique. In order to prove the last claim, set
\[
\mathrm d\lambda_1(x) = \frac{f_a(x)}{\sqrt{a^2-x^2}}\,\frac{\mathrm dx}{2\pi},  \quad x\in(-a,a),
\]
where the function $f_a$ is positive and was defined in \eqref{dif2},
\[
\mathrm d\lambda_2(x) = \left\{
\begin{array}{rl}
\displaystyle f_1(x)\sqrt{\frac{x-b}{1-x}}\,\frac{\mathrm dx}{2\pi}, & x\in(b,1), \medskip \\
\displaystyle -f_{-1}(x)\sqrt{\frac{-x-b}{1+x}}\,\frac{\mathrm dx}{2\pi}, & x\in(-1,-b),
\end{array}
\right.
\]
where the function $f_1$ is positive and was defined in \eqref{dif1} while and $f_{-1}$ is negative and was defined in \eqref{dif3}, and finally
\[
\mathrm d\lambda_0(x) = \pm(h_{1+} - h_{1-})(x)\,\frac{\mathrm dx}{2\pi\mathrm i}, \quad \pm x\in[\mathrm ic,\mathrm i\infty).
\]
The measures $\lambda_1$, $\lambda_2$, and $\lambda_0$ are positive (for $\lambda_{12}$ this claim follows from the discussion after \eqref{Delta0cycle} where one needs to recall that we deformed $\Delta_0$ to be as in \eqref{A.7}). It further follows from \eqref{dif1}, \eqref{dif2}, \eqref{dif3}, and Privalov's lemma \cite[Sec.~III.2]{Privalov} that
\[
h_i(z) = \int\frac{\mathrm d\lambda_i(x)}{x-z} - (-1)^i\int\frac{\mathrm d\lambda_0(x)}{x-z}, \quad \text{and} \quad h_0(z) = - \int\frac{\mathrm d(\lambda_1+\lambda_2)(x)}{x-z}.
\]
Since $h_k(z)=2\partial_z N_k(z)$, we deduce that the branches $N_i$, $i\in\{1,2\}$, have the following global representation
\[
N_i(z)=V_i(z)-(-1)^iV_0(z)-c_i,
\]
where $V_k(z)=-\int\log|z-t|\,\mathrm d\lambda_k(t)$, $k\in\{0,1,2\}$. Notice that the symmetry of $\RS$ implies that $\lambda_2$ is an even measure. Therefore, in a complex neighborhood of $[-b,b]$ (which is the gap between two connected components of $\supp(\lambda_{2})$), the potential $V_{2}$ has the form of a saddle such that on the imaginary axis $\mathrm i\R$ it is an even, concave and decreasing function. Hence,
\[
\frac{\partial^{2}V_{2}(x)}{\partial x^{2}}>0,\quad x\in[-b,b], \quad \mbox{and} \quad \frac{\partial^{2}V_{2}(iy)}{\partial y^{2}}<0,\quad y\in\R.
\]
The same is true for $V_0$, that is,
\[
\frac{\partial^{2}V_0(x)}{\partial x^{2}}<0,\quad x\in\R, \quad \mbox{and} \quad \frac{\partial^{2}V_0(iy)}{\partial y^{2}}>0,\quad y\in[-c,c].
\]
Thus, $N_2$ is a convex function on $[-b,b]$. As $2N_{1}=-N_{2}$ on $[-a,a]$, $N_1$ is a concave function on $[-a,a]$. Thus, inside $[-a,a]$, the inequality $N_{1}>N_{2}$ can be true on the connected set only, i.e., on the interval $[-x_{a},x_{a}]$, so \eqref{A.5} and the uniqueness claim are proved and we obtain
\begin{equation}
\label{A.9}
N_0\,=\,N_1\,>\,N_2 \quad \mbox{on}\quad (-x_a,\, x_a).
\end{equation}

Now, consider the trajectory $\Gamma_{12}$ (i.e., the set where $N_{1}=N_{2}$) emanating from the point $\mathrm ic$ into $\{\re\,z>0\}$. It cannot cross the set $\mathrm i\R\cup[1,\infty]\cup[-\infty,-1]$ as it would contradict the maximum principle for harmonic functions. Hence, the only possibility for $\Gamma_{12}$ to cross $[0,1]$ is to cross it at $x_{a}$. Analogous considerations lead to the conclusion that the trajectory of $\Gamma_{02}$ emanating from the point $b$ to the lower half-plane arrives at $x_{a}$ and coincides with the continuation of the considered trajectory $\Gamma_{12}$. Thus, we have that the three subarcs of $\Gamma_{12}$ emanating from the point $\mathrm ic$ (resp. $-\mathrm ic$) terminate at the points $\mathrm i\infty$ (resp. $-\mathrm i\infty$) and $\pm x_{a}$; three subarcs of $\Gamma_{02}$ emanating from the point $b$ (resp. $-b$) terminate at the points $1$ (resp. $-1$) and $\pm x_{a}$; the trajectory $\Gamma_{01}$ joins the points $-a$ and $+a$. That is, we have shown that the set $\Gamma$ has indeed the form as in Figure~\ref{Fig.3A}.

\begin{figure}[ht!]
\includegraphics[scale=.5]{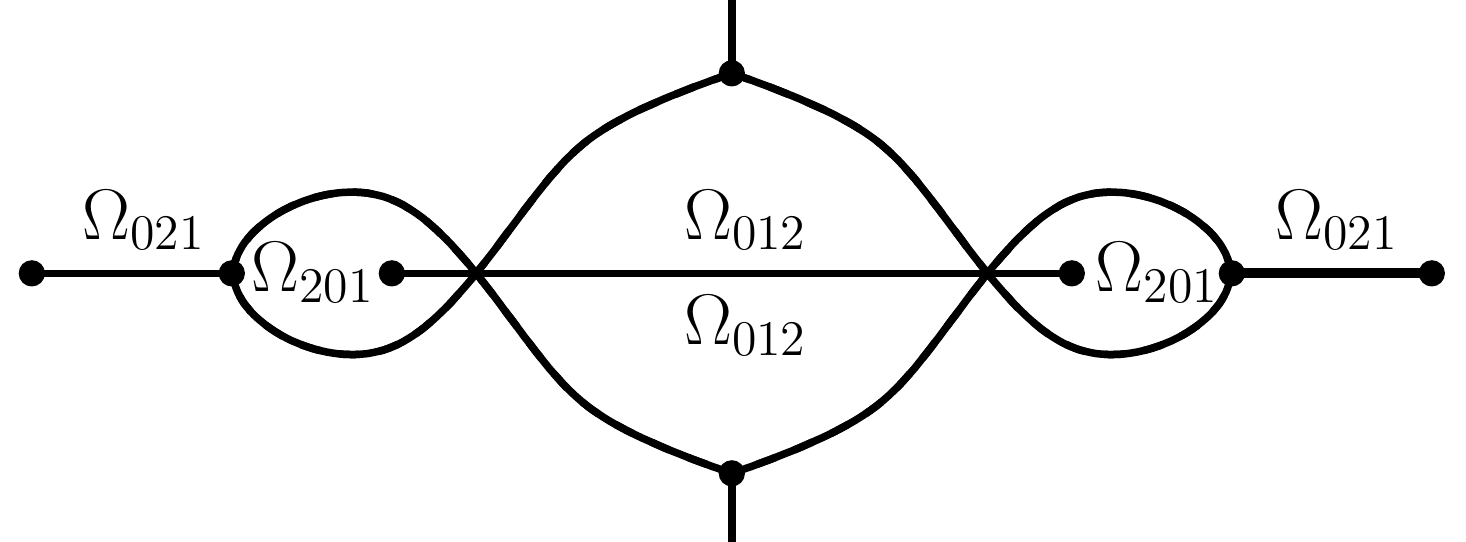}
\caption{\small Case I: The sets $\Omega_{ijk}$ when $\Delta_{0}=[-\mathrm ic,-\mathrm i\infty]\cup[\mathrm ic,\mathrm i\infty]$.}
\label{Fig.3AA}
\end{figure}

The structure of $\Gamma$ and the order of the branches $\{N_k\}_{k=0}^{2}$ on $\R$ (see Figure~\ref{Fig.3D} and \eqref{A.9}) allow us to identify the decomposition of $\mathbb{C}\backslash\Gamma$ into $\cup \, \Omega_{ijk}$ like in Figure~\ref{Fig.3AA}.
It remains to deform the cut \eqref{A.7} to the one in Figure~\ref{Oms}(a) while simultaneously interchanging indices 1 and 2 in the subscripts of $\Omega_{ijk}$ bounded by $\Delta_{0}$ and $[-\mathrm ic,-\mathrm i\infty]\cup[\mathrm ic,\mathrm i\infty]$. This finishes the proof of the Case I of Theorem~\ref{thm:Omegas}.

\begin{figure}[ht!]
\parbox{1.3in}{\centering
\unitlength 1pt
\linethickness{0.5pt}
\begin{picture}(120,140)
\put(20,87){$\bullet$}
\put(20,95){$0$}
\put(40,87){$\bullet$}
\put(40,95){$a^{2}$}
\put(20,90){\line(1,0){20}}
\put(60,87){$\bullet$}
\put(63,78){$b^{2}$}
\put(90,87){$\bullet$}
\put(90,78){$1$}
\put(63,90){\line(1,0){30}}
\thinlines
\qbezier(15,90)(20,150)(60,90)
\qbezier(15,90)(20,30)(60,90)
\put(20,75){$\Omega_{201}$}
\put(50,110){$\Omega_{021}$}
\put(80,100){$w$}
\put(20,17){$\bullet$}
\put(20,8){$0$}
\put(40,17){$\bullet$}
\put(40,8){$a^{2}$}
\put(20,20){\line(1,0){20}}
\put(60,17){$\bullet$}
\put(63,8){$b^{2}$}
\put(90,17){$\bullet$}
\put(90,8){$1$}
\put(63,20){\line(1,0){30}}
\qbezier(13,20)(15,50)(30,20)
\qbezier(13,20)(15,-10)(30,20)
\qbezier(30,20)(45,50)(60,20)
\qbezier(30,20)(45,-10)(60,20)
\put(10,40){$\Omega_{012}$}
\qbezier(9,39)(5,29)(19,25)
\put(50,35){$\Omega_{021}$}
\put(80,30){$w$}
\put(20,130){(I)}
\qbezier(30,70)(35,50)(40,20)
\end{picture}}
\parbox{1.7in}{\centering
\unitlength 1pt
\begin{picture}(110,140)
\put(20,87){$\bullet$}
\put(20,78){$0$}
\put(70,87){$\bullet$}
\put(73,78){$1/2$}
\put(20,90){\line(1,0){100}}
\put(120,87){$\bullet$}
\put(120,78){$1$}
\put(72,40){\line(0,1){100}}
\put(20,120){$\Omega_{012}$}
\put(20,60){$\Omega_{012}$}
\put(100,120){$\Omega_{021}$}
\put(100,60){$\Omega_{021}$}
\put(20,30){II}
\end{picture}}
\parbox{1.7in}{\centering
\unitlength 1pt
\begin{picture}(110,140)
\put(20,87){$\bullet$}
\put(20,78){$0$}
\put(40,87){$\bullet$}
\put(37,78){$b^{2}$}
\put(20,90){\line(1,0){20}}
\put(60,87){$\bullet$}
\put(63,78){$a^{2}$}
\put(90,87){$\bullet$}
\put(90,78){$1$}
\put(63,90){\line(1,0){30}}
\thinlines
\qbezier(42,90)(90,149)(100,90)
\qbezier(42,90)(90,31)(100,90)
\put(75,100){$\Omega_{102}$}
\put(20,17){$\bullet$}
\put(20,8){$0$}
\put(40,17){$\bullet$}
\put(37,8){$b^{2}$}
\put(20,20){\line(1,0){20}}
\put(60,17){$\bullet$}
\put(57,8){$a^{2}$}
\put(90,17){$\bullet$}
\put(90,8){$1$}
\put(63,20){\line(1,0){30}}
\qbezier(40,20)(60,50)(70,20)
\qbezier(40,20)(60,-10)(70,20)
\qbezier(70,20)(94,70)(98,20)
\qbezier(70,20)(94,-20)(98,20)
\put(20,40){$\Omega_{012}$}
\put(77,25){$\Omega_{021}$}
\qbezier(65,45)(60,35)(55,25)
\put(20,120){$\Omega_{012}$}
\put(60,50){$\Omega_{102}$}
\put(100,120){$w$}
\put(100,60){$w$}
\put(120,10){III}
\end{picture}}
\caption{\small Domains $\Omega_{ijk}$ for $\widetilde N_k$. Case (I): $a^{2}<1/2$; Case II: $a^{2}=1/2$; Case III: $a^{2}>1/2$} \label{Fig.3E}
\end{figure}
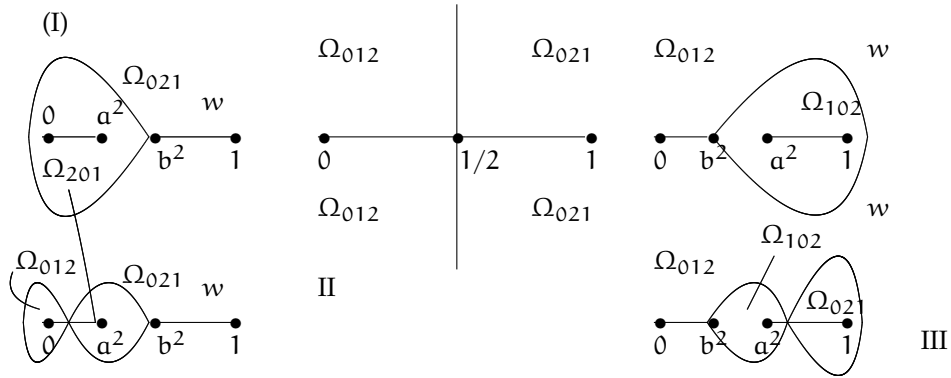

Now we prove Cases II and III of Theorem~\ref{thm:Omegas}. Recall that we put $p=\sqrt{(1+a^2)/3}$ in \eqref{h} and therefore the imaginary branch points $\pm \mathrm ic$ annihilate at infinity. Similar to Case I, we can use \eqref{dif4}, \eqref{dif5}, \eqref{dif7}, and \eqref{dif6}, to show that
\[
N_i(z)=V_i(z)-c_i,\quad V_i(z):=-\int\log|z-t|\,\mathrm d\lambda_i(t),
\]
where $\lambda_i$ is a positive measure supported inside of $\Delta_i$, $i\in\{1,2\}$. Then the choice of the additive constant in \eqref{2.6_1} and the ordering of the branches $N$ along $\R$ give us
\begin{equation}
\label{A.10}
\left\{\begin{array}{l}
(2V_{1}+V_{2})(t)-\gamma_{1}\left\{\begin{array}{ll}
=0, & t\in\supp(\lambda_{1}), \medskip \\
\geqslant 0, & t\in[-a,a],
\end{array}\right.\\
\\
(2V_{2}+V_{1})(t)-\gamma_{2}\left\{\begin{array}{ll}
=0, & t\in\supp(\lambda_{2}) \medskip \\
\geqslant 0, & t\in[-1,-a]\cup[a,1],
\end{array}\right.
\end{array}\right. \qquad \left\{\begin{array}{l}
\gamma_{1}=2c_{1}+c_{2},\medskip \\
\gamma_{2}=2c_{2}+c_{1}.\\
\end{array}\right.
\end{equation}
This type of vector-potential equilibrium problem is well-known in the theory of Hermite-Pad\'e approximants and systems of measures
$\{\lambda_{1},\lambda_{2}\}$ which satisfy such equilibrium conditions are called Angelesco systems, \cite{Ang19, Nik79, GRakh81}. Using the
symmetry of the problem (\ref{A.10}) with respect to the imaginary axis we can transform it to the following one. Set
\begin{eqnarray*}
\widetilde{V}_{1}(w) &:=& \int_{0}^{a^{2}}\log\frac{1}{|w-x|}\,\mathrm d\lambda_{1}(\sqrt{x}) \\
  &=& \frac12\int_{-a}^{a}\log\frac{1}{|z^{2}-t^{2}|}\,\mathrm d\lambda_{1}(t)
= (V_{1}(z)+V_{1}(-z))/2 = V_{1}(z),
\end{eqnarray*}
where we used the change of variables $w=z^2$ and $x=t^2$.
Analogously, we put
\[
\widetilde{V}_{2}(w):=\int_{a^{2}}^{1}\ln\frac{1}{|w-x|}\,\mathrm d\lambda_{2}(\sqrt{x}) = V_{2}(z), \quad w=z^2.
\]
Thus, (with the same constants $\gamma_i$ as in \eqref{A.10}) we
have
\begin{equation}
\label{A.11'}
\left\{\begin{array}{l}
(2\widetilde{V}_{1}+\widetilde{V}_{2})(x)-\gamma_{1}
\left\{\begin{array}{ll}
=0 & x\in\supp(\lambda_{1}(\sqrt{x}))\cap[0,1], \medskip \\
\geqslant 0 & x\in[0,a^{2}], \\
\end{array}\right. \\
\\
(2\widetilde{V}_{2}+\widetilde{V}_{1})(x)-\gamma_{2}
\left\{\begin{array}{ll}
=0 & x\in\supp(\lambda_{2}(\sqrt{x}))\cap[0,1], \medskip \\
\geqslant 0 & x\in[a^{2},1].\\
\end{array}\right. \\
\end{array}\right.
\end{equation}
Therefore,
\begin{equation}
\label{A.12'}
N_i(z)=\widetilde{N}_i(w):=\widetilde{V}_i(w)-c_i, \quad w=z^2, \quad i\in\{1,2\}.
\end{equation}
Angelesco systems \eqref{A.11'} for two touching intervals $[0,a^{2}]$ and $[a^{2},1]$ (for any $a\in(0,1)$) are well studied, see \cite{Kal79, GRakhS97, Ap88,ApKalLysToul09}. In particular, one knows that the regions $\Omega_{ijk}$ constructed for $\widetilde N_k(w)$ decompose $\overline{\mathbb{C}}$ as in Figure~\ref{Fig.3E}.

\noindent Thus, if we plot images of the sets $\Omega_{ijk}$ from Figure~\ref{Fig.3E}
(Cases II and III) under the transformation $z=\sqrt{w}$, then (due to
\eqref{A.12'}) we obtain the corresponding sets from  Figure~\ref{Oms} (b,c,d).
This finishes the proof Theorem~\ref{thm:Omegas}.

\begin{figure}[ht!]
\parbox{2.1in}{
\unitlength 1pt
\linethickness{0.5pt}
\begin{picture}(-200,80)
\put(10,47){$\bullet$}
\put(7,40){$-1$}
\put(50,47){$\bullet$}
\put(40,40){$-b$}
\put(10,49){\line(1,0){40}}
\put(70,47){$\bullet$}
\put(65,40){$-a$}
\put(100,47){$\bullet$}
\put(100,40){$a$}
\put(70,49){\line(1,0){30}}
\put(120,47){$\bullet$}
\put(120,40){$b$}
\put(160,47){$\bullet$}
\put(160,40){$1$}
\put(120,49){\line(1,0){40}}
\thinlines
\qbezier(52,50)(87,115)(122,50)
\qbezier(52,50)(87,-15)(122,50)
\put(80,60){$\Omega_{201}$}
\put(135,60){$\Omega_{021}$}
\end{picture}}
\parbox{3in}{\centering
\unitlength 1pt
\linethickness{0.5pt}
\begin{picture}(200,80)
\put(10,47){$\bullet$}
\put(7,40){$-1$}
\put(50,47){$\bullet$}
\put(40,40){$-b$}
\put(10,49){\line(1,0){40}}
\put(70,47){$\bullet$}
\put(65,40){$-a$}
\put(170,47){$\bullet$}
\put(170,40){$a$}
\put(70,49){\line(1,0){100}}
\put(190,47){$\bullet$}
\put(190,40){$b$}
\put(230,47){$\bullet$}
\put(230,40){$1$}
\put(190,49){\line(1,0){40}}
\thinlines
\qbezier(52,50)(70,92)(90,50)
\qbezier(52,50)(70,8)(90,50)
\qbezier(90,50)(120,99)(150,50)
\qbezier(90,50)(120,1)(150,50)
\qbezier(150,50)(170,92)(192,50)
\qbezier(150,50)(170,8)(192,50)
\put(30,70){$\Omega_{201}$}
\put(110,60){$\Omega_{012}$}
\put(180,20){$\Omega_{201}$}
\put(200,60){$\Omega_{021}$}
\qbezier(50,70)(50,65)(70,55)
\put(71,54){\vector(1,-1){1}}
\qbezier(180,30)(160,40)(165,45)
\put(166,46){\vector(1,1){1}}
\end{picture}}
\caption{\small Sets $\Omega_{ijk}$ for \eqref{h} with $p=\sqrt{(1+a^2)/3}$ and $a < 1/\sqrt2$.}
\label{Fig.3F}
\end{figure}
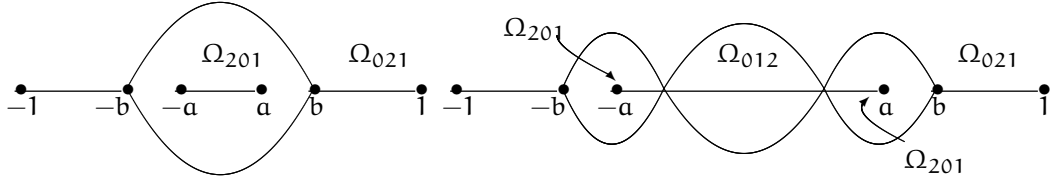

It remains to show that \eqref{sHWRIG} defines single-valued harmonic function on the Riemann surface of $h(\cdot,p^*)$, $p^*=\sqrt{(1+a^2)/3}$, which satisfies \eqref{sdvgads}. This follows from the fact we can obtain $N(\cdot;p^*)$ from $\widetilde N(\cdot;p^*)$ through the transformation $z^2\to w$ and the latter can be computed via \eqref{A.11'} and \eqref{A.12'}. The partition of $\C$ by the domains $\Omega_{ijk}$ for $\widetilde N$ is presented in Figure~\ref{Fig.3E} [Case (I)]. Hence, the corresponding domains for $N$ are distributed as in Figure~\ref{Fig.3F}, which finishes the proof of Theorem~\ref{thm:RS} and justifies the proof of Theorem~\ref{thm:Omegas}.

\section{Nuttall-Szeg\H{o} Functions}
\label{sec:N-S}

This section is devoted to the construction of the Nuttall-Szeg\H{o} functions of Theorems~\ref{thm:N-S} and~\ref{thm:N-Sasymp}. Along the way we prove Propositions~\ref{prop:genus1} and~\ref{prop:genus2}. For brevity, we shall use the following notation:
\[
\RS_{\ualpha} := \RS\setminus\bigcup_{i=1}^g\ualpha_i \quad \text{and} \quad \RS_{\ualpha,\ubeta} := \RS\setminus\bigcup_{i=1}^g\left(\ualpha_i\cup\ubeta_i\right).
\]

\subsection{Abelian Differentials}
\label{ssec:AD}

A holomorphic differential on $\RS$ is a differential of the form
\[
\mathrm{d}\Omega(\z) = f(\z)\,\mathrm{d}z,
\]
where $f$ is a rational function on $\RS$ whose \emph{principal divisor} is of the form
\begin{multline}
\label{(f)}
(f) = \mathcal{D}_f +2\left(\infty^{(0)}+\infty^{(1)}+\infty^{(2)}\right)\\
 -\sum_{\e\in\{\text{branch points of }\RS\}}\big(\{\text{order of branching at }\e\}-1\big) \e,
\end{multline}
for some integral divisor $\mathcal{D}_f$ of order $2g-2$. It is known that such rational functions (integrands of holomorphic differentials) form a subspace of dimension $g$. Hence, there exists precisely one such integrand (up to a multiplicative constant) when $g=1$, and there exist $2$ linearly independent ones when $g=2$. It is further known that if $f_1$ and $f_2$ are distinct integrands of holomorphic differentials, then $\mathcal{D}_{f_1}$ and $\mathcal{D}_{f_2}$ have no element in common. Moreover, any $\z\in\RS$ has a unique complementary point, say $\hat\z$, such that $\mathcal{D}_{f_\z}=\z+\hat\z$ for some holomorphic integrand $f_\z$. When $g=2$, the integral divisors $\mathcal{D}_f$ in \eqref{(f)} are exactly the special divisors mentioned before Proposition~\ref{prop:genus2}. Clearly, if $\mathcal{D}_{f_1}$ and $\mathcal{D}_{f_2}$ are any two special divisors, then $\mathcal{D}_{f_1}-\mathcal{D}_{f_2}$ is principal.

A point $\z$ is a Weierstrass point of a genus 2 surface $\RS$ if there exists a rational function on $\RS$ with a double pole at $\z$ and no other poles. Hence, $\z$ is a Weierstrass point if and only if $\hat\z=\z$. To find Weierstrass points it is enough to find a two-sheeted Riemann surface conformally equivalent to $\RS$. Indeed, Weierstrass points are mapped into Weierstrass points and the Weierstrass points of a two sheeted surface of genus 2 are precisely the branch points. It was shown in \cite[Theorem~1.1]{uApToulVA} that $\RS$ is conformally equivalent to
\begin{equation}
\label{bigZ}
Z^2 = A^\prime(t)^2 -4A(t)\big(3t^2-(1+a^2)+3p^2\big),
\end{equation}
where $A$ is the same polynomial as in \eqref{h} and  the correspondence between the surfaces is given by\footnote{One needs to replace $h$ by $-h$ in \eqref{h} to get the correct correspondence with \cite[Theorem~1.1]{uApToulVA}.}
\[
\left\{
\begin{array}{ll}
t = z+1/h \medskip \\
Z = 2\big(3t^2-(1+a^2)+3p^2\big)/h + A^\prime(t)
\end{array}
\right.
\]
for the uniformizing parameter $t\in\overline\C$. It follows from \eqref{foromega1} and \eqref{foromega2} that
\[
t\big(\infty^{(2)}\big) = \sqrt{\frac{1+a^2}3} \quad \text{and} \quad Z^2\big(\infty^{(2)}\big) = A^\prime\left(t\big(\infty^{(2)}\big)\right)^2.
\]
As $Z\big(\infty^{(2)}\big)\neq0$, it is not a branch point of \eqref{bigZ} and hence $\infty^{(2)}$ is not a Weierstrass point of $\RS$, which is a conclusion we shall need later.

As just explained, $\RS$ has $g$ linearly independent holomorphic differentials. Denote by
\[
\mathrm{d}\vec\Omega := \left\{
\begin{array}{rl}
\left(\mathrm{d}\Omega_1,\mathrm{d}\Omega_2\right)^\mathsf{T}, & g=2, \smallskip \\
\mathrm{d}\Omega_1, & g=1,
\end{array}
\right.
\]
the column vector of such differentials or the unique differential normalized so
\[
\oint_{\ualpha_k}\mathrm{d}\Omega_j = \delta_{jk},
\]
where $\delta_{jk}$ is the usual Kronecker symbol. The differentials $\mathrm{d}\Omega_k$ form a basis for the space of holomorphic differentials on $\RS$. The vector of analytic integrals $\vec\Omega$ defined before \eqref{Omega-jump} is given by
\[
\vec\Omega(\z) := \int_{\e_1}^\z\mathrm{d}\vec\Omega,
\]
where $\e_1$ is the branch point of $\RS$ such that $\pi(\e_1)=1$ (we could have chosen any other point to be the initial point for integration). The significance of this vector lies in Abel's theorem: if $\mathcal{D}_1$ and $\mathcal{D}_2$ are integral divisors, the divisor $\mathcal{D}_1-\mathcal{D}_2$ is principal if and only if the orders of $\mathcal{D}_i$ coincide and
\begin{equation}
\label{Abel}
\vec\Omega(\mathcal{D}_1) \equiv \vec\Omega(\mathcal{D}_2)  \quad \left(\mdp \mathrm{d}\vec\Omega\right),
\end{equation}
where $\vec\Omega(\mathcal{D}):=\sum_i n_i\vec\Omega(\z_i)$ for $\mathcal{D}=\sum_in_i\z_i$.

In what follows we shall also use differentials $\mathrm{d}\Omega_{\z,\w}$ that are the normalized (that is, $\oint_{\ualpha_i}\mathrm{d}\Omega_{\z,\w}=0$, $1\leq i\leq g$) abelian differentials of the third kind with simple poles at $\z$ and $\w$ of respective residues $+1$ and $-1$, and holomorphic otherwise.

\subsection{Cauchy Kernel on $\RS$}

To solve the boundary value problem stated in Theorem~\ref{thm:N-S}, we need a (discontinuous) Cauchy kernel suited for our purposes. Define $\mathrm{d}C_\z$ to be the normalized third kind differential with three simple poles, say  $\z,\z_1,\z_2$, assuming  $\pi^{-1}(z)=\{\z,\z_1,\z_2\}$, which have respective residues $2,-1,-1$. It is quite simple to see that for any fixed $\w$ we have
\[
\mathrm{d}C_\z = 2\mathrm{d}\Omega_{\z,\w}-\mathrm{d}\Omega_{\z_1,\w}-\mathrm{d}\Omega_{\z_2,\w}.
\]

Let $\ugamma$ be a chain on $\RS$ possessing projective involution. For each $\tr\in\ugamma$ which is not a branch point of $\RS$, we shall denote by $\tr^*$ another point on $\ugamma$ having the same canonical projection, i.e., $\pi(\tr)=\pi(\tr^*)$. When $\tr\in\ugamma$ is a branch point of the surface, we simply set $\tr^*=\tr$.

Henceforth, fix $\ugamma$ as above and let $\lambda$ be a H\"older continuous function on $\ugamma$. For simplicity, we shall also assume that the set $\pi(\ugamma)\cap\bigcup_{i=1}^g\pi\left(\ualpha_i\right)$ consists of finitely many points. Define
\[
\Lambda(\z) := \frac1{6\pi\mathrm{i}}\oint_{\ugamma}\lambda\,\mathrm{d}C_\z, \qquad \z\in\RS\setminus\pi^{-1}\bigg(\pi\bigg(\ugamma\cup\bigcup_{i=1}^g\ualpha_i\big)\bigg).
\]
The function $\Lambda$ is holomorphic in the domain of its definition. Moreover, when $\z\to\tr\in\ualpha_k^{\pm}$, one of the points $\z_i$, say $\z_1$, approaches $\tr^*\in\ualpha_k^{\mp}$, while $\z_2$ does not approach any point on $\ugamma\cup\ualpha_k$ (projective involution property). Hence,
\begin{eqnarray*}
\Lambda^+(\tr) - \Lambda^-(\tr) &=& \frac13\left[\frac2{2\pi\mathrm{i}}\oint_{\ugamma}\lambda\left(\mathrm{d}\Omega_{\tr,\w}^+-\mathrm{d}\Omega_{\tr,\w}^-\right) - \frac1{2\pi\mathrm{i}}\oint_{\ugamma}\lambda\left(\mathrm{d}\Omega_{\tr^*,\w}^- - \mathrm{d}\Omega_{\tr^*,\w}^+\right) \right] \\
 &=& \displaystyle -\oint_{\ugamma}\lambda\,\mathrm{d}\Omega_k,
\end{eqnarray*}
since, according to \cite[Eq. (1.9)]{Zver71}, one has
\[
\mathrm{d}\Omega_{\z,\w}^+-\mathrm{d}\Omega_{\z,\w}^- = -2\pi\mathrm{i}\mathrm{d}\Omega_k, \quad \z\in\ualpha_k.
\]
On the other hand, if $\z\to\tilde\tr$,  $\z_1\to\tr\in\ualpha_k^\pm$, and $\z_2\to\tr^*\in\ualpha_k^\mp$, then
\[
\Lambda^+\big(\tilde\tr\big) - \Lambda^-\big(\tilde\tr\big) = \frac13\left[-\frac1{2\pi\mathrm{i}}\oint_{\ugamma}\lambda\left(\mathrm{d}\Omega_{\tr,\w}^+ - \mathrm{d}\Omega_{\tr,\w}^-\right) - \frac1{2\pi\mathrm{i}}\oint_{\ugamma}\lambda\left(\mathrm{d}\Omega_{\tr^*,\w}^- - \mathrm{d}\Omega_{\tr^*,\w}^+\right) \right] = 0.
\]
Furthermore, if $\z\to\tr\in\ugamma^{\pm}\setminus\bigcup_{i=1}^g\ualpha_i$,  $\z_1\to\tr^*\in\ugamma^{\mp}$, while $\z_2$ does not approach any point on $\ugamma\cup\bigcup_{i=1}^g\ualpha_i$, then
\[
\Lambda^+(\tr) - \Lambda^-(\tr) = \frac{2\lambda(\tr) + \lambda(\tr^*)}3,
\]
according to \cite[Eq. (2.8)]{Zver71}. Finally, if $\z\to\tilde\tr$,  $\z_1\to\tr\in\ugamma^\pm$, and $\z_2\to\tr^*\in\ugamma^\mp$, then
\[
\Lambda^+\big(\tilde\tr\big) - \Lambda^-\big(\tilde\tr\big) = \frac{\lambda(\tr^*) - \lambda(\tr)}3.
\]
Thus, if we additionally require that $\lambda(\tr)=\lambda(\tr^*)$, then $\Lambda$ is a holomorphic function in $\RS_\ualpha\setminus\ugamma$ such that
\[
\left(\Lambda^+ - \Lambda^-\right)(\tr) = \left\{
\begin{array}{rl}
\displaystyle -\oint_{\ugamma}\lambda\,\mathrm{d}\Omega_k, & \tr\in\ualpha_k\setminus\big(\ugamma\cup\ualpha_{3-k}\big), \medskip \\
\lambda(\tr), & \tr\in\ugamma\setminus\bigcup_{i=1}^g\ualpha_i.
\end{array}
\right.
\]
It also can readily be verified that
\begin{equation}
\label{Lambda-norm}
\Lambda(\z) + \Lambda(\z_1) + \Lambda (\z_2) \equiv 0 \quad \text{on} \quad \RS ,
\end{equation}
since the cycles $\ugamma$ and $\ualpha_i$ possess projective involutions.

\subsection{Logarithmic Discontinuities}

It is known that the continuity of $\Lambda^\pm$, in fact, H\"older continuity, depends on the H\"older continuity of $\lambda$ only locally. That is, if $\lambda$ is H\"older continuous on some open subarc of $\ugamma$, so are the traces $\Lambda^\pm$ on this subarc irrespective of the smoothness of $\lambda$ on the remaining part of $\ugamma$. Of particular importance for us are the logarithmic discontinuities at the branch points of $\RS$.

More precisely, let $\e$ be a branch point of $\RS$ and $\mathfrak{U}$ be its circular neighborhood. Given a cycle on $\RS$  that possesses a projective involution and passes through $\e$, denote by $\ugamma$ its part that belongs to $\mathfrak{U}$. In this case $\pi(\ugamma)$ is a Jordan arc ending at $e=\pi(\e)$. Recall \cite[Sec. I.8.6]{Gakhov} that
\begin{equation}
\label{Gakhov}
\frac1{2\pi\mathrm{i}}\int_{\pi(\ugamma)}\frac{\alpha\log(t-e)}{(t-e)^{\nu+}}\frac{\mathrm{d}t}{t-z} = \pm \frac{\alpha e^{\pm\nu\pi\mathrm{i}}}{2\mathrm{i}\sin(\nu\pi)}\frac{\log(z-e)}{(z-e)^\nu} + \mathcal{O}\left((z-e)^{-\nu}\right)
\end{equation}
as $z\to e$, where $\alpha\log(z-e)$, $\alpha\in\R$, and $(z-e)^\nu$, $\nu\neq0$, are holomorphic outside of $\pi(\ugamma)$, $\log(t-e)$ is any continuous determination of the logarithm on $\pi(\ugamma)$, $(t-e)^{\nu+}$ is the trace of $(z-e)^\nu$ on the positive side of $\pi(\ugamma)$, and the choice of the sign in $\pm$ is determined by the orientation of $\pi(\ugamma)$: $+$ if $\pi(\ugamma)$ is oriented away from $e$ and $-$ if it is oriented towards $e$.

Assuming that $\RS$ has double branching at $\e$, define $u(\z):=\pm(z-e)^{1/2}$, where the choice of the sign depends on whether $\z$ lies to the left ($+$) or to the right ($-$) of $\ugamma$. Set
\[
\mathrm{d}\mathcal{C}_\z(\tr) := \frac12\left(1+3\frac{u(\z)}{u(\tr)}\right)\frac{\mathrm{d}t}{t-z}.
\]
The differential $\mathrm{d}\mathcal{C}_\z$ is holomorphic in $\mathfrak{U}\setminus\{\z,\z^*\}$ with residues $+2$ at $\z$ and $-1$ at $\z^*$, where $\z^*$ is the unique point in $\mathfrak{U}$ which is different from $\z$ and has the same natural projection. Clearly, the difference $\mathrm{d}C_\z-\mathrm{d}\mathcal{C}_\z$ is a holomorphic differential in $\mathfrak{U}$. Therefore, if $\lambda(\tr)$ is of the form $\alpha\log(t-e) + \big\{\text{H\"older continuous part}\big\}$ satisfying $\lambda(\tr)=\lambda(\tr^*)$, then
\begin{eqnarray}
\Lambda(\z) & = & \frac1{6\pi\mathrm{i}} \int_\ugamma\alpha\log(t-e)\,\mathrm{d}\mathcal{C}_\z(\tr) + \mathcal{O}(1) \nonumber \\
\label{Lambda-Log}
& = & \frac1{2\pi\mathrm{i}}\int_{\pi(\ugamma)}\alpha\log(t-e)\frac{u(\z)}{\sqrt{t-e}^+} \frac{\mathrm{d}t}{t-z}  + \mathcal{O}(1) = \pm \frac\alpha2\log(z-e) + \mathcal{O}(1)
\end{eqnarray}
as $\z\to\e\in\ugamma^\pm$ by \eqref{Gakhov} applied with $\nu=1/2$ (observe that the second equality above is valid independently of the orientation of $\pi(\ugamma)$).

Assume now that $\RS$ has triple branching at $\e$.  Since $\ugamma$ possesses projective involution there is exactly one component of $\mathfrak{U}\setminus\ugamma$ that lies schlicht over $\pi(\mathfrak{U}\setminus\ugamma)$. Orient $\ugamma$ so that this component lies to the right when $\ugamma$ is traversed in the positive direction (it is bordered by $\ugamma^-$). Define
\[
\mathrm{d}\mathcal{C}_\z(\tr) := \left(\frac{u(\z)}{u(\tr)}+\frac{u^2(\z)}{u^2(\tr)}\right)\frac{\mathrm{d}t}{t-z},
\]
where $u(\z)$ is equal to the lift of a fixed determination of $(z-e)^{1/3}$ to the ``schlicht'' component of $\mathfrak{U}$ and then is continued holomorphically to the whole $\mathfrak{U}$. Since the values of $u$ at the points with the same natural projection  differ by $e^{\pm2\pi\mathrm{i}/3}$, one has that $\mathrm{d}\mathcal{C}_\z$ is a holomorphic differential in $\mathfrak{U}\setminus\pi^{-1}(z)$, with 3 simple poles at the elements of $\pi^{-1}(z)$ having residues $+2$ at $\z$ and $-1$ at the two elements of $\pi^{-1}(z)\setminus\{\z\}$. Again, $\mathrm{d}C_\z-\mathrm{d}\mathcal{C}_\z$ is a holomorphic differential in $\mathfrak{U}$ and therefore
\begin{eqnarray}
\Lambda(\z) & = & \frac1{6\pi\mathrm{i}} \int_\ugamma\alpha\log(t-e)\,\mathrm{d}\mathcal{C}_\z(\tr) + \mathcal{O}(1) \nonumber \\
& = & \frac{\sqrt3\mathrm{i}}{6\pi\mathrm{i}} \int_{\pi(\ugamma)}\alpha\log(t-e)\left(\pm e^{\pm\pi\mathrm{i}/3} \frac{u(\z)}{(t-e)^{1/3+}} \pm e^{\pm2\pi\mathrm{i}/3}\frac{u^2(\z)}{(t-e)^{2/3+}} \right) \frac{\mathrm{d}t}{t-z}, \nonumber
\end{eqnarray}
where the sign $+$ is used when $\pi(\ugamma)$ is oriented towards $e$ and the sign $-$ is used if it is oriented away from $e$, and we used the fact that
\[
\int_\ugamma\frac{\mathrm{d}t}{u(\tr)} = \int_{\pi(\ugamma)}\frac{\mathrm{d}t}{(t-e)^{1/3-}} -  \int_{\pi(\ugamma)}\frac{\mathrm{d}t}{(t-e)^{1/3+}} =  \pm\sqrt3\mathrm{i}e^{\pm\pi\mathrm{i}/3}\int_{\pi(\ugamma)}\frac{\mathrm{d}t}{(t-e)^{1/3-}}
\]
with the same choice of signs. Then \eqref{Gakhov} applied with $\nu=1/3$ and $\nu=2/3$ gives that
\begin{eqnarray}
\Lambda(\z) & = & -\frac\alpha3\log(z-e)\left(\frac{u(\z)}{(z-e)^{1/3}} + \frac{u^2(\z)}{(z-e)^{2/3}} \right) + \mathcal{O}(1) \nonumber \\
\label{Lambda-Log1}
& = & \mathcal{O}(1) + \frac\alpha3\log(z-e)\left\{ \begin{array}{rl} 1,  & \z\to\e\in\ugamma^+, \smallskip \\ -2, & \z\to\e\in\ugamma^-, \end{array} \right.
\end{eqnarray}
where the determination of $(z-e)^{1/3}$ is precisely the one of $u(\z)$ within the component of $\mathfrak{U}$ that lies to the right of $\ugamma$.

\subsection{Szeg\H{o} Function $S_\rho$}

Let $\boldsymbol\Delta$ be the chain introduced in Theorem~\ref{thm:N-S}. We assume that the cycles forming $\boldsymbol\Delta$ and separating $\RS^{(0)}$ and $\RS^{(i)}$, $i\in\{1,2\}$, are oriented so that $\RS^{(0)}$ lies on the left (positive) side as each of the cycles is traversed in the positive direction, while the cycle separating $\RS^{(1)}$ and $\RS^{(2)}$ is oriented so that $\RS^{(1)}$ remains on the left when it is traversed in the positive direction.

To define a density whose Cauchy transform on $\boldsymbol\Delta$ we shall need, let us introduce auxiliary ``trigonometric'' weights. Namely,
\begin{equation}
\label{trig1}
\left\{
\begin{array}{rcl}
w_1(z) &:=& \sqrt{z^2-a^2}, \medskip \\
w_2(z) &:=& \sqrt{(z^2-1)(z^2-b^2)(z^2+c^2)},
\end{array}
\right.
\end{equation}
in Case I with branch cuts along $\Delta_1$, $\Delta_2\cup\Delta_0$, respectively (all are positive for positive reals large enough), and
\begin{equation}
\label{trig2}
\left\{
\begin{array}{rcl}
w_1(z) &:=& \sqrt{z^2-b^2}, \medskip \\
w_2(z) &:=& \sqrt{(z^2-1)(z^2-a^2)},
\end{array}
\right.
\end{equation}
in Cases II and III with branch cuts along $\Delta_1$ and $\Delta_2$, respectively (all are positive for positive reals large enough). We set
\begin{equation}
\label{lambda-rho}
\lambda_\rho(\tr) :=
\left\{
\begin{array}{ll}
-\log\big(\rho_1w_1^+\big)(t), & t\in\Delta_1^\circ, \medskip \\
-\log\big(-\rho_2w_2^+\big)(t), & t\in\Delta_{21}^\circ, \medskip \\
-\log\big(\rho_2w_2^+\big)(t), & t\in\Delta_{22}^\circ, \medskip \\
-\log\big(\rho_2w_2^+/\rho_1w_1\big)(t), & t\in\Delta_0^\circ,
\end{array}
\right.
\quad \tr\in\boldsymbol\Delta,
\end{equation}
where we choose continuous branches of the logarithms, which is possible as the functions $\rho_i$ and the trigonometric weights are non-vanishing except at the endpoints. Put
\begin{equation}
\label{Srho}
S_\rho (\z) := \exp\big\{\Lambda_\rho(\z)\big\}, \quad \vec c_\rho := -\frac1{2\pi\mathrm{i}}\oint_{\boldsymbol\Delta}\lambda_\rho\,\mathrm{d}\vec\Omega.
\end{equation}
Then $S_\rho$ is a holomorphic and non-vanishing function on $\RS_{\ualpha}\setminus\boldsymbol\Delta$ with continuous traces except at the branch points that satisfy
\begin{equation}
\label{Srho-jumps}
S_\rho^+ = S_\rho^- \left\{
\begin{array}{lll}
\displaystyle \exp\left\{2\pi\mathrm{i}\big(\vec c_\rho\big)_k\right\} & \text{on} & \ualpha_k \medskip \\
\displaystyle \exp\left\{\lambda_\rho\right\} & \text{on} & \boldsymbol\Delta.
\end{array}
\right.
\end{equation}
Moreover, $S_\rho$ satisfies \eqref{normal} with $\Phi$ replaced by $S_\rho$ as follows from \eqref{Lambda-norm}, and, excluding the case $a=b=1/\sqrt2$, one has
\begin{equation}
\label{Srho-BU}
\left\{
\begin{array}{lll}
\big|S_\rho^{(0)}(z)\big| \sim |z-e|^{-1/4} & \text{as} & z\to e\in\{\pm1,\pm a, \pm b\}, \medskip \\
\big|S_\rho^{(1)}(z)\big| \sim |z-e|^{1/4} & \text{as} & z\to e\in\{\pm e_1\}, \medskip \\
\big|S_\rho^{(1)}(z)\big| \sim |z-e|^{-1/4} & \text{as} & z\to e\in\left\{\pm\mathrm{i}c\right\}, \medskip \\
\big|S_\rho^{(2)}(z)\big| \sim |z-e|^{1/4} & \text{as} & z\to e\in\{\pm1,\pm e_2,\pm\mathrm{i}c\},
\end{array}
\right.
\end{equation}
by \eqref{Lambda-Log}, where $e_1=a$ and $e_2=b$ in Case I while $e_1=b$ and $e_2=a$ in Cases II and III, and the points $\pm\mathrm{i}c$ appear only in Case I. When $e_1=e_2=1/\sqrt2$, we need to use \eqref{Lambda-Log1} instead of \eqref{Lambda-Log}. Now there are two cycles that pass through a branch point: one that separates $\RS^{(0)}$ and $\RS^{(1)}$ ($\RS^{(1)}$ lies to the right of this cycle), and one that separates $\RS^{(0)}$ and $\RS^{(2)}$ ($\RS^{(2)}$ lies to the right of this cycle).  Combing the contributions from the integrals over both cycles we get that
\begin{equation}
\label{Srho-BU1}
\left\{
\begin{array}{lll}
\big|S_\rho^{(0)}(z)\big| \sim |z-e|^{-1/3} & \text{as} & z\to e\in\big\{\pm1/\sqrt2\big\}, \medskip \\
\big|S_\rho^{(i)}(z)\big| \sim |z-e|^{1/6} & \text{as} & z\to e\in\big\{\pm1/\sqrt2\big\}, \quad i\in\{1,2\}.
\end{array}
\right.
\end{equation}

\subsection{Riemann Theta Function}

The theta function associated with $\mathbf B$ is an entire transcendental function of $g$ complex variables defined by
\[
\theta\left(\vec u\right) := \sum_{\vec n\in\Z^g}\exp\bigg\{\pi\mathrm{i}\vec n^\mathsf{T}\mathbf B\vec n + 2\pi\mathrm{i}\vec n^\mathsf{T}\vec u\bigg\}, \quad \vec u\in\C^g.
\]
As shown by Riemann, the symmetry of $\mathbf B$ and positive definiteness of its imaginary part ensures the convergence of the series for any $\vec u$. It is known that
\begin{equation}
\label{theta-zeros}
\theta\left(\vec u\right)=0 \quad \Leftrightarrow \quad \vec u\equiv \left\{
\begin{array}{rl}
\vec\Omega\left(\w\right) + \vec K, & g=2, \smallskip \\
\vec K, & g=1,
\end{array}
\right. \quad \left(\mdp \mathrm{d}\vec\Omega\right)
\end{equation}
for some $\w\in\RS$, where the equivalence $\equiv$ was defined in \eqref{equivalence} and $\vec K$ is the vector of Riemann constants defined by
\[
\left\{
\begin{array}{lcll}
\big(\vec K\big)_k &:=&\left([\mathbf B]_{kk}-1\right)/2 -\oint_{\ualpha_{3-k}}\Omega_k^-\mathrm{d}\Omega_{3-k}, & g=2, \smallskip \\
\vec K &:=& \left(\mathbf B-1\right)/2, & g=1.
\end{array}
\right.
\]
Let $\mathcal{D}$ be an integral divisor of order $g$. Then $\theta\big( \vec\Omega(\z) - \vec\Omega(\mathcal{D})-\vec K\big)$ is a multi-valued holomorphic function on $\RS$ if $\mathcal{D}$ is not special and is identically zero otherwise. Indeed, in the later case $\Omega(\mathcal{D})=\Omega(\z+\hat\z)+\vec j_\z+\mathbf{B}\vec m_\z$, $\vec j_\z,\vec m_\z\in \Z^g$ and therefore
\[
\theta\left( \vec\Omega(\z) - \vec\Omega(\mathcal{D})-\vec K\right) = \theta\left( -\vec\Omega(\hat\z) - \vec K\right) e^{-\pi \mathrm{i}\vec m_\z^\mathsf{T}\mathbf B \vec m_\z - 2\pi \mathrm{i}\vec m_\z^\mathsf{T}\left(\vec\Omega(\hat\z)+\vec K\right)} = 0
\]
by \eqref{theta-zeros}, where we used the fact $\theta(-\vec u)=\theta(\vec u)$ and the periodicity property of theta functions:
\begin{equation}
\label{theta-periods}
\theta\left(\vec u + \vec j + \mathbf B\vec m\right) = \exp\bigg\{-\pi \mathrm{i}\vec m^\mathsf{T}\mathbf B \vec m - 2\pi \mathrm{i}\vec m^\mathsf{T}\vec u\bigg\}\theta\big(\vec u\big), \quad \vec j,\vec m\in\Z^g.
\end{equation}

Set $\vec v_k := \vec\Omega\big(\infty^{(k)} - \infty^{(0)}\big)$, $k\in\{0,1,2\}$, and define
\begin{equation}
\label{thetan}
\Theta_{n,k}(\z) := \frac{\theta\left( \vec\Omega(\z) - \vec\Omega\big(g\infty^{(2)}\big) - \vec K -\vec v_k - \vec c_\rho - \vec x_n - \mathbf{B}\vec y_n \right)}{\theta\left( \vec\Omega(\z) - \vec\Omega\big(g\infty^{(2)}\big) - \vec K \right)}
\end{equation}
(recall that $2\infty^{(2)}$ is not a special divisor, see \eqref{bigZ} and after, and therefore $\Theta_{n,k}$ is always well defined), where $\vec x_n,\vec y_n\in[0,1)^g$ are such that
\begin{equation}
\label{xnyn}
\vec x_n + \mathbf{B}\vec y_n \equiv n\left(\vec \omega + \mathbf{B}\vec\tau\right)  \quad \left(\mdp \mathrm{d}\vec\Omega\right).
\end{equation}
If the numerator is not identically zero, $\Theta_{n,k}$ is a multiplicatively multi-valued meromorphic function on $\RS$ with $g$ simple poles at $\infty^{(2)}$, $g$ zeros that we shall describe by an integral divisor $\mathcal{D}_{n,k}$, and otherwise non-vanishing and finite (it is also possible that zeros could cancel poles). Moreover, it follows from \eqref{Omega-jump} and \eqref{theta-periods} that $\Theta_{n,k}$ is meromorphic and single-valued in $\RS_\ualpha$ and
\begin{equation}
\label{jumpThetan}
\Theta_{n,k}^+  =  \Theta_{n,k}^- \exp\left\{-2\pi \mathrm{i}\left(\vec v_k + \vec c_\rho + \vec x_n + \mathbf{B}\vec y_n\right)_i\right\} \qquad \text{on} \quad \ualpha_i.
\end{equation}

Together with the functions $\Theta_{n,k}$ we shall need two more auxiliary theta functions. We set
\[
\Theta_k(\z) =  \frac{\theta\left( \vec\Omega(\z) - \vec\Omega\big(\infty^{(0)}+(g-1)\w\big) - \vec K \right)}{\theta\left( \vec\Omega(\z) - \vec\Omega\big(\infty^{(k)}+(g-1)\w\big) - \vec K \right)}, \quad k\in\{0,1,2\},
\]
for any fixed $\w$ such that $\w\neq\hat\infty^{(k)}$ for all $k\in\{0,1,2\}$. Clearly, $\Theta_0\equiv1$ and $\Theta_i$, $i\in\{1,2\}$, is a meromorphic function in $\RS_\ualpha$, with a simple zero at $\infty^{(0)}$, a simple pole at $\infty^{(i)}$, and otherwise non-vanishing and finite. Moreover, its traces on $\ualpha_k$ satisfy
\begin{equation}
\label{jumpThetai}
\Theta_i^+  =  \Theta_i^- \exp\left\{2\pi \mathrm{i}\left(\vec v_i\right)_k\right\} \qquad \text{on} \quad \ualpha_k.
\end{equation}

\subsection{Jacobi Inversion Problem}

It follows from \eqref{theta-zeros} and \eqref{xnyn} that the divisor $\mathcal{D}_{n,k}$ is a solution of the Jacobi inversion problem
\begin{equation}
\label{JIP}
\vec\Omega\big(\mathcal{D}_{n,k}\big) \equiv \vec\Omega\big( g\infty^{(2)} \big) + \vec v_k + \vec c_\rho + n\big(\vec\omega+\mathbf B\vec\tau\big) \quad \left(\mdp \mathrm{d}\vec\Omega\right).
\end{equation}
Clearly, $\mathcal{D}_{n,0}$ are precisely the divisors $\mathcal{D}_n$ defined in \eqref{one-jip}.

As we already mentioned after \eqref{one-jip}, any Jacobi inversion problem is uniquely solvable on surfaces of genus 1.  Thus, $\mathcal{D}_{n,k}$ are well defined in this case. Since
\begin{equation}
\label{2div}
\vec\Omega\big(\mathcal{D}_{n,k}\big) \equiv \vec\Omega\big(\mathcal{D}_{n-1,k}\big) + \vec\omega+\mathbf B\vec\tau \quad \left(\mdp \mathrm{d}\vec\Omega\right),
\end{equation}
Proposition~\ref{prop:periods} and the unique solvability property imply that
\[
\mathcal{D}_{2m,k} = \mathcal{D}_{0,k} \neq \mathcal{D}_{1,k} = \mathcal{D}_{2m+1,k}
\]
for all $m\geq0$. This finishes the proof of Proposition~\ref{prop:genus1}. Observe also that
\begin{equation}
\label{lp}
\vec\Omega\big(\mathcal{D}_{n,k}\big) \equiv \vec\Omega\big(\mathcal{D}_{n,j}-\infty^{(j)} + \infty^{(k)}\big) \quad \left(\mdp \mathrm{d}\vec\Omega\right).
\end{equation}
The latter equivalence together with the unique solvability property imply that
\[
\mathcal{D}_{n,j} = \infty^{(j)} \quad \text{for some} \quad j\in\{0,1,2\} \quad  \Rightarrow \quad \mathcal{D}_{n,k} = \infty^{(k)} \quad \text{for all} \quad k\in\{0,1,2\}.
\]
Hence, if $\mathcal{D}_{n,0}\neq\infty^{(0)}$, then $\mathcal{D}_{n,i}\neq\infty^{(i)}$, $i\in\{1,2\}$, which will be important in what follows.

As to Proposition~\ref{prop:genus2}, observe that its first claim follows from \eqref{2div} and Proposition~\ref{prop:periods} as $\vec\Omega$ assumes equivalent values at special divisors by Abel's theorem \eqref{Abel}.  To prove the second claim, assume to the contrary that all the limit points of $\big\{\mathcal{D}_{n,0}\big\}_{n\in\N}$ are either special or of the form $\infty^{(0)}+\w$ for some $\w\in\RS$. Using compactness, we always can select a subsequence $\N_1$ such that
\[
\mathcal{D}_{n+j,0} \to \mathcal{D}^j \quad \text{as} \quad \N_1\ni n\to\infty
\]
for some divisors $\mathcal{D}^j$ simultaneously for all $j\in\{0,1,2,3\}$. According to our assumption,
one has
\begin{equation}
\label{limit-step1}
\vec\Omega\big(\mathcal{D}^j\big) \equiv \vec\Omega\big(\infty^{(0)}+\w_j\big),
\end{equation}
where $\w_j=\hat\infty^{(0)}$ when $\mathcal{D}^j$ is special. Observe also that continuity of $\vec \Omega$ and \eqref{2div} imply that
\begin{equation}
\label{limit-step2}
\vec\Omega\big(\mathcal{D}^{j+1}\big) \equiv \vec\Omega\big(\mathcal{D}^j\big) + \vec\omega + \mathbf{B}\vec \tau \quad \left(\mdp \mathrm{d}\vec\Omega\right)
\end{equation}
for $j\in\{0,1,2\}$. By combining \eqref{limit-step2} for different $j$ and then recalling \eqref{limit-step1},  we deduce that
\[
\vec\Omega\big(2\w_j\big) \equiv \vec\Omega\big(\w_{j+1}+\w_{j-1}\big) \quad \left(\mdp \mathrm{d}\vec\Omega\right),
\]
 $j\in\{1,2\}$. The latter equivalence and Abel's theorem imply that $2\w_1$ and $\w_3+\w_1$ are special divisors, and therefore necessarily $\w_3=\w_1$. The above conclusion and \eqref{limit-step2} imply that
 \[
 2\vec\omega + 2\mathbf{B}\vec\tau \equiv \vec 0  \quad \left(\mdp \mathrm{d}\vec\Omega\right),
 \]
 which is impossible since $2\vec\omega=(2\omega,4(1-\omega))$ and $2\omega\in(1,2)$ by Proposition~\ref{prop:periods}. This contradiction proves the second claim of Proposition~\ref{prop:genus2}.

Finally, observe that the subsequence $\N_*$ is such that no limit point of $\big\{\mathcal{D}_{n,i}\big\}_{n\in\N_*}$ is special or of the form $\infty^{(i)}+\w$ for $i\in\{1,2\}$. Indeed, otherwise there would exist a limit point $\mathcal{D}$ of one those sequences such that
\[
\vec\Omega\big(\mathcal{D}\big) \equiv \vec\Omega\big(\infty^{(i)}+\w\big),
\]
where $\w=\hat\infty^{(i)}$ if $\mathcal{D}$ is special. Choose a limit point $\mathcal{D}_*$ of $\big\{\mathcal{D}_{n,0}\}$ along the same subsequence. The continuity of $\vec\Omega$ and \eqref{lp} then would imply that
\[
\vec\Omega\big(\mathcal{D}_*\big) \equiv \vec\Omega\big(\mathcal{D}-\infty^{(i)} + \infty^{(0)}\big) \equiv  \vec\Omega\big(\infty^{(0)} + \w \big) \quad \left(\mdp \mathrm{d}\vec\Omega\right),
\]
which means that either $\mathcal{D}_*$ is special or contains $\infty^{(0)}$. As both options are impossible by the choice of $\N_*$, such a limit point $\mathcal{D}$ does not exist.

\subsection{Nuttall-Szeg\H{o} Functions}

Finally, we are ready to prove Theorem~\ref{thm:N-S}. Let $\vec y\in\R^g$. Then it follows from \eqref{Omega-jump} that
\[
S_{\vec y}(\z) := \exp\big\{-2\pi\mathrm{i}\vec y^\mathsf{T}\vec\Omega(\z)\big\}, \quad \z\in\RS_{\ualpha,\ubeta},
\]
is a holomorphic function on $\RS_{\ualpha,\ubeta}$ with continuous traces that satisfy
\begin{equation}
\label{jumpSy}
S_{\vec y}^+ = S_{\vec y}^-\left\{
\begin{array}{lll}
\displaystyle \exp\left\{2\pi\mathrm{i}\big(\mathbf{B}\vec y\big)_i\right\} & \text{on} & \ualpha_i \medskip \\
\displaystyle \exp\left\{-2\pi\mathrm{i}\big(\vec y\big)_i\right\} & \text{on} & \ubeta_i.
\end{array}
\right.
\end{equation}
Moreover, \eqref{Lambda-norm} implies that \eqref{normal} holds with $\Phi$ replaced by $S_{\vec y}$. Finally, set
\[
W\big(z^{(0)}\big) \equiv 1 \quad \text{and} \quad W\big(z^{(i)}\big) = w_i(z), \quad i\in\{1,2\},
\]
where the functions $w_1$ and $w_2$ are defined by \eqref{trig1} or \eqref{trig2}, depending on the considered case. The Nuttall-Szeg\H{o} functions we are looking for are defined by
\begin{equation}
\label{Psink}
\Psi_{n,k} := \Phi^nS_\rho S_{\vec y_n}\Theta_{n,k}\Theta_kW^{-1},
\end{equation}
where $\vec y_n$ was defined in \eqref{xnyn} and we assume $n$ is such that the corresponding Jacobi inversion problem \eqref{JIP} is uniquely solvable. Indeed, it can readily be checked, using \eqref{Phi-jumps}, \eqref{Srho-jumps}, \eqref{xnyn}, \eqref{jumpThetan},  \eqref{jumpThetai}, and \eqref{jumpSy}, that each $\Psi_{n,k}$ is holomorphic in $\RS\setminus\boldsymbol\Delta$ and its traces satisfy
\[
(W\Psi_{n,k})^+ = (W\Psi_{n,k})^-\exp\left\{\lambda_\rho\right\} \quad \text{on} \quad \boldsymbol\Delta.
\]
Using the definition of $W$ and \eqref{lambda-rho}, we see that the pullbacks of $\Psi_{n,k}$ solve the boundary value problem \eqref{BVP}. Further, the functions $\Psi_{n,k}$ have the local behavior around the branch points of $\RS$ as stated in \eqref{around-branching} by the very definition of $W$ and \eqref{Srho-BU} or \eqref{Srho-BU1}. Finally, the zero/pole multi-set of $\Psi_{n,k}$ in $\RS\setminus\boldsymbol\Delta$ is described by the divisor
\begin{equation}
\label{Psi-Divisor}
\mathcal{D}_{n,k} + (n+1)\left(\infty^{(1)} + \infty^{(2)}\right) - \infty^{(k)} - (2n-1)\infty^{(0)}
\end{equation}
as follows from the properties of $\Phi^n$, $\Theta_{n,k}$, $\Theta_k$, and $W$. In particular, the functions satisfying conditions of Theorem~\ref{thm:N-S} are given by $\Psi_n:=\Psi_{n,0}$.

It only remains to show the uniqueness of $\Psi_{n,k}$. To this end, consider $\Psi_{n,k}/\Psi$, where $\Psi$ satisfies \eqref{around-branching}, \eqref{BVP}, and \eqref{Psi-Divisor} with $\mathcal{D}_{n,k}$ replaced by any other integral divisor, say $\mathcal{D}$. Then the analytic continuation property implies that $\Psi_{n,k}/\Psi$ is a rational function on $\RS$ whose divisor is given by $\mathcal{D}_{n,k}-\mathcal{D}$. As $\mathcal{D}_{n,k}$ is not special, $\mathcal{D}=\mathcal{D}_{n,k}$, that is $\Psi_{n,k}/\Psi$ is a constant.

\subsection{Auxiliary Observation}

In Section~\ref{sec:CaseII}, we shall need the following observation. Assuming we are in Case II, let $\mathcal D_{n,k}(\rho)$ be the solution of \eqref{JIP} for the weights $\rho_1$ and $\rho_2$ and $\mathcal D_{n,k}(1/\rho)$ be the solution of \eqref{JIP} for the weights $1/\rho_1$ and $1/\rho_2$. Then it can easily be seen from \eqref{JIP}, \eqref{lambda-rho}, and \eqref{not-rational2} that
\[
\vec\Omega\big(\mathcal D_{n,k}(\rho)+\mathcal D_{n,k}(1/\rho)\big) \equiv 2\vec\Omega\big(\infty^{(2)}\big) + 2\vec v_k + 2\vec c_1 \equiv 2\vec\Omega(\mathcal D_{n,k}(1)).
\]
Moreover, it also follows from \eqref{not-rational2} and \eqref{BVP} that $\big(S_1S_{\vec y_n}\Theta_{n,k}\Theta_k\big)^2$ is a rational function on $\RS$ with the divisor equal to the sum of the divisors \eqref{(f)} and $2\mathcal D_{n,k}(1)-2\infty^{(k)}$. As the divisor \eqref{(f)} is a divisor of a rational function, we get from Abel's theorem \eqref{Abel} that
\[
\vec\Omega\big(\mathcal D_{n,k}(\rho)+\mathcal D_{n,k}(1/\rho)\big) \equiv 2\vec\Omega\big(\infty^{(k)}\big).
\]
In particular, it follows from the unique solvability of the Jacobi inversion problem that
\begin{equation}
\label{simult}
\mathcal D_{n,k}(\rho) = \infty^{(k)} \quad \Leftrightarrow \quad \mathcal D_{n,k}(1/\rho) = \infty^{(k)}.
\end{equation}
That is, $n\in\N_*(\rho)$ if and only if $n\in\N_*(1/\rho)$.

\subsection{Asymptotic Behavior}

Below we prove Theorem~\ref{thm:N-Sasymp}. In fact, we shall show that a more general statement is true. Given $\N_*$ and $\varepsilon>0$ small, there exists a constant $1<C_\varepsilon(\N_*)<\infty$ such that
\begin{equation}
\label{Psi-Estimates}
\left\{
\begin{array}{lll}
\left|\Psi_{n,k}\right| \leq C_\varepsilon(\N_*)\left|\Phi\right|^n & \text{in} & \RS_\varepsilon, \medskip \\
\left|\Psi_{n,k}\right| \geq C_\varepsilon(\N_*)^{-1}\left|\Phi\right|^n & \text{in} & \mathfrak{U}^{(k)}_\varepsilon,
\end{array}
\right.
\end{equation}
for all $n\in\N_*$ large, where $\mathfrak{U}^{(k)}_\varepsilon:=\RS^{(k)}\cap\pi^{-1}\big(\{|z|>1/\varepsilon\}\big)$.

It follows from the continuity of the boundary values of $S_\rho$ as well as from \eqref{Srho-BU} that $|S_\rho|$ is bounded from above and away from zero in $\RS_\varepsilon$ (that is, including the boundary values on $\boldsymbol\Delta$ and the $\ualpha$-cycles). Moreover, since the image of $\RS_{\ualpha,\ubeta}$ under $\vec\Omega$ is bounded and $\vec y_n\in[0,1)^g$ by the very definition, see \eqref{xnyn}, the moduli $|S_{\vec y_n}|$ are uniformly bounded from above and away from zero. Thus, we only need to estimate the function $\Theta_{n,k}\Theta_kW^{-1}$.

It follows from \eqref{Psi-Divisor} that the zero/pole multi-set of $\Theta_{n,k}\Theta_kW^{-1}$ in $\RS_\ualpha\setminus\boldsymbol\Delta$ is equal to $\mathcal{D}_{n,k} + \infty^{(0)} + \infty^{(1)} + \infty^{(2)} - \infty^{(k)}$ and therefore these functions are holomorphic there as well as non-vanishing in $\mathfrak{U}^{(k)}_\varepsilon$ for all $\varepsilon>0$ small and $n$ large by the very choice of $\N_*$. Hence,
\begin{equation}
\label{Psi-Est1}
\left\{
\begin{array}{lll}
\left|\Theta_{n,k}\Theta_kW^{-1}\right| \leq C_{n,\varepsilon}(\N_*) & \text{in} & \RS_\varepsilon, \medskip \\
\left|\Theta_{n,k}\Theta_kW^{-1}\right| \geq C_{n,\varepsilon}(\N_*)^{-1} & \text{in} & \mathfrak{U}^{(k)}_\varepsilon.
\end{array}
\right.
\end{equation}
Thus, if we show that $\big\{\Theta_{n,k}\Theta_kW^{-1}\big\}$ is a normal family and every limit point satisfies the estimates as in \eqref{Psi-Est1}, a standard compactness argument will finish the proof of \eqref{Psi-Estimates}. The only varying part of this family is given by
\begin{equation}
\label{Psi-Est2}
\theta\left(\vec\Omega(\z) - \vec\Omega\big(g\infty^{(2)}\big) - \vec K - \vec v_k - \vec c_\rho - \vec x_n - \mathbf{B}\vec y_n \right).
\end{equation}
Boundedness of its image follows from the continuity of $\theta(\cdot)$, boundedness of the image of $\vec\Omega(\z)$ for $\z\in\RS_{\ualpha,\ubeta}$, and the fact that $\vec x_n, \vec y_n\in[0,1)^g$. Thus, $\big\{\Theta_{n,k}\Theta_kW^{-1}\big\}$ is a normal family and its elements are indexed by pairs $(\vec x_n,\vec y_n)\in[0,1)^{2g}$. The continuity of $\theta(\cdot)$ and \eqref{thetan} imply that any limit point of the functions in \eqref{Psi-Est2} has the same form with $(\vec x_n,\vec y_n)$ replaced by a limit point $(\vec x,\vec y)\in[0,1]^{2g}$. Denote by $\mathcal{D}$ the integral divisor of order $g$ describing the zeros of this limit point. Then
\[
\vec\Omega(\mathcal{D}) \equiv \vec\Omega\big(g\infty^{(2)}\big) + \vec v_k  + \vec c_\rho + \vec x + \mathbf B\vec y \quad \left(\mdp \mathrm{d}\vec\Omega\right).
\]
The continuity of $\vec\Omega$ implies that the right-hand side of the equivalence above is a limit point of the right-hand sides in \eqref{JIP}.  By the definition of $\N_*$, it determines $\mathcal{D}$ uniquely as it is a limit point of $\{\mathcal{D}_{n,k}\}_{n\in\N_*}$.  Moreover, all such limit points are uniformly bounded away from containing $\infty^{(k)}$ and therefore from containing an element from $\mathfrak{U}^{(k)}_\varepsilon$ for all $\varepsilon>0$ small enough, which is sufficient to get the required lower estimate.

\section{Riemann-Hilbert Analysis: Case I}
\label{sec:CaseI}

\subsection{Global Lenses}

Set
\begin{equation}
\label{global-matr}
\boldsymbol G_1(u) := \left(
\begin{array}{ccc}
1 & 0 & 0 \smallskip \\
0 & 1 & 0  \smallskip \\
0 & u & 1
\end{array}
\right)
\quad \text{and} \quad
\boldsymbol G_2(v) := \left(
\begin{array}{ccc}
1 & 0 & 0 \smallskip \\
0 & 1 & v  \smallskip \\
0 & 0 & 1
\end{array}
\right).
\end{equation}
Observe that $\boldsymbol G_i(u)^{-1} = \boldsymbol G_i(-u)$, $i\in\{1,2\}$. Moreover, with $J$ as in \eqref{J},
\begin{equation}
\label{global-matr-modif}
\left\{
\begin{array}{l}
\boldsymbol G_1^{-1}(u)\boldsymbol J(x,y)\boldsymbol G_1(u) = \boldsymbol J(x+uy,y) \medskip \\
\boldsymbol G_2^{-1}(v)\boldsymbol J(x,y)\boldsymbol G_2(v) = \boldsymbol J(x,y+vx).
\end{array}
\right.
\end{equation}
That is, conjugating $\boldsymbol J$ by $\boldsymbol G_i$ allows one to modify the $i$-th variable of the matrix $\boldsymbol J$ using the other variable.

\begin{figure}[!ht]
\centering
\includegraphics[scale=.5]{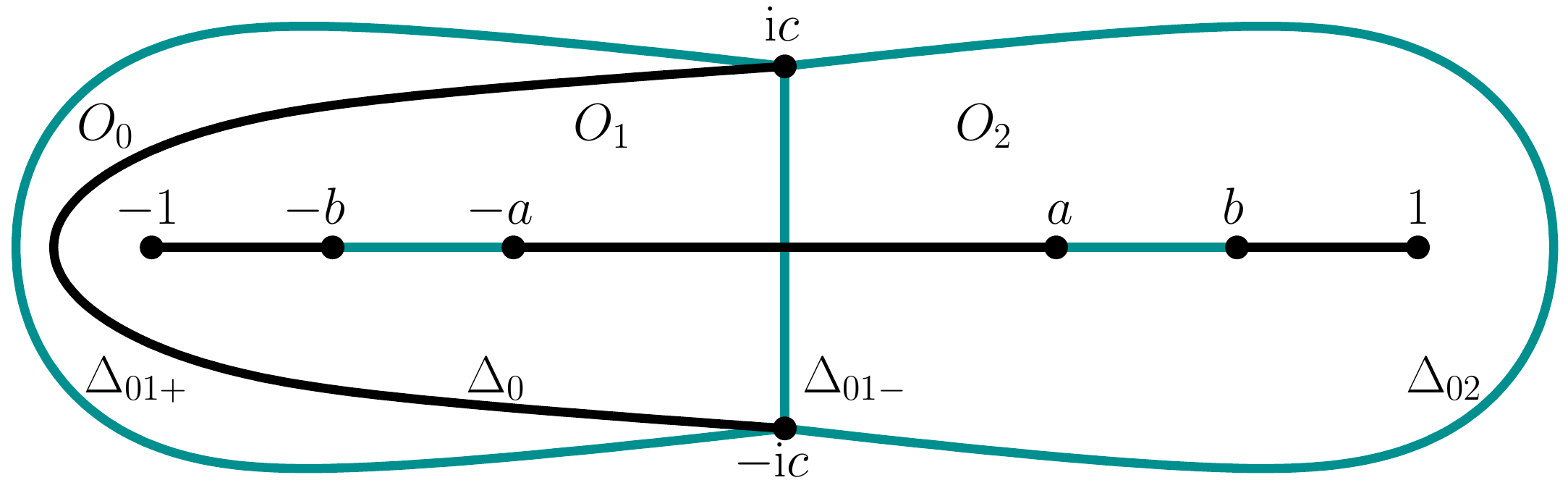}
\caption{\small The lens $\Sigma_{\boldsymbol S}$ and the domains $O_0$, $O_1$, and $O_2$.}
\label{fig:Global-I}
\end{figure}

Let, as before, $\pm\mathrm{i} c$, $c>0$, be the projections of the branch points of $\RS$ and $\Gamma$ be the cut along which the sheets $\RS^{(1)}$ and $\RS^{(2)}$ are glued to each other. Pick arcs $\Delta_{01\pm}$ and $\Delta_{02}$, all oriented from $-\mathrm{i}c$ to $\mathrm{i}c$, as in Figure~\ref{fig:Global-I}. These arcs delimit three bounded domains that we label from left to right as $ O_0$, $O_1$, and $O_2$.

Set $\boldsymbol S:=\boldsymbol Y$ outside of the closure of $O_0\cup O_1\cup O_2$ and inside put
\begin{equation}
\label{S}
\boldsymbol S := \boldsymbol Y\left\{
\begin{array}{lcl}
\boldsymbol G_1(-\rho_1/\rho_2) & \text{in} & O_0, \medskip \\
\boldsymbol G_2(-\rho_2/\rho_1)\boldsymbol G_1(\rho_1/\rho_2) & \text{in} & O_1, \medskip \\
\boldsymbol G_2(-\rho_2/\rho_1) & \text{in} & O_2,
\end{array}
\right.
\end{equation}
where the $\rho_i$ are the weight functions in \eqref{appr-fun}.
This way $\boldsymbol S$ has jumps on $\Delta_{01\pm}$, $\Delta_{02}$, and $\Delta_0$. In particular, the  jump on $\Delta_0$ is equal to
\[
\boldsymbol J_0 := \boldsymbol G_1(-\rho_1/\rho_2)\boldsymbol G_2(\rho_2/\rho_1)\boldsymbol G_1(-\rho_1/\rho_2) = \left(
\begin{array}{ccc}
1 & 0 & 0 \smallskip \\
0 & 0 & \rho_2/\rho_1  \smallskip \\
0 & -\rho_1/\rho_2 & 0
\end{array}
\right).
\]
For brevity, let us write
\begin{equation}
\label{rhostar}
\rho_2^*(z) : =\left\{
\begin{array}{rl}
\rho_2(z), & \re(z)>0, \medskip \\
-\rho_2(z), & \re(z)<0.
\end{array}
\right.
\end{equation}

Put $\Sigma_{\boldsymbol S}:=[-1,1]\cup\Delta_{01+}\cup\Delta_{01-}\cup\Delta_{02}\cup\Delta_0$. Then $\boldsymbol S$ solves \rhs:
\begin{itemize}
\label{rhs}
\item[(a)] $\boldsymbol S$ is analytic in $\C\setminus\Sigma_{\boldsymbol S}$ and $\lim_{z\to\infty} {\boldsymbol S}(z)\diag\left(z^{-2n},z^n,z^n\right) = \boldsymbol I$;
\item[(b)] ${\boldsymbol S}$ has continuous traces on $\Sigma_{\boldsymbol S}^\circ:=\Sigma_{\boldsymbol S}\setminus\{\pm1,\pm a,\pm\mathrm{i}c\}$ that satisfy
\[
{\boldsymbol S}_+ = {\boldsymbol S}_-
\left\{
\begin{array}{lll}
{\boldsymbol J}(0,\rho_2^*) & \text{on} & (-1,-a)\cup(a,1), \smallskip \\
{\boldsymbol J}(\rho_1,0) & \text{on} & (-a,a)\setminus\{0\}, \smallskip \\
\boldsymbol J_0, & \text{on} & \Delta_0, \smallskip \\
{\boldsymbol G}_1(\rho_1/\rho_2) & \text{on} & \Delta_{01\pm}, \smallskip \\
{\boldsymbol G}_2(-\rho_2/\rho_1) & \text{on} & \Delta_{02};
\end{array}
\right.
\]
\item[(c)] $\boldsymbol S$ satisfies \hyperref[rhy]{\rhy}(c) (see Section \ref{MOP}) with $[-1,1]$ replaced by $\Sigma_{\boldsymbol S}$.
\end{itemize}
If  \hyperref[rhs]{\rhs} is solvable, then so is  \hyperref[rhy]{\rhy}, and the solutions are connected via \eqref{S}.

\subsection{Local Lenses}

\begin{figure}[!ht]
\centering
\includegraphics[scale=.5]{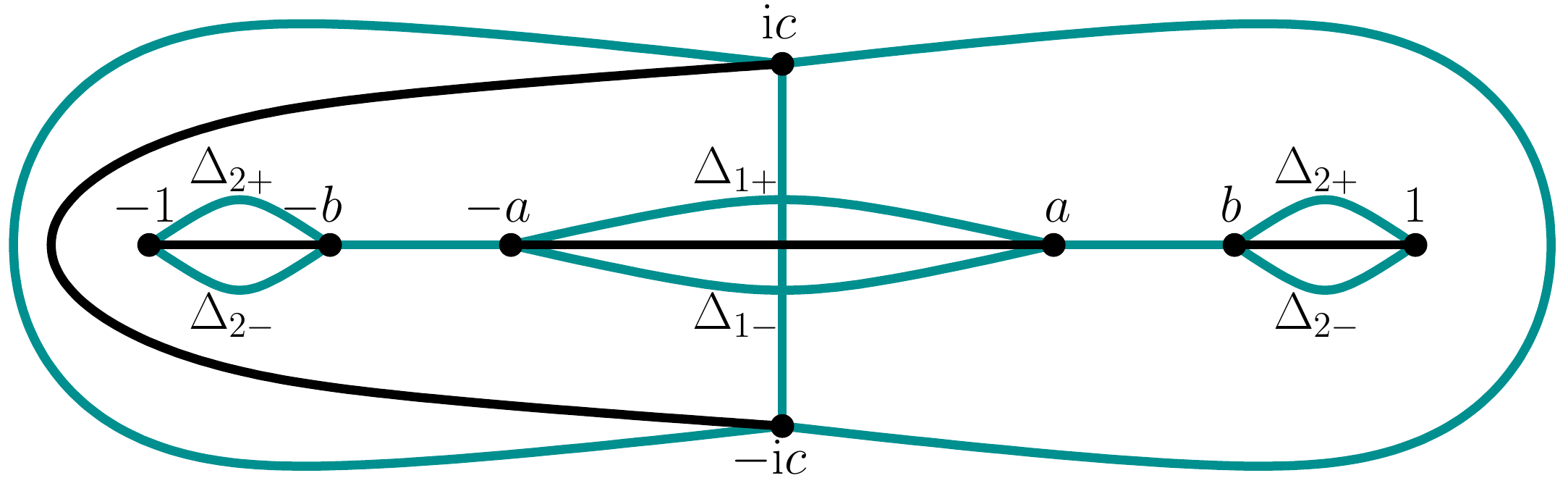}
\caption{\small The lens $\Sigma_{\boldsymbol Z}$: specifically, the arcs $\Delta_{1\pm}$ and systems of arcs $\Delta_{2\pm}$. The black curves constitute the system $\Sigma_{\boldsymbol N}$.}
\label{fig:Local-I}
\end{figure}

For the next step we introduce additional arcs $\Delta_{1\pm}$ and systems of two arcs $\Delta_{2\pm}$ as on Figure~\ref{fig:Local-I}, all oriented from left to right.
 Denote further by $O_{1\pm}$ the domains bounded by $(-a,a)$ and the arcs $\Delta_{1\pm}$, respectively, and by $O_{2\pm}$ the open sets bounded by $\Delta_{2\pm}$ and $(-1,-b)\cup(b,1)$. Set
\begin{equation}
\label{Z}
\boldsymbol Z:= \boldsymbol S\boldsymbol L_i^{\mp1} \quad \text{in} \quad O_{i\pm},
\end{equation}
$i\in\{1,2\}$, where
\begin{equation}
\label{local-matr}
\boldsymbol L_1 := \left(
\begin{array}{ccc}
1 & 0 & 0 \smallskip \\
\displaystyle 1/\rho_1 & 1 & 0  \smallskip \\
0 & 0 & 1
\end{array}
\right)
\quad \text{and} \quad
{\boldsymbol L}_2 := \left(
\begin{array}{ccc}
1 & 0 & 0 \smallskip \\
0 & 1 & 0  \smallskip \\
\displaystyle 1/\rho_2^* & 0 & 1
\end{array}
\right).
\end{equation}
It can readily be checked that ${\boldsymbol Z}_+ = {\boldsymbol Z}_-\boldsymbol J_i$ on $\Delta_i$, $i\in\{1,2\}$, see \eqref{chains}, where
\begin{equation}
\label{local-jumps}
{\boldsymbol J}_1 :=
\left(
\begin{array}{ccc}
0 & \rho_1 & 0 \smallskip \\
-1/\rho_1 & 0 & 0  \smallskip \\
0 & 0 & 1
\end{array}
\right)
\quad \text{and} \quad
{\boldsymbol J}_2 :=
\left(
\begin{array}{ccc}
0 & 0 & \rho_2^* \smallskip \\
0 & 1 & 0  \smallskip \\
\displaystyle -1/\rho_2^* & 0 & 0
\end{array}
\right).
\end{equation}

Put $\Sigma_{\boldsymbol Z} := \Sigma_{\boldsymbol S}\cup\Delta_{1+}\cup\Delta_{1-}\cup\Delta_{2+}\cup\Delta_{2-}$. Then ${\boldsymbol Z}$ solves \rhz:
\begin{itemize}
\label{rhz}
\item[(a)] ${\boldsymbol Z}$ is analytic in $\C\setminus\Sigma_{\boldsymbol Z}$ and $\lim_{z\to\infty} {\boldsymbol Z}(z)\diag\left(z^{2n},z^{-n},z^{-n}\right) = \boldsymbol I$;
\item[(b)] ${\boldsymbol Z}$ has continuous traces on each side of $\Sigma_{\boldsymbol Z}^\circ:=\Sigma_{\boldsymbol Z}\setminus\{\pm1,\pm a,\pm b,\pm\mathrm{i}c\}$ that satisfy
\[
{\boldsymbol Z}_+ = {\boldsymbol Z}_-
\left\{
\begin{array}{lll}
{\boldsymbol J}_k & \text{on} & \Delta_k^\circ\setminus\{0\}, \quad k\in\{0,1,2\}, \medskip \\
\boldsymbol J_{\boldsymbol Z} &\text{on} & \Sigma_{\boldsymbol Z}^\circ\setminus\big(\Delta_0\cup\Delta_1\cup\Delta_2\big),
\end{array}
\right.
\]
where
\begin{equation}
\label{J-Z}
\boldsymbol J_{\boldsymbol Z} :=
\left\{
\begin{array}{lll}
{\boldsymbol J}(0,\rho_2^*) & \text{on} & (-b,-a)\cup(a,b), \medskip \\
{\boldsymbol G}_1(\rho_1/\rho_2) & \text{on} & \Delta_{01\pm}, \medskip \\
{\boldsymbol G}_2(-\rho_2/\rho_1) & \text{on} & \Delta_{02}, \medskip \\
{\boldsymbol L}_i & \text{on} & \Delta_{i\pm}, \quad i\in\{1,2\};
\end{array}
\right.
\end{equation}
\item[(c)] $\boldsymbol Z$ satisfies \hyperref[rhs]{\rhs}(c) with $\Sigma_{\boldsymbol S}$ replaced by $\Sigma_{\boldsymbol Z}$.
\end{itemize}
As before, if \hyperref[rhz]{\rhz} is solvable, then so is \hyperref[rhs]{\rhs}, and the solutions are connected via \eqref{Z}.

\subsection{Global Parametrix}
\label{ssec:GP-I}

Let $\Sigma_{\boldsymbol N}:=\Delta_0\cup\Delta_1\cup\Delta_2$, see Figure~\ref{fig:Local-I}. In this section we are seeking the solution of the following Riemann-Hilbert problem (\rhn):
\begin{itemize}
\label{rhn}
\item[(a)] ${\boldsymbol N}$ is analytic in $\C\setminus\Sigma_{\boldsymbol N}$ and $\lim_{z\to\infty} {\boldsymbol N}(z)\diag\left(z^{-2n},z^n,z^n\right) = \boldsymbol I$;
\item[(b)] ${\boldsymbol N}$ has continuous traces on each side of $\Sigma_{\boldsymbol N}\setminus\{\pm1,\pm a,\pm b,\pm\mathrm{i}c\}$ that satisfy $\boldsymbol N_+ = \boldsymbol N_- \boldsymbol J_k$ on $\Delta_k^\circ$, $k\in\{0,1,2\}$.\end{itemize}

We solve \hyperref[rhn]{\rhn} only along the subsequence $\N_*$ defined in Proposition~\ref{prop:genus2}. Let $\Psi_{n,k}$ be the $n$-th Nuttall-Szeg\H{o} functions constructed in \eqref{Psink}. As the functions $\Psi_{n,k}$ satisfy \eqref{BVP} with the zero/pole sets described by \eqref{Psi-Divisor}, it readily follows that a solution of \hyperref[rhn]{\rhn} is given by
\begin{equation}
\label{N}
\boldsymbol N := \boldsymbol C_n\left(
\begin{array}{ccc}
\Psi_{n,0}^{(0)} & \Psi_{n,0}^{(1)} & \Psi_{n,0}^{(2)} \smallskip \\
\Psi_{n,1}^{(0)} & \Psi_{n,1}^{(1)} & \Psi_{n,1}^{(2)} \smallskip \\
\Psi_{n,2}^{(0)} & \Psi_{n,2}^{(1)} & \Psi_{n,2}^{(2)}
\end{array}
\right) =: \boldsymbol C^{-n}\boldsymbol {MD}^n,
\end{equation}
where $n\in\N_*$, $\boldsymbol C_n$ is a diagonal matrix of constants chosen to fulfill the normalization condition in \hyperref[rhn]{\rhn}(a), $\boldsymbol C  := \diag(C_0,C_1,C_2)$, see \eqref{Phi_jatInfty}, and $\boldsymbol D:=\diag\left(\Phi^{(0)}, \Phi^{(1)}, \Phi^{(2)}\right)$.

Observe that $\det(\boldsymbol N)$ is an entire function in $\C \setminus \{\pm1,\pm a,\pm b,\pm\mathrm{i}c\}$ since the determinants of the jump matrices in \hyperref[rhn]{\rhn}(b) are all equal to 1. Moreover, the normalization at infinity implies that $\det(\boldsymbol N)(\infty)=1$ and therefore $\det(\boldsymbol N)$ is a rational function. It also follows from \eqref{around-branching} that near any of the points in $\{\pm1,\pm a,\pm b,\pm\mathrm{i}c\}$ the entries of one of the columns of $\boldsymbol N$ are bounded and the other entries behave like   $\mathcal{O}(|z-e|^{-1/4})$. Hence,
\[
\det(\boldsymbol N)(z) = \mathcal{O}\left(|z-e|^{-1/2}\right) \quad \text{as} \quad z\to\left\{\pm1,\pm a,\pm b,\pm\mathrm{i}c\right\},
\]
and therefore $\det(\boldsymbol N)\equiv1$.

\subsection{Local Parametrices}
\label{ssec:local-param}

Denote by $U_e$, $e\in\{\pm1,\pm a,\pm b,\pm\mathrm{i}c\}$, a disk centered at $e$ of small enough radius so that $\Sigma_e:=U_e\cap\Sigma_{\boldsymbol Z}$ consists of disjoint, except at $e$, analytic arcs. We are seeking a solution of the following \rhp$_e$:
\begin{itemize}
\label{rhp}
\item[(a)] $\boldsymbol P_e$ is analytic in $U_e\setminus\Sigma_e$;
\item[(b)] $\boldsymbol P_e$ has continuous traces on each side of $\Sigma_e^\circ$ that satisfy \hyperref[rhz]{\rhz}(b) within $U_e$;
\item[(c)] $\boldsymbol P_e$ is either bounded or has the behavior near $e$ within $U_e$ described by \hyperref[rhz]{\rhz}(c);
\item[(d)] $\boldsymbol P_e=\boldsymbol M\left(\boldsymbol I+\boldsymbol{\mathcal{O}}(1/n)\right)\boldsymbol D^n$ uniformly on $\partial U_e\setminus\Sigma_{\boldsymbol Z}$, where $\boldsymbol M$ and $\boldsymbol D$ are defined by \eqref{N}.
\end{itemize}
We solve each \hyperref[rhp]{\rhp$_e$} only for $n\in\N_*$. For these indices the above problem is well-posed as $\det(\boldsymbol N)\equiv1$ and therefore $\boldsymbol N^{-1}$ is an analytic matrix function in $\C\setminus\Sigma_{\boldsymbol N}$. In fact, the solution does not depend on the actual choice of $\N_*$, however, the term $\mathcal{O}(1/n)$ in \hyperref[rhp]{\rhp$_e$}(d) may depend on the choice of this subsequence.

In solving \hyperref[rhp]{\rhp$_e$}, it will be convenient to use the notation $\sigma_3=\left(\begin{matrix} 1 & 0 \\ 0 & -1 \end{matrix}\right)$ and introduce the following matrix transformations $\mathsf{T}_1\boldsymbol A$, $\mathsf{T}_2\boldsymbol A$, and $\mathsf{T}_3\boldsymbol A$ given by
\[
\left(
\begin{array}{ccc}
1 & 0 & 0 \medskip \\
0 & [\boldsymbol A]_{11} & [\boldsymbol A]_{12} \medskip \\
0 & [\boldsymbol A]_{21} & [\boldsymbol A]_{22}
\end{array}
\right),
\quad
\left(
\begin{array}{ccc}
~[\boldsymbol A]_{11} & 0 & [\boldsymbol A]_{12} \medskip \\
0 & 1 & 0 \medskip \\
~[\boldsymbol A]_{21} & 0 & [\boldsymbol A]_{22}
\end{array}
\right),
\quad \text{and} \quad
\left(
\begin{array}{ccc}
~[\boldsymbol A]_{11} & [\boldsymbol A]_{12} & 0 \medskip \\
~[\boldsymbol A]_{21} & [\boldsymbol A]_{22} & 0 \medskip \\
0 & 0 & 1
\end{array}
\right),
\]
respectively, where $[\boldsymbol A]_{ik}$ is the $(i,k)$-th entry of $\boldsymbol A$. Observe that $\mathsf{T}_j(\boldsymbol{AB})=\mathsf{T}_j(\boldsymbol{A})\mathsf{T}_j(\boldsymbol{B})$ for all $2\times2$ matrices $\boldsymbol A,\boldsymbol B$.

\subsubsection{\hyperref[rhp]{\rhp$_{\pm1}$}}

In \cite{KMcLVAV04}, a $2\times2$ matrix function was constructed out of Bessel and Hankel functions that solves \rhpsi:
\begin{itemize}
\label{rhpsi}
\item[(a)] $\boldsymbol\Psi$ is holomorphic in $\C\setminus(I_+\cup I_-\cup(-\infty,0])$, where $I_\pm := \big\{\zeta:\arg(\zeta)=\pm2\pi/3\big\}$, all the rays are oriented towards the origin;
\item[(b)] $\boldsymbol\Psi$ has continuous traces on $I_+\cup I_-\cup(-\infty,0)$ that satisfy
\[
\boldsymbol\Psi_+ = \boldsymbol\Psi_-
\left\{
\begin{array}{lll}
\left(\begin{matrix} 1 & 0 \\ 1 & 1 \end{matrix}\right) & \text{on} & I_\pm, \medskip \\
\left(\begin{matrix} 0 & 1 \\ -1 & 0 \end{matrix}\right) & \text{on} & (-\infty,0);
\end{array}
\right.
\]
\item[(c)] $\boldsymbol\Psi(\zeta)=\mathcal{O}(\log|\zeta|)$ as $\zeta\to0$;
\item[(d)] $\boldsymbol\Psi$ has the following behavior near $\infty$:
\[
\boldsymbol\Psi(\zeta) = \left(2\pi\zeta^{1/2}\right)^{-\sigma_3/2}\frac1{\sqrt2}\left(\begin{matrix} 1 & \mathrm{i} \\ \mathrm{i} & 1 \end{matrix}\right)\left(\boldsymbol I+\boldsymbol{\mathcal{O}}\left(\zeta^{-1/2}\right)\right)\exp\left\{2\zeta^{1/2}\sigma_3\right\}
\]
uniformly in $\C\setminus(I_+\cup I_-\cup(-\infty,0])$.
\end{itemize}
Furthermore, $\sigma_3\boldsymbol\Psi\sigma_3$ solves the same R-H problem only with the reversed orientation of all the rays. Notice also that \hyperref[rhpsi]{\rhpsi}(d) should be replaced in this case by
\[
\sigma_3\boldsymbol\Psi(\zeta)\sigma_3 = \left(2\pi\zeta^{1/2}\right)^{-\sigma_3/2}\frac1{\sqrt2}\left(\begin{matrix} 1 & -\mathrm{i} \\ -\mathrm{i} & 1 \end{matrix}\right)\left(\boldsymbol I+\boldsymbol{\mathcal{O}}\left(\zeta^{-1/2}\right)\right)\exp\left\{2\zeta^{1/2}\sigma_3\right\}.
\]

To carry $\boldsymbol\Psi$ from the $\zeta$-plane to $U_e$, $e\in\{\pm1\}$, we need to introduce local conformal maps. To this end, set
\[
g_e(z) := \frac12\int_e^z\left(h_0-h_2\right)(t)\,\mathrm{d}t = \frac12\log\left(\Phi^{(0)}(z)/\Phi^{(2)}(z)\right),
\]
$z\in U_e\setminus[-1,1]$, where the second equality follows from \eqref{GreenDiff} and \eqref{Phi}. As mentioned after \eqref{likearoot}, $g_e$ has purely imaginary traces on $(-1,1)\cap U_e$ that differ by a sign. Moreover, since $g_e$ vanishes at $e$ as a square root, $g_e^2$ is conformal in $U_e$. Furthermore,
\[
\left\{
\begin{array}{rcl}
g_e^2\big(\{x:~\sgn(e)x>1\}\cap U_e\big) &\subset& \big\{z:z>0\big\}, \medskip \\
g_e^2\big((-1,1)\cap U_e\big) &\subset& \big\{z:z<0\big\}, \medskip \\
g_e^2\big(\Delta_{2\pm}\cap U_e\big) &\subset& \big\{z:~\sgn(e)\arg(z)=\pm2\pi/3\big\}.
\end{array}
\right.
\]
Indeed, the first property follows from \eqref{label} and \eqref{arrange3} while the second is a consequence of the fact that $g_e$ has purely imaginary traces there. The last property is the requirement we impose on the arcs $\Delta_{2\pm}$. Choosing the branch of $g_e^{1/2}$ which is positive on $\{x:~\sgn(e)x>1\}\cap U_e$, one has on $(-1,1)\cap U_e$ that
\begin{equation}
\label{ge1-root}
g_{e+}^{1/2} = \sgn(e)\mathrm{i}g_{e-}^{1/2}.
\end{equation}

Now, it can readily be verified that the matrix function
\begin{equation}
\label{sol-1}
\boldsymbol P_e := \boldsymbol E_e\mathsf{T}_2\boldsymbol\Psi_e\left(n^2g_e^2/4\right)\boldsymbol W_e,
\end{equation}
satisfies \hyperref[rhp]{\rhp$_e$}(a,b,c) for any holomorphic matrix function $\boldsymbol E_e$, where $\boldsymbol\Psi_1:=\boldsymbol\Psi$, $\boldsymbol\Psi_{-1}:=\sigma_3\boldsymbol\Psi\sigma_3$, and
\begin{equation}
\label{W_e}
\boldsymbol W_e := \diag\left(\big(\Phi^{(0)}\Phi^{(2)}\big)^{n/2} / \sqrt{\rho_2^*}, \big(\Phi^{(1)}\big)^n, \big(\Phi^{(0)}\Phi^{(2)}\big)^{n/2}\sqrt{\rho_2^*}\right).
\end{equation}
Moreover, one has on $\partial U_e$ that
\begin{equation}
\label{solution-1}
\boldsymbol P_e = \boldsymbol E_e \mathsf{T}_2\left( \rho_2^{*-\sigma_3/2} \left(\pi ng_e\right)^{-\sigma_3/2}\frac1{\sqrt2}\left(\begin{matrix} 1 & \sgn(e)\mathrm{i} \\ \sgn(e)\mathrm{i} & 1 \end{matrix}\right)\right)\big(\boldsymbol I+\boldsymbol{\mathcal{O}}(1/n)\big)\boldsymbol D^n.
\end{equation}
Thus, it only remains to choose $\boldsymbol E_e$ so that \hyperref[rhp]{\rhp$_e$}(d) is fulfilled. Direct computation using \hyperref[rhn]{\rhn}(b) and \eqref{ge1-root} shows that the matrix
\begin{equation}
\label{E_e}
\boldsymbol E_e := \boldsymbol M\mathsf{T}_2\left(\rho_2^{*\sigma_3/2}\left(\pi ng_e\right)^{\sigma_3/2} \frac1{\sqrt2}\left(\begin{matrix} 1 & -\sgn(e)\mathrm{i} \\ -\sgn(e)\mathrm{i} & 1 \end{matrix}\right)\right)
\end{equation}
is holomorphic in $U_e\setminus\{e\}$. As the entries of $\boldsymbol M$ behave like $|z-e|^{-1/4}$ as $z\to e$ by \eqref{around-branching}, the entries of $\boldsymbol E_e$ can have at most a square root singularity there. Thus, $\boldsymbol E_e$ is holomorphic throughout $U_e$ as desired.

\subsubsection{\hyperref[rhp]{\rhp$_{\pm a}$}}

To solve \hyperref[rhp]{\rhp$_{\pm a}$}, we again use the matrix $\boldsymbol\Psi$. In this case we have an additional complication coming from the jumps on $(-b,-a)\cup(a,b)$. To circumvent it, we shall need the following fact about the matrix $\boldsymbol\Psi$:
\begin{equation}
\label{bessel}
[\boldsymbol\Psi]_{11}(\zeta) = I_0(2\zeta^{1/2}) \quad \text{and} \quad [\boldsymbol\Psi]_{21}(\zeta) = 2\pi\mathrm{i}\zeta^{1/2}I_0^\prime(2\zeta^{1/2})
\end{equation}
within $|\arg(z)|<2\pi/3$, where $I_0$ is the modified Bessel function of order $0$. Observe that both functions above are in fact entire in the whole complex plane.

Define
\begin{equation}
\label{connect-a}
g_e(z) := \frac12\int_e^z\left(h_0-h_1\right)(t)\mathrm{d}t = \frac12\log\left(\Phi^{(0)}(z)/\Phi^{(1)}(z)\right),
\end{equation}
$z\in U_e\setminus[-a,a]$. As mentioned after \eqref{likearoot1}, $g_e$ has purely imaginary traces on $(-a,a)\cap U_e$ that differ by a sign. Moreover, since $g_e$  vanishes at $e$ as a square root, $g_e^2$ is conformal in $U_e$. Furthermore,
\[
\left\{
\begin{array}{rcl}
g_e^2\big(\{x:~\sgn(e)x>a\}\cap U_e\big) &\subset& \big\{z:z>0\big\}, \medskip \\
g_e^2\big((-a,a)\cap U_e\big) &\subset& \big\{z:z<0\big\}, \medskip \\
g_e^2\big(\Delta_{1\pm}\cap U_e\big) &\subset& \big\{z:~\sgn(e)\arg(z)=\pm2\pi/3\big\},
\end{array}
\right.
\]
where the first property follows from \eqref{arrange1} and \eqref{arrange2}, the second is a consequence of the fact that $g_e$ has purely imaginary traces there, and the third is a requirement we impose on the arcs $\Delta_{1\pm}$. Choosing the branch of $g_e^{1/2}$ which is positive on $\{x:~\sgn(e)x>a\}\cap U_e$, we see that \eqref{ge1-root} holds on $(-a,a)\cap U_e$.

We further define
\[
\left\{
\begin{array}{lll}
F_{n1}(z) &:=& \rho_2(z)\big(\Phi^{(2)}(z)\big)^{-3n/2} I_0(ng_e(z))\frac1{2\pi\mathrm{i}}\log\frac{z-b}{z-a},\medskip \\
F_{n2}(z) &:=& \rho_2(z)\big(\Phi^{(2)}(z)\big)^{-3n/2} \pi\mathrm{i}nI_0^\prime(ng_e(z))\frac1{2\pi\mathrm{i}}\log\frac{z-b}{z-a},
\end{array}
\right.
\]
where $I_0$ is the modified Bessel function of order 0, see \eqref{bessel}. The above functions are holomorphic in $U_a\setminus[a,b]$ and
\[
F_{ni+} - F_{ni-} = [\boldsymbol\Psi]_{i1}\big(n^2g_a^2/4\big)\rho_2\big(\Phi^{(2)}\big)^{-3n/2}
\]
on $(a,b)$. According to \cite[Eq. 10.40.5]{DLMF},
\begin{eqnarray*}
(ng_a)^{1/2} F_{n1} &=& \left(e^{ng_a}\mathcal{O}(1) + e^{-ng_a}\mathcal{O}(1)\right)\big(\Phi^{(2)}\big)^{-3n/2} \\
&=& \left(\frac{\Phi^{(0)}}{\Phi^{(2)}}\right)^n\mathcal{O}(1) + \left(\frac{\Phi^{(1)}}{\Phi^{(2)}}\right)^n\mathcal{O}(1) = o(1)
\end{eqnarray*}
uniformly on $\partial U_a$ by \eqref{connect-a} and since $\partial U_a\subset\Omega_{201}$ by Theorem~\ref{thm:Omegas} and the choice of the radius of $U_a$. Clearly, $o(1)$ in the above equality is geometric. Moreover, a completely analogous estimate holds for $F_{n2}$. Given these two functions, we can set
\[
\boldsymbol\Psi_a := \left(
\begin{matrix}
~[\boldsymbol\Psi]_{11}\big(n^2g_a^2/4\big) & [\boldsymbol\Psi]_{12}\big(n^2g_a^2/4\big) & F_{n1} \medskip \\
~[\boldsymbol\Psi]_{21}\big(n^2g_a^2/4\big) & [\boldsymbol\Psi]_{22}\big(n^2g_a^2/4\big) & F_{n2} \medskip \\
0 & 0 & 1
\end{matrix}
\right).
\]
Then the solution of \hyperref[rhp]{\rhp$_a$} is given by $\boldsymbol P_a := \boldsymbol E_a\boldsymbol\Psi_a\boldsymbol W_a$, where
\[
\left\{
\begin{array}{lll}
\boldsymbol W_a & := & \diag\left( \big(\Phi^{(0)}\Phi^{(1)}\big)^{n/2}/\sqrt{\rho_1}, \big(\Phi^{(0)}\Phi^{(1)}\big)^{n/2}\sqrt{\rho_1}, \big(\Phi^{(2)}\big)^n  \right) \medskip \\
\boldsymbol E_a & := & \boldsymbol M\mathsf{T}_3\left(\rho_1^{\sigma_3/2}\left(\pi ng_a\right)^{\sigma_3/2} \frac1{\sqrt2} \left(\begin{matrix} 1 & -\mathrm{i} \\ -\mathrm{i} & 1 \end{matrix}\right)\right).
\end{array}
\right.
\]
To verify \hyperref[rhp]{\rhp$_a$}(d), observe that
\[
\mathsf{T}_3\left(\frac1{\sqrt2}\left( \begin{matrix} 1 & -\mathrm{i} \smallskip \\ -\mathrm{i} & 1 \end{matrix} \right) \big(\pi ng_a\big)^{\sigma_3/2}\right)\boldsymbol\Psi_a \mathsf{T}_3\left(\Phi^{(0)}/\Phi^{(1)}\right)^{-n\sigma_3/2} = \boldsymbol I + \boldsymbol{\mathcal{O}}(1/n)
\]
by the asymptotic properties of $F_{ni}$. The matrix $\boldsymbol P_{-a}$ solving \hyperref[rhp]{\rhp$_{-a}$} can be constructed analogously.

\subsubsection{\hyperref[rhp]{\rhp$_{\pm b}$}}

In \cite{DKMLVZ99b}, a $2\times2$ matrix was constructed out of Airy functions that solves \rhphi:
\begin{itemize}
\label{rhphi}
\item[(a)] $\boldsymbol\Phi$ is holomorphic in $\C\setminus(I_+\cup I_-\cup(-\infty,\infty))$, where the real line is oriented from left to right;
\item[(b)] $\boldsymbol\Phi$ has continuous traces on $I_+\cup I_-\cup(-\infty,0)\cup(0,\infty)$ that satisfy
\[
\boldsymbol\Phi_+ = \boldsymbol\Phi_-
\left\{
\begin{array}{lll}
\left(\begin{matrix} 1 & 0 \\ 1 & 1 \end{matrix}\right) & \text{on} & I_\pm, \medskip \\
\left(\begin{matrix} 0 & 1 \\ -1 & 0 \end{matrix}\right) & \text{on} & (-\infty,0), \medskip \\
\left(\begin{matrix} 1 & 1 \\ 0 & 1 \end{matrix}\right) & \text{on} & (0,\infty);
\end{array}
\right.
\]
\item[(c)] $\boldsymbol\Phi$ is bounded around the origin;
\item[(d)] $\boldsymbol\Phi$ has the following behavior near $\infty$:
\[
\boldsymbol\Phi(\zeta) = \zeta^{-\sigma_3/4}\frac1{\sqrt2}\left(\begin{matrix} 1 & \mathrm{i} \\ \mathrm{i} & 1 \end{matrix}\right)\left(\boldsymbol I+\boldsymbol{\mathcal{O}}\left(\zeta^{-3/2}\right)\right)\exp\left\{-\frac23\zeta^{3/2}\sigma_3\right\}
\]
uniformly in $\C\setminus(I_+\cup I_-\cup(-\infty,\infty))$.
\end{itemize}
Again, $\sigma_3\boldsymbol\Phi\sigma_3$ solves the same R-H problem only with the reversed orientation of all the rays. As in the case of $\boldsymbol\Psi$,
\[
\sigma_3\boldsymbol\Phi(\zeta)\sigma_3 = \zeta^{-\sigma_3/4}\frac1{\sqrt2}\left(\begin{matrix} 1 & -\mathrm{i} \\ -\mathrm{i} & 1 \end{matrix}\right)\left(\boldsymbol I+\boldsymbol{\mathcal{O}}\left(\zeta^{-3/2}\right)\right)\exp\left\{-\frac23\zeta^{3/2}\sigma_3\right\}.
\]

To map $\boldsymbol\Phi$ into $U_e$, $e\in\{\pm b\}$, define
\[
g_e(z) := -\frac12\int_e^z\left(h_0-h_2\right)(t)\mathrm{d}t = -\frac12\log\left( \Phi^{(0)}(z)/\Phi^{(2)}(z) \right),
\]
$z\in U_e\setminus\big([-1,-b]\cup[b,1]\big)$. As before, $g_e$ has purely imaginary traces on $\big((-1,-b)\cup(b,1)\big)\cap U_e$ that differ by a sign. Moreover, since $g_e$  vanishes at $e$ as $(z-e)^{3/2}$, $g_e^{2/3}$ is conformal in $U_e$. Furthermore,
\[
\left\{
\begin{array}{rcl}
g_e^{2/3}\big(\{x:~\sgn(e)x<b\}\cap U_e\big) &\subset& \big\{z:z>0\big\}, \medskip \\
g_e^{2/3}\big(\big((-1-b)\cup(b,1)\big)\cap U_e\big) &\subset& \big\{z:z<0\big\}, \medskip \\
g_e^{2/3}\big(\Delta_{2\pm}\cap U_e\big) &\subset& \big\{z:~-\sgn(e)\arg(z)=\pm2\pi/3\big\},
\end{array}
\right.
\]
where the first property follows from \eqref{arrange1} and \eqref{arrange2}, the second is a consequence of the fact that $g_e$ has purely imaginary traces there, and the third is a requirement we impose on the arcs $\Delta_{2\pm}$. Choosing the branch of $g_e^{1/6}$ which is positive on $\{x:~\sgn(e)x<b\}\cap U_e$, one has on $\big((-1,-b)\cup(b,1)\big)\cap U_e$ that
\begin{equation}
\label{geb-root}
g_{e+}^{1/6} = -\sgn(e)\mathrm{i}g_{e-}^{1/6}.
\end{equation}

As in the previous cases, one can verify that the solution of \hyperref[rhp]{\rhp$_e$} is given by
\begin{equation}
\label{sol-b}
\boldsymbol P_e := \boldsymbol E_e\mathsf{T}_2\boldsymbol\Phi_e\left((3n/2)^{2/3}g_e^{2/3}\right)\boldsymbol W_e,
\end{equation}
where $\boldsymbol\Phi_{-b}:=\boldsymbol\Phi$ and $\boldsymbol\Phi_b:=\sigma_3\boldsymbol\Phi\sigma_3$, $\boldsymbol W_e$ is defined by \eqref{W_e}, and
\begin{equation}
\label{aux-b}
\boldsymbol E_e := \boldsymbol M\mathsf{T}_2\left(\rho_2^{*\sigma_3/2}\left(3ng_e/2\right)^{\sigma_3/6} \frac1{\sqrt2} \left(\begin{matrix} 1 & -\sgn(e)\mathrm{i} \\ -\sgn(e)\mathrm{i} & 1 \end{matrix}\right)\right),
\end{equation}
whose holomorphy can be checked as in the previous cases using \eqref{geb-root}.

\subsubsection{\hyperref[rhp]{\rhp$_{\pm \mathrm{i}c}$}}

To map $\boldsymbol\Phi$ into $U_e$, $e\in\{\pm\mathrm{i}c\}$, set $\sgn(\pm\mathrm{i}c)=\pm$ and define
\[
g_e(z) := -\frac12\int_e^z\left(h_1-h_2\right)(t)\mathrm{d}t = -\frac12\log\left(\Phi^{(1)}(z)/\Phi^{(2)}(z)\right),
\]
$z\in U_e\setminus\Gamma$. Since $\Gamma$ is the branch cut for $h_1$ and $h_2$, the traces of $g_e$ on $\Gamma\cap U_e$ differ by a sign. Moreover, since $g_e$  vanishes at $e$ as $(z-e)^{3/2}$, $g_e^{2/3}$ is conformal in $U_e$. The following are conditions we impose on the arcs $\Gamma$, $\Gamma_2$, and $\Gamma_{1\pm}$:
\[
\left\{
\begin{array}{rcl}
g_e^{2/3}\big(\Gamma_2\cap U_e\big) &\subset& \big\{z:z>0\big\}, \medskip \\
g_e^{2/3}\big(\Gamma\cap U_e\big) &\subset& \big\{z:z<0\big\}, \medskip \\
g_e^{2/3}\big(\Gamma_{1\pm}\cap U_e\big) &\subset& \big\{z:~\sgn(e)\arg(z)=\pm2\pi/3\big\}.
\end{array}
\right.
\]
Choosing the branch of $g_e^{1/6}$ which is positive on $\Gamma_2\cap U_e$, we see that \eqref{geb-root} holds on $\Gamma\cap U_e$ with $-\sgn(e)$ replaced by $\sgn(e)$. Then the solution of \hyperref[rhp]{\rhp$_e$} is given by
\begin{equation}
\label{sol-ic}
\boldsymbol P_e := \boldsymbol E_e\mathsf{T}_1\boldsymbol\Phi_e\left((3n/2)^{2/3}g_e^{2/3}\right)\boldsymbol W_e,
\end{equation}
where $\boldsymbol\Phi_{\mathrm{i}c}:=\boldsymbol\Phi$ and $\boldsymbol\Phi_{-\mathrm{i}c}:=\sigma_3\boldsymbol\Phi\sigma_3$, and
\begin{equation}
\label{aux-ic}
\left\{
\begin{array}{lll}
\boldsymbol W_e & := & \diag\left( \big(\Phi^{(0)}\big)^n,  \big(\Phi^{(1)}\Phi^{(2)}\big)^{n/2}\sqrt{\rho_1/\rho_2},
\big(\Phi^{(1)}\Phi^{(2)}\big)^{n/2}\sqrt{\rho_2/\rho_1} \right), \medskip \\
\boldsymbol E_e & := & \boldsymbol M\mathsf{T}_1\left((\rho_2/\rho_1)^{\sigma_3/2}\left(3ng_e/2\right)^{\sigma_3/6} \frac1{\sqrt2}\left(\begin{matrix} 1 & -\sgn(e)\mathrm{i} \\ -\sgn(e)\mathrm{i} & 1 \end{matrix}\right)\right).
\end{array}
\right.
\end{equation}

\subsection{Final R-H Problem}
\label{ssec:FR-H}

\begin{figure}[!ht]
\centering
\includegraphics[scale=.5]{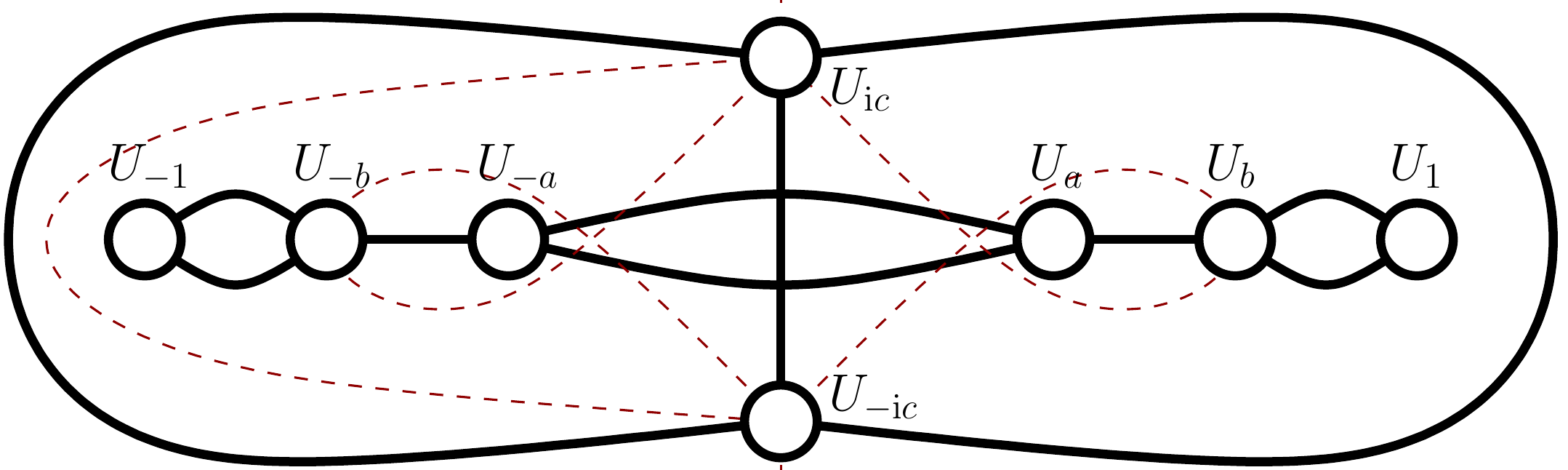}
\caption{\small The contour $\Sigma_{\boldsymbol R}$: solid lines. The dashed lines represent the relevant borders of the domains $\Omega_{ijk}$, see Figure~\ref{Oms}(a).}
\label{fig:Sigma_R}
\end{figure}

Consider the following Riemann-Hilbert problem (\rhr):
\begin{itemize}
\label{rhr}
\item[(a)] $\boldsymbol R$ is analytic in $\overline\C\setminus\Sigma_{\boldsymbol R}$, where $\Sigma_{\boldsymbol R}$ is the contour shown in Figure~\ref{fig:Sigma_R}, and $\boldsymbol R(\infty)=\boldsymbol I$;
\item[(b)] $\boldsymbol R$ has continuous traces on each side of $\Sigma_{\boldsymbol R}^\circ$ that satisfy
\[
{\boldsymbol R}_+ = {\boldsymbol R}_-
\left\{
\begin{array}{lll}
\big(\boldsymbol {MD}^n\big){\boldsymbol J}_{\boldsymbol Z}\big(\boldsymbol {MD}^n\big)^{-1} & \text{on} & \Sigma_{\boldsymbol R}^\circ\cap\Sigma_{\boldsymbol Z}, \medskip \\
\boldsymbol P_e\big(\boldsymbol {MD}^n\big)^{-1} & \text{on} & \partial U_e, \quad e\in\left\{\pm1,\pm a,\pm b,\pm\mathrm{i}c\right\},
\end{array}
\right.
\]
where ${\boldsymbol J}_{\boldsymbol Z}$ was defined in \eqref{J-Z}, while $\boldsymbol M$ and $\boldsymbol D$ were introduced after \eqref{N}.
\end{itemize}

Let us prove that the jumps of $\boldsymbol R$ are uniformly close to $\boldsymbol I$ for $n\in\N_*$. In Section~\ref{ssec:GP-I} we have shown that $\det(\boldsymbol N)=$, while $\det(\boldsymbol D)=1$ by \eqref{normal} and $\det(\boldsymbol C)=1$ by \eqref{Phi_jatInfty}. Hence, $\det(\boldsymbol M)=1$ and therefore $\boldsymbol M^{-1}=\mathsf{adj}(\boldsymbol M)$, where $\mathsf{adj}(\boldsymbol M)$ is the adjoint matrix of $\boldsymbol M$. Thus,
\begin{equation}
\label{small-jump}
\boldsymbol M\big(\boldsymbol I+\boldsymbol{\mathcal{O}}(\cdot)\big)\boldsymbol M^{-1} = \boldsymbol I+ \boldsymbol M\boldsymbol{\mathcal{O}}(\cdot)\mathsf{adj}(\boldsymbol M) = \boldsymbol I + \boldsymbol{\mathcal{O}}(\cdot)
\end{equation}
uniformly away from $\{\pm1,\pm a,\pm b,\pm\mathrm{i}c\}$ by \eqref{Psi-Estimates}. Therefore, it follows from \hyperref[rhp]{\rhp$_e$}(d) that
\begin{equation}
\label{small-jump-1}
\boldsymbol P_e\big(\boldsymbol {MD}^n\big)^{-1} = \boldsymbol I+\boldsymbol M\boldsymbol{\mathcal{O}}(1/n)\boldsymbol M^{-1} = \boldsymbol I +\boldsymbol{\mathcal{O}}(1/n)
\end{equation}
uniformly on each $\partial U_e$. Furthermore,  we can write
\[
{\boldsymbol J}_{\boldsymbol Z} = \boldsymbol I + \rho_{j,k}\boldsymbol E_{j,k}, \quad j\neq k, \quad j,k\in\{0,1,2\},
\]
where $\boldsymbol E_{j,k}$ is the matrix with all zero entries except for the $(j+1,k+1)$-st one, which is 1, and $\rho_{j,k}$ is always a combination of $\rho_1$ and $\rho_2$ (particular values of $j,k$ and the value of the entry depend on the arc in question). Thus,
\[
\boldsymbol D^n{\boldsymbol J}_{\boldsymbol Z}\boldsymbol D^{-n} = \boldsymbol I + \left(\frac{\Phi^{(j)}}{\Phi^{(k)}}\right)^n\rho_{j,k}\boldsymbol E_{j,k} = \boldsymbol I + \boldsymbol{\mathcal{O}}\big(C^{-n}_{\boldsymbol R}\big)
\]
uniformly on $\Sigma_{\boldsymbol Z}\cap\Sigma_{\boldsymbol R}$ for some constant $C_{\boldsymbol R}>1$ by Theorem~\ref{thm:Omegas} (it is a simple examination of the five relevant cases). Therefore, we get from \eqref{small-jump} that
\begin{equation}
\label{small-jump-2}
\big(\boldsymbol {MD}^n\big){\boldsymbol J}_{\boldsymbol Z}\big(\boldsymbol M\boldsymbol D^n\big)^{-1} = \boldsymbol I + \boldsymbol{\mathcal{O}}\big(C^{-n}_{\boldsymbol R}\big)
\end{equation}
uniformly on $\Sigma_{\boldsymbol Z}\cap\Sigma_{\boldsymbol R}$. The relations \eqref{small-jump-1} and \eqref{small-jump-2} together with \cite[Corollary~7.108]{Deift} imply that \hyperref[rhr]{\rhr} is solvable for all $n\in\N_*$ large enough and satisfies
\begin{equation}
\label{R-asymp}
\boldsymbol R = \boldsymbol I + \boldsymbol{\mathcal{O}}(1/n), \quad \boldsymbol R(\infty) = \boldsymbol I,
\end{equation}
uniformly in $\overline\C$, that is, including the boundary values on $\Sigma_{\boldsymbol R}$.

\subsection{Asymptotics of Hermite-Pad\'e Approximants}

Inverting \eqref{Z} and \eqref{S}, we get from \eqref{Y} that
\begin{equation}
\label{Qn-Final}
Q_n(z) = [\boldsymbol Z]_{11}(z) + \frac1{\rho_i^*(z)}\left\{
\begin{array}{rl}
\pm[\boldsymbol Z]_{1i+1}(z), &  z\in O_{i\pm}, \medskip \\
0, & \text{otherwise},
\end{array}
\right.
\end{equation}
where $\rho^*_1=\rho_1$;
\begin{equation}
\label{Rn1-Final}
R_{\vec n}^{(1)}(z) = [\boldsymbol Z]_{12}(z) - \frac{\rho_1(z)}{\rho_2(z)}\left\{
\begin{array}{rl}
-[\boldsymbol Z]_{13}(z), & z\in O_0, \medskip \\
~[\boldsymbol Z]_{13}(z), & z\in O_1, \medskip \\
0, & \text{otherwise},
\end{array}
\right.
\end{equation}
for $z\notin F_1$; and
\begin{equation}
\label{Rn2-Final}
R_{\vec n}^{(2)}(z) = \left\{
\begin{array}{rl}
0, & z\in O_1 \medskip \\
~[\boldsymbol Z]_{13}(z), & \text{otherwise}
\end{array}
\right\} + \frac{\rho_2(z)}{\rho_1(z)}\left\{
\begin{array}{rl}
~[\boldsymbol Z]_{12}(z), & O_1\cup O_2 \medskip \\
0, & \text{otherwise}
\end{array}
\right\}
\end{equation}
for $z\notin F_2$.

Let $\boldsymbol R$ be the solution of \hyperref[rhr]{\rhr}. Then
\begin{equation}
\label{Z-solution}
\boldsymbol Z = \boldsymbol C^{-n}\boldsymbol R
\left\{
\begin{array}{rll}
\boldsymbol M\boldsymbol D^n & \text{in} & \C\setminus\bigcup_e U_e, \medskip \\
\boldsymbol P_e & \text{in} & U_e, \quad e\in\{\pm1,\pm a,\pm b,\pm\mathrm{i}c\},
\end{array}
\right.
\end{equation}
solves \hyperref[rhz]{\rhz} for all $n\in\N_*$ large enough. Denote the first row of $\boldsymbol R$ by $(\begin{matrix} 1+\upsilon_{n,0} & \upsilon_{n,1} & \upsilon_{n,2} \end{matrix})$. It follows from \eqref{R-asymp} that
\begin{equation}
\label{upsilons}
|\upsilon_{n,k}(\infty)|=0 \quad \text{and} \quad |\upsilon_{n,k}(z)|=\mathcal{O}(1/n) \quad \text{uniformly for} \quad z\in\C,
\end{equation}
meaning that $\upsilon_{n,k}(z)$ is replaced by $\upsilon_{n,k}^\pm(z)$ when $z\in\Sigma_{\boldsymbol R}$. If $z\not\in U_e$, $e\in\{\pm1,\pm a,\pm b,\pm\mathrm{i}c\}$, then
\begin{equation}
\label{Z-Exact}
[\boldsymbol Z]_{1k+1}(z) = C_n\Psi_{n,0}^{(k)}(z) + C_n\sum_{j=0}^2\upsilon_{n,j}(z)\Psi_{n,j}^{(k)}(z), \quad k\in\{0,1,2\},
\end{equation}
according to \eqref{Cn}, where the functions are replaced by their traces when necessary.

Formula \eqref{Z-Exact} is valid for $z\in U_{\mathrm{i}c}\cup U_{-\mathrm{i}c}$ when $k=0$ and for $z\in U_e$, $e\in\{\pm b,\pm1\}$, when $k=1$. Indeed, in these cases $\boldsymbol P_e=\boldsymbol M\mathsf{T}_{k+1}(\cdot)\boldsymbol W_e$, see \eqref{sol-ic} and \eqref{aux-ic}, \eqref{sol-1} and \eqref{E_e}, \eqref{sol-b} and \eqref{aux-b}. As the $(k+1)$-st column of $\boldsymbol W_e$ is the same as the $(k+1)$-st column of $\boldsymbol D^n$, the $(k+1)$-st column of $\boldsymbol P_e$ is the same as the $(k+1)$-st column of $\boldsymbol{MD}^n$ from which the claim follows.

To estimate the sum on the right-hand side of \eqref{Z-Exact}, recall that
\begin{equation}
\label{nj-n0}
\Psi_{n,j} = \Psi_{n,0}\Theta_j\Theta_{n,j}\Theta_{n,0}^{-1},
\end{equation}
see \eqref{Psink}, where $\Theta_j\Theta_{n,j}\Theta_{n,0}^{-1}$ is a rational function on $\RS$ with the divisor $\mathcal{D}_{n,j}+\infty^{(0)}-\mathcal{D}_{n,0}-\infty^{(j)}$. Then \eqref{upsilons}, the maximum modulus principle, and the same normal family argument as in \eqref{Psi-Est1} and \eqref{Psi-Est2} imply that
\begin{equation}
\label{atlast}
\big|\upsilon_{n,j}(z)\Theta_j\big(z^{(k)}\big)\Theta_{n,j}\big(z^{(k)}\big)\Theta_{n,0}^{-1}\big(z^{(k)}\big)\big| = \mathcal{O}\big(1/n;\mathcal{D}_{n,0}\cap\RS^{(k)}\big)
\end{equation}
uniformly in $\overline\C$, where the function on the left-hand side of \eqref{atlast} needs to be replaced by its traces when $z\in\Sigma_{\boldsymbol R}\cup\pi\big(\partial\RS^{(k)}\big)$. By combining \eqref{nj-n0} with \eqref{atlast} we get that
\begin{equation}
\label{bootstrap}
\upsilon_{n,j}(z)\Psi_{n,j}^{(k)}(z) = \Psi_{n,0}^{(k)}(z)\mathcal{O}\big(1/n;\mathcal{D}_{n,0}\cap\RS^{(k)}\big), \quad j,k\in\{0,1,2\},
\end{equation}
where $\mathcal{O}(\cdot)$ is uniform in $\overline\C$.

The first relation in \eqref{main-A} follows immediately from the first line of \eqref{Qn-Final}, \eqref{Z-Exact}, and \eqref{bootstrap}. Moreover, \eqref{Z-Exact} implies that the second line of \eqref{Qn-Final} can be rewritten as
\[
Q_{\vec n}(z) = C_n \left(\Psi_{n,0+}^{(0)}(z) + \Psi_{n,0-}^{(0)}(z) \right) + C_n\sum_{j=0}^2\upsilon_{n,j}(z)\left(\Psi_{n,j+}^{(0)}(z) + \Psi_{n,j-}^{(0)}(z) \right)
\]
for $z\in O_{1+}\cup O_{1-}\cup O_{2+}\cup O_{2-}$ and $z\not\in\overline U_e$, $e\in\{\pm b,\pm a,\pm 1\}$, where
\begin{equation*}
\Psi_{n,j\pm}^{(0)}(z) = \left\{
\begin{array}{ll}
\Psi_{n,j}^{(0)}(z), & z\in O_{i\pm}, \medskip \\
\mp\Psi_{n,j}^{(i)}(z)/\rho_i^*(z), & z\in O_{i\mp},
\end{array}
\right. \quad i\in\{1,2\}.
\end{equation*}
Clearly, each $\Psi_{n,j\pm}^{(0)}$ extends to a holomorphic function in $O_{i+}\cup O_{i-}\cup \Delta_i^\circ$ by \eqref{BVP} (recall that $\Psi_{n,1}$ and $\Psi_{n,2}$ also satisfy these relations). The first part of \eqref{main-B} now follows from \eqref{bootstrap}.

Furthermore, \eqref{Rn1-Final}, \eqref{Z-Exact}, and \eqref{bootstrap} imply the second line of \eqref{main-A} outside of $O_0\cup O_1$ for $i=1$. In the spirit of \eqref{PsiHat}, define
\[
\widehat\Psi_{n,j+}^{(1)}(z) = \left\{
\begin{array}{ll}
\Psi_{n,j}^{(1)}(z), & z\in O_0, \medskip \\
\Psi_{n,j}^{(2)}(z)(-\rho_1/\rho_2)(z), & z\in O_1,
\end{array}
\right.
\]
which is holomorphic in $(O_0\cup O_1\cup \Delta_0^\circ)\setminus\Delta_{21}$ by \eqref{BVP}, and
\[
\widehat\Psi_{n,j-}^{(1)}(z) = \left\{
\begin{array}{ll}
\Psi_{n,j}^{(1)}(z), & z\in O_1,  \medskip \\
\Psi_{n,j}^{(2)}(z)(\rho_1/\rho_2)(z), & z\in O_0,
\end{array}
\right.
\]
which is holomorphic in $(O_0\cup O_1\cup\Delta_0^\circ)\setminus\Delta_1$, again, by \eqref{BVP}. Then \eqref{Rn1-Final} and \eqref{Z-Exact} imply that
\[
R_{\vec n}^{(1)}(z) =  C_n\left( \widehat\Psi_{n,0+}^{(1)}(z) + \widehat\Psi_{n,0-}^{(1)}(z) \right)  + C_n\sum_{j=0}^2\upsilon_{n,j}(z)\left( \widehat\Psi_{n,j+}^{(1)}(z) + \widehat\Psi_{n,j-}^{(1)}(z) \right).
\]
for $z\in (O_0\cup O_1\cup\Delta_0^\circ)\setminus F_1$ and $z\not\in U_e$, $e\in\{-1,-b,-a,\pm\mathrm{i}c\}$. Clearly, formula \eqref{main-B} is an immediate consequence of \eqref{bootstrap}. To finish the proof of \eqref{main-A} for $i=1$, let us show that
\begin{equation*}
\widehat\Psi_{n,j+}^{(1)} = \widehat\Psi_{n,0-}^{(1)}\mathcal{O}\big(1/n;\mathcal{D}_{n,0}\cap\RS^{(1)}\big)
\end{equation*}
locally uniformly in $\big(\mathsf{int}(\Gamma_1\cup\Delta_{01-})\cup \Delta_{01-}\big)\setminus \big(F_1\cup \overline U_{-\mathrm{i}c}\cup \overline U_{\mathrm{ic}}\big)$. Indeed, we have that
\[
\frac{\widehat\Psi_{n,j+}^{(1)}}{\widehat\Psi_{n,0-}^{(1)}} = -\frac{\rho_1w_2}{\rho_2w_1}\left(\frac{\Phi^{(2)}}{\Phi^{(1)}}\right)^n\frac{S_\rho^{(2)}S_{\vec y_n}^{(2)}\Theta_{n,j}^{(2)}\Theta_j^{(2)}}{S_\rho^{(1)}S_{\vec y_n}^{(1)}\Theta_{n,0}^{(1)}}
\]
see \eqref{Psink}. The claim now follows from Theorem~\ref{thm:Omegas} and the normal family argument along the lines of \eqref{Psi-Est1} and \eqref{Psi-Est2}. This proves \eqref{main-A} outside of $\mathsf{int}(\Gamma_1\cup\Delta_{01+})$. Analogously we can argue that
\begin{equation*}
\widehat\Psi_{n,j-}^{(1)} = \widehat\Psi_{n,0+}^{(1)}\mathcal{O}\big(1/n;\mathcal{D}_{n,0}\cap\RS^{(2)}\big)
\end{equation*}
locally uniformly in $\big(\mathsf{int}(\Gamma_1\cup\Delta_{01+})\cup \Delta_{01+}\big)\setminus \big(F_1\cup \overline U_{-\mathrm{i}c}\cup \overline U_{\mathrm{ic}}\big)$, which finishes the proof of \eqref{main-A} for $R_{\vec n}^{(1)}$.  The proof of \eqref{main-A} and \eqref{main-B} for $R_{\vec n}^{(2)}$ can be completed analogously starting with \eqref{Rn2-Final}.

\section{Riemann-Hilbert Analysis: Case II}
\label{sec:CaseII}

\subsection{Global Lenses}

\begin{figure}[!ht]
\centering
\includegraphics[scale=.5]{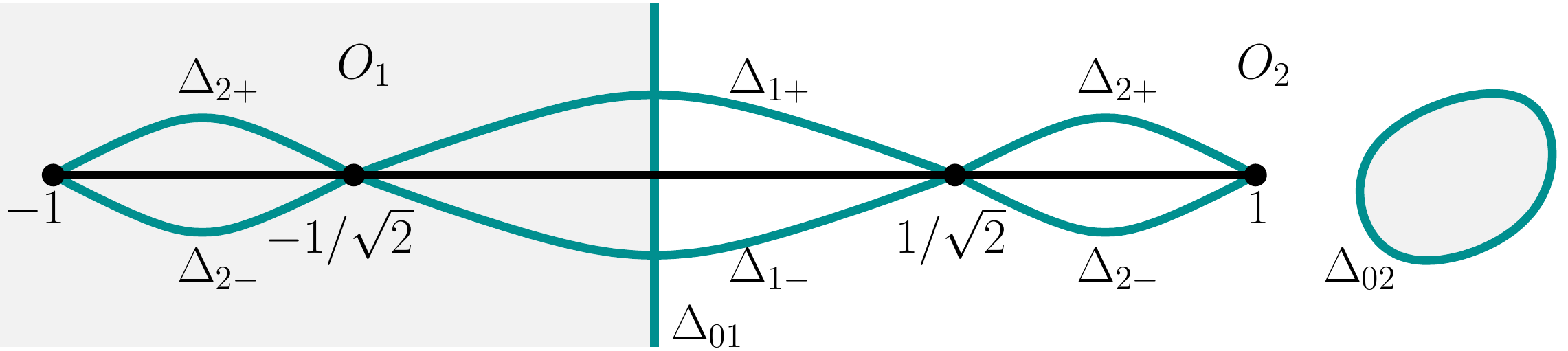}
\caption{\small The lens $\Sigma_{\boldsymbol Z}$ without the circle $\{|z|=R\}$, the domains $O_1$ (shaded regions on the left) and $O_2$ (unshaded region), local lenses $\Delta_{1\pm}$ and $\Delta_{2\pm}$.}
\label{fig:Global-II}
\end{figure}

Let $\boldsymbol G_1(u)$ and $\boldsymbol G_2(v)$ be defined by \eqref{global-matr}. Further, let $O_1:=\{z:~\re(z)<0\}$ be the left half-plane. We orient the boundary of $O_1$, say $\Delta_{01}$ (the imaginary axis), so that $O_1$ lies to the left of $\Delta_{01}$ when the latter is traversed in the positive direction. Denote by $\Delta_{02}$ a simple Jordan curve lying within the right component of $\Omega_{021}$ containing all the singularities of $\rho_2/\rho_1$, see Figure~\ref{fig:Global-II}
(this is always possible because of Condition~\ref{cond}). We orient $\Delta_{02}$ counter-clockwise and set $O_2$ to be the intersection of the exterior domain of $\Delta_{02}$ and the right half-plane $\{\re(z)>0\}$. Put
\begin{equation}
\label{S2}
\boldsymbol S := \boldsymbol G_1(-\varrho/2)\boldsymbol G_2(1/\varrho) \boldsymbol Y\left\{
\begin{array}{lcl}
\boldsymbol G_2(-\rho_2/\rho_1)\boldsymbol G_1(\rho_1/\rho_2) & \text{in} & O_1, \medskip \\
\boldsymbol G_2(-\rho_2/\rho_1) & \text{in} & O_2, \medskip \\
\boldsymbol I & \text{in} & \C\setminus\left(\overline{O_1\cup O_2}\right),
\end{array}
\right.
\end{equation}
where $\varrho:=\rho_1(\infty)/\rho_2(\infty)$. Put $\Sigma_{\boldsymbol S}:=[-1,1]\cup\Delta_{01}\cup\Delta_{02}$. Then $\boldsymbol S$ solves \rhs:
\begin{itemize}
\label{rhs2}
\item[(a)] $\boldsymbol S$ is analytic in $\C\setminus\Sigma_{\boldsymbol S}$ and $\lim_{z\to\infty,\pm\re(z)>0} \boldsymbol S(z) \diag\left(z^{-2n},z^n,z^n\right) = \boldsymbol G_1(\mp\varrho/2)$;
\item[(b)] ${\boldsymbol S}$ has continuous traces on $\Sigma_{\boldsymbol S}^\circ:=\Sigma_{\boldsymbol S}\setminus\big\{\pm1,\pm1/\sqrt2,0\big\}$ that satisfy
\[
{\boldsymbol S}_+ = {\boldsymbol S}_-
\left\{
\begin{array}{lll}
{\boldsymbol J}(\rho_1,0) & \text{on} & \Delta_1\setminus\{0\}, \medskip \\
{\boldsymbol J}(0,\rho_2^*) & \text{on} & \Delta_2, \medskip \\
{\boldsymbol G}_1(\rho_1/\rho_2) & \text{on} & \Delta_{01}, \medskip \\
{\boldsymbol G}_2(\rho_2/\rho_1) & \text{on} & \Delta_{02},
\end{array}
\right.
\]
where $\rho_2^*$ is defined by \eqref{rhostar} and $\Delta_i$ in \eqref{chains};
\item[(c)] $\boldsymbol S$ satisfies \hyperref[rhy]{\rhy}(c) (see Section \ref{MOP}) with $[-1,1]$ replaced by $\Sigma_{\boldsymbol S}$.
\end{itemize}
If  \hyperref[rhs]{\rhs} is solvable, then so is  \hyperref[rhy]{\rhy}, and the solutions are connected via \eqref{S2}.

\subsection{Local Lenses}

As in Case I, we introduce additional arcs $\Delta_{1\pm}$ and systems of two arcs $\Delta_{2\pm}$ as in Figure~\ref{fig:Global-II}, all oriented from left to right.  We further denote by $O_{i\pm}$ the open sets bounded by $\Delta_i$ and the arcs $\Delta_{i\pm}$, $i\in\{1,2\}$. Set
\begin{equation}
\label{Z2}
\boldsymbol Z:= \boldsymbol S\boldsymbol L_i^{\mp1} \quad \text{in} \quad O_{i\pm},
\end{equation}
where the matrices $\boldsymbol L_i$ are defined by \eqref{local-matr}. Put $\Sigma_{\boldsymbol Z} := \Sigma_{\boldsymbol S}\cup\Delta_{1+}\cup\Delta_{1-}\cup\Delta_{2+}\cup\Delta_{2-}$. Then it can readily be checked that ${\boldsymbol Z}$ solves \rhz:
\begin{itemize}
\label{rhz2}
\item[(a)] ${\boldsymbol Z}$ is analytic in $\C\setminus\Sigma_{\boldsymbol Z}$ and $\lim_{z\to\infty,\pm\re(z)>0} \boldsymbol Z(z) \diag\left(z^{-2n},z^n,z^n\right) = \boldsymbol G_1(\mp\varrho/2)$;
\item[(b)] ${\boldsymbol Z}$ has continuous traces on each side of $\Sigma_{\boldsymbol Z}^\circ:=\Sigma_{\boldsymbol Z}\setminus\big\{\pm1,\pm 1/\sqrt2,0\big\}$ that satisfy
\[
{\boldsymbol Z}_+ = {\boldsymbol Z}_-
\left\{
\begin{array}{lll}
{\boldsymbol J}_i & \text{on} & \Delta_i^\circ\setminus\{0\}, \quad i\in\{1,2\}, \medskip \\
\boldsymbol J_{\boldsymbol Z} &\text{on} & \bigcup_{i=1}^2(\Delta_{0i}\cup\Delta_{i+}\cup\Delta_{i-}),
\end{array}
\right.
\]
where $\boldsymbol J_i$ are defined by \eqref{local-jumps} and
\begin{equation}
\label{J-Z2}
\boldsymbol J_{\boldsymbol Z} :=
\left\{
\begin{array}{lll}
{\boldsymbol G}_1(\rho_1/\rho_2) & \text{on} & \Delta_{01}, \medskip \\
{\boldsymbol G}_2(\rho_2/\rho_1) & \text{on} & \Delta_{02}, \medskip \\
{\boldsymbol L}_i & \text{on} & \Delta_{i\pm}, \quad i\in\{1,2\};
\end{array}
\right.
\end{equation}
\item[(c)] $\boldsymbol Z$ satisfies \hyperref[rhs3]{\rhs}(c) with $\Sigma_{\boldsymbol S}$ replaced by $\Sigma_{\boldsymbol Z}$.
\end{itemize}
As before, if \hyperref[rhz3]{\rhz} is solvable, then so is \hyperref[rhs3]{\rhs}, and the solutions are connected via \eqref{Z2}.

\subsection{Global Parametrix}

Let $\boldsymbol N=\boldsymbol C^{-n}\boldsymbol {MD}^n$ be given by \eqref{N}. Then it is a solution of the following Riemann-Hilbert problem (\rhn):
\begin{itemize}
\label{rhn2}
\item[(a)] ${\boldsymbol N}$ is analytic in $\C\setminus[-1,1]$ and $\lim_{z\to\infty} {\boldsymbol N}(z)\diag\left(z^{-2n},z^n,z^n\right) = \boldsymbol I$;
\item[(b)] ${\boldsymbol N}$ has continuous traces on each side of $\Delta_i^\circ$, $i\in\{1,2\}$, that satisfy  $\boldsymbol N_+ = \boldsymbol N_- \boldsymbol J_i$.
\end{itemize}

We cannot argue that $\det(\boldsymbol M)\equiv1$ as in Section~\ref{ssec:GP-I} since every entry of $\boldsymbol M$ behaves like $(z\mp1/\sqrt2)^{-1/3}$ as $z\to\pm1/\sqrt2$. However, we can construct $\boldsymbol M^{-1}$ explicitly. Denote by $\widetilde{\boldsymbol N} =\boldsymbol C^{-n}\widetilde{\boldsymbol M}\boldsymbol D^n$ the matrix that satisfies \hyperref[rhn2]{\rhn} as above with $\rho_1$ and $\rho_2^*$ replaced by $1/\rho_1$ and $1/\rho_2^*$. It follows from \eqref{simult} that the construction \eqref{N} of the matrices $\boldsymbol N$ and $\widetilde{\boldsymbol N}$ is simultaneously applicable or not applicable for each index $n$. Observe that
\[
\left(\boldsymbol M\widetilde{\boldsymbol M}^\mathsf{T}\right)_-^{-1}\left(\boldsymbol M\widetilde{\boldsymbol M}^\mathsf{T}\right)_+ = \left(\widetilde{\boldsymbol M}_+^\mathsf{T}\right)^{-1}\left(\boldsymbol D_-/\boldsymbol D_+\right)^{2n}\widetilde{\boldsymbol M}^\mathsf{T}_+ = \left(\widetilde{\boldsymbol M}_+^\mathsf{T}\right)^{-1}\widetilde{\boldsymbol M}^\mathsf{T}_+ = \boldsymbol I
\]
on $\pi(\ualpha_1)$ and $\pi(\ubeta_1)$ by \eqref{Phi-jumps} and~\eqref{not-rational2}. Moreover,
\[
\left(\boldsymbol M\widetilde{\boldsymbol M}^\mathsf{T}\right)_-^{-1}\left(\boldsymbol M\widetilde{\boldsymbol M}^\mathsf{T}\right)_+ = \left(\widetilde{\boldsymbol M}_-^\mathsf{T}\right)^{-1}\boldsymbol J_i\widetilde{\boldsymbol M}^\mathsf{T}_+ = \left(\widetilde{\boldsymbol M}_+\boldsymbol J_i^\mathsf{T}\widetilde{\boldsymbol M}_-^{-1}\right)^\mathsf{T}= \boldsymbol I
\]
on $\Delta_i^\circ$ since $\boldsymbol J_i^\mathsf{T}=\widetilde{\boldsymbol J}_i^{-1}$. It also follows from \eqref{around-branching} that the entries of $\boldsymbol M$ and $\widetilde{\boldsymbol M}$ have at most $1/4$ and $1/3$-root singularities at $\pm1$ and $\pm1/\sqrt2$, respectively (the functions $\Psi_{n,1}$ and $\Psi_{n,2}$ possess exactly the same behavior around those points as $\Psi_n=\Psi_{n,0}$). Hence, the product $\boldsymbol M\widetilde{\boldsymbol M}^\mathsf{T}$ is holomorphic in the entire complex plane. Since $\boldsymbol M(\infty)=\widetilde{\boldsymbol M}(\infty)=\boldsymbol I$, we deduce that
\begin{equation}
\label{Minverse}
\boldsymbol I = \widetilde{\boldsymbol M}^\mathsf{T}\boldsymbol M.
\end{equation}
In particular, it follows from \eqref{around-branching} that
\begin{equation}
\label{M-local}
\left\{
\begin{array}{lll}
\boldsymbol M(z) &=& \big(z\mp1/\sqrt2\big)^{-1/3}\boldsymbol M_\pm + \boldsymbol{\mathcal O}(1), \medskip \\
\widetilde{\boldsymbol M}^\mathsf{T}(z) &=& \big(z\mp1/\sqrt2\big)^{-1/3}\widetilde{\boldsymbol M}_\pm^\mathsf{T} + \boldsymbol{\mathcal O}(1),
\end{array}
\quad \text{as} \quad z\to\pm1/\sqrt2,
\right.
\end{equation}
for some constant matrices $\boldsymbol M_\pm$ and $\widetilde{\boldsymbol M}_\pm^\mathsf{T}$. Then it follows immediately from \eqref{Minverse} that
\begin{equation}
\label{zero-product}
\widetilde{\boldsymbol M}_\pm^\mathsf{T} \boldsymbol M_\pm = \boldsymbol 0.
\end{equation}

\subsection{Local Parametrices}

Again, we need to solve \hyperref[rhz2]{\rhz} locally, this time around $\pm1$, $\pm1/\sqrt2$, and $\infty$. The local problems \hyperref[rhp]{\rhp$_{\pm1}$} are exactly the same as in Case I and therefore their solutions are given by \eqref{W_e}--\eqref{E_e}.

\subsubsection{\hyperref[rhpi]{\rhp$_\infty$}}

Let $R>0$ be large enough so that $\Delta_{02}\subset \big\{|z|<R\big\}$. In this section, we are looking for a solution to the following Riemann-Hilbert problem (\rhp$_\infty$):
\begin{itemize}
\label{rhpi}
\item[(a)] $\boldsymbol P_\infty$ is holomorphic in $\{|z|>R\}\setminus\Delta_{01}$ and
\[
\lim_{z\to\infty,\pm\re(z)>0} \boldsymbol C^{-n}\boldsymbol P_\infty(z) \diag\left(z^{-2n},z^n,z^n\right) = \boldsymbol G_1(\mp\varrho/2);
\]
\item[(b)] $\boldsymbol P_\infty$ has continuous traces on each side of $\Delta_{01}\cap \{|z|>R\}$ that satisfy $\boldsymbol P_{\infty+}=\boldsymbol P_{\infty-}\boldsymbol G_1(\rho_1/\rho_2)$;
\item[(c)] $\boldsymbol P_\infty=\boldsymbol M\left(\boldsymbol I+\boldsymbol{\mathcal{O}}\big(n^{-1/2}\big)\right)\boldsymbol D^n$ uniformly on $\{|z|=R\}$.
\end{itemize}

Let us show that \hyperref[rhpi]{\rhp$_\infty$} is solved by
\begin{equation}
\label{Pinfty}
\boldsymbol P_\infty := \boldsymbol M\boldsymbol G_1\left((\rho_1/\rho_2)(C_2/C_1)^nu\big(\sqrt{n/2}\zeta\big)\right)\boldsymbol D^n,
\end{equation}
where the constants $C_k$ are defined in \eqref{Phi_jatInfty},  the function $u$ is given by
\[
u(\zeta) := \frac12e^{2\zeta^2}\left\{
\begin{array}{rl}
\mathrm{erfc}\big(-\sqrt 2\zeta\big), & \re(\zeta)<0, \medskip \\
-\mathrm{erfc}\big(\sqrt 2\zeta\big), & \re(\zeta)>0,
\end{array}
\right.
\]
and $\zeta$ is defined by
\[
\zeta(z) := \mathrm{i}\sqrt{\log\left(\Phi^{(1)}(z)C_2/\Phi^{(2)}(z)C_1\right)}, \quad |z|>R.
\]
Indeed, $\zeta(z)$ is a conformal function in $\{|z|>R\}$ that vanishes at infinity by \eqref{case2hi} (make $R$ larger if necessary). Here, we choose the branch of the square root so that $z\zeta(z)$ tends to $(12)^{-1/4}$ when $z\to\infty$. Hence, we can deform $\Delta_{01}$ in $\{|z|>R\}$ so that $\zeta(\Delta_{01})\subset\mathrm{i}\R$. Thus, the right-hand side of \eqref{Pinfty} is holomorphic in $\{|z|>R\}\setminus\Delta_{01}$. As it follows from \cite[Eq. (7.2.2)]{DLMF} that
\[
\lim_{z\to0,\pm\re(\zeta)>0} u(\zeta) = \mp1/2,
\]
\hyperref[rhpi]{\rhp$_\infty$}(a) is indeed satisfied. To verify \hyperref[rhpi]{\rhp$_\infty$}(b), notice that $\boldsymbol G_1(f_+)=\boldsymbol G_1(f_-)\boldsymbol G_1(f_+-f_-)$ and that $u_+(x) - u_-(x) = e^{2x^2}$ for $x\in\mathrm{i}\R$. Since
\[
\exp\left\{2(n/2)\zeta^2(z)\right\} = \left(\Phi^{(2)}(z)C_1/\Phi^{(1)}(z)C_2\right)^n
\]
\hyperref[rhpi]{\rhp$_\infty$}(b) follows. Finally, \cite[Eq. (7.12.1)]{DLMF} implies that
\[
u(\zeta) \sim \frac1{\sqrt{2\pi}}\sum_{k=0}^\infty(-1)^{k+1}\frac{\Gamma(k+1/2)}{2^{k+1}\Gamma(1/2)}\zeta^{-(2k+1)}
\]
uniformly in the left and right half-planes. Hence,
\[
\boldsymbol P_\infty = \boldsymbol M\left(\boldsymbol I+ (\rho_1/\rho_2)(C_2/C_1)^n\mathcal{O}\big(1/\sqrt{n}\big)\boldsymbol E_{3,2}\right)\boldsymbol D^n
\]
uniformly on $|z|=R$. As $|C_1|=|C_2|$ by Theorem~\ref{thm:Omegas}, we see that \hyperref[rhpi]{\rhp$_\infty$}(c) holds as well.

\subsubsection{\hyperref[rhp]{\rhp$_{\pm1/\sqrt2}$}}

Denote by $U_{\pm1/\sqrt2}:=\big\{z:~\big|z\mp1/\sqrt2\big|<n^{-1/2}\big\}$. In this section, we construct a solution to the following Riemann-Hilbert problem:
\begin{itemize}
\label{rhp2}
\item[(a)] $\boldsymbol P_{\pm1/\sqrt2}$ is analytic in $U_{\pm1/\sqrt2}\setminus\Sigma_{\boldsymbol Z}$;
\item[(b)] $\boldsymbol P_{\pm1/\sqrt2}$ has continuous traces on each side of $\Sigma_{\boldsymbol Z}^\circ\cap U_{\pm1/\sqrt2}$ that satisfy \hyperref[rhz]{\rhz}(b);
\item[(c)] $\boldsymbol P_{\pm1/\sqrt2}$ has the behavior near ${\pm1/\sqrt2}$ within $U_{\pm1/\sqrt2}$ described by \hyperref[rhz]{\rhz}(c);
\item[(d)] $\boldsymbol P_{\pm1/\sqrt2}=\left(\boldsymbol I+\boldsymbol{\mathcal{O}}\big(n^{-1/6}\big)\right)\boldsymbol M\boldsymbol D^n$ uniformly on $\partial U_{\pm1/\sqrt2}\setminus\Sigma_{\boldsymbol Z}$.
\end{itemize}

In \cite{DeK11} (an alternative approach to asymptotics of multiple orthogonal polynomials around cubic branch point was developed in \cite{Toul09}), a $3\times3$ matrix function was constructed out of solutions to $zy^{\prime\prime\prime}(z)-\tau y^\prime(z)+y(z)=0$ that solves \rhup:
\begin{itemize}
\label{rhup}
\item[(a)] $\boldsymbol\Upsilon$ is holomorphic in $\C\setminus\big(L_+\cup L_-\cup(-\infty,\infty)\big)$, where $L_\pm := \big\{\zeta:\re(\zeta)=\pm\im(\zeta)\big\}$ and the positive direction on all the lines goes from the left half-plane to the right half-plane;
\item[(b)] $\boldsymbol\Upsilon$ has continuous traces on $\big(L_+\cup L_-\cup(-\infty,\infty)\big)\setminus\{0\}$ that satisfy
\[
\boldsymbol\Upsilon_+ = \boldsymbol\Upsilon_-
\left\{
\begin{array}{lll}
\mathsf{T}_2\left(\begin{matrix} 1 & 0 \\ 1 & 1 \end{matrix}\right) & \text{on} & L_\pm\cap \big\{\zeta:\re(\zeta)>0\big\}, \medskip \\
\mathsf{T}_2\left(\begin{matrix} 0 & 1 \\ -1 & 0 \end{matrix}\right) & \text{on} & (0,\infty),
\end{array}
\right.
\]
and
\[
\boldsymbol\Upsilon_+ = \boldsymbol\Upsilon_-
\left\{
\begin{array}{lll}
\mathsf{T}_3\left(\begin{matrix} 1 & 0 \\ 1 & 1 \end{matrix}\right) & \text{on} & L_\pm\cap \big\{\zeta:\re(\zeta)<0\big\}, \medskip \\
\mathsf{T}_3\left(\begin{matrix} 0 & 1 \\ -1 & 0 \end{matrix}\right) & \text{on} & (-\infty,0);
\end{array}
\right.
\]
\item[(c)] $\boldsymbol\Upsilon(\zeta)=\boldsymbol{\mathcal{O}}(\log|\zeta|)$ as $\zeta\to0$;
\item[(d)] $\boldsymbol\Upsilon$ has the following behavior near $\infty$:
\[
\boldsymbol\Upsilon(\zeta;\tau) = \boldsymbol A(\zeta;\tau) \left(\boldsymbol I+\boldsymbol\Upsilon_1(\tau)\zeta^{-1/3} + \boldsymbol{\mathcal{O}}\left(\zeta^{-2/3}\right)\right)\exp\left\{ -\frac32\zeta^{2/3}\boldsymbol B^2-\tau\zeta^{1/3}\boldsymbol B\right\}
\]
uniformly for $\zeta\in\C\setminus\big(L_+\cup L_-\cup(-\infty,\infty)\big)$ and $\tau$ on bounded sets with
\[
\boldsymbol \Upsilon_1(\tau) = -\frac\tau3 \left(\frac{\tau^2}9+1\right)\boldsymbol B^2 -\frac\tau9\boldsymbol C,
\]
where
\[
\boldsymbol A(\zeta;\tau) := \sqrt{\frac{2\pi}3}e^{\tau^2/6} \mathsf{T}_2\left(\zeta^{\sigma_3/3}\right) \left\{
\begin{array}{ll}
\left( \begin{matrix} -e^{4\pi\mathrm{i}/3}  & 1 & e^{2\pi\mathrm{i}/3} \smallskip \\ 1 & -1 & -1 \smallskip \\ -e^{2\pi\mathrm{i}/3} & 1 & e^{4\pi\mathrm{i}/3} \end{matrix} \right), & \im(\zeta)>0, \medskip \\
\left( \begin{matrix} e^{2\pi\mathrm{i}/3}  & 1 & e^{4\pi\mathrm{i}/3} \smallskip \\ -1 & -1 & -1 \smallskip \\ e^{4\pi\mathrm{i}/3} & 1 & e^{2\pi\mathrm{i}/3} \end{matrix} \right), & \im(\zeta)<0,
\end{array}
\right.
\]
\[
\boldsymbol B := \left\{
\begin{array}{ll}
\diag\left(e^{4\pi\mathrm{i}/3}, 1, e^{2\pi\mathrm{i}/3}\right), & \im(\zeta)>0, \medskip \\
\diag\left(e^{2\pi\mathrm{i}/3}, 1, e^{4\pi\mathrm{i}/3}\right), & \im(\zeta)<0,
\end{array}
\right.
\]
and
\[
\boldsymbol C := \left\{
\begin{array}{ll}
\left( \begin{matrix} 0 & e^{4\pi\mathrm{i}/3}-e^{2\pi\mathrm{i}/3} & 1- e^{2\pi\mathrm{i}/3} \smallskip \\ e^{4\pi\mathrm{i}/3}-1 & 0 & 1-e^{2\pi\mathrm{i}/3} \smallskip \\ 1-e^{4\pi\mathrm{i}/3} & e^{4\pi\mathrm{i}/3}-e^{2\pi\mathrm{i}/3} & 0 \end{matrix} \right), & \im(\zeta)>0, \medskip \\
\left( \begin{matrix} 0 & e^{4\pi\mathrm{i}/3}-e^{2\pi\mathrm{i}/3} & e^{4\pi\mathrm{i}/3}-1 \smallskip \\ 1-e^{2\pi\mathrm{i}/3} & 0 & 1-e^{4\pi\mathrm{i}/3} \smallskip \\ e^{2\pi\mathrm{i}/3}-1 & e^{2\pi\mathrm{i}/3}-e^{4\pi\mathrm{i}/3} & 0 \end{matrix} \right) & \im(\zeta)<0.
\end{array}
\right.
\]
\end{itemize}

Below, we explain how to solve \hyperref[rhp2]{\rhp$_{1/\sqrt2}$} using $\boldsymbol\Upsilon$, a solution of \hyperref[rhp2]{\rhp$_{-1/\sqrt2}$} can be constructed analogously.

To carry \hyperref[rhup]{\rhup} into a neighborhood of $1/\sqrt2$, define
\[
\zeta(z) := \left(\frac23\int_{1/\sqrt2}^z\left(x-\frac1{\sqrt 2}\right)^{-1/3}\frac{H^{1/3}(x)}{\sqrt{1-x^2}}\, \mathrm{d}x\right)^{3/2},
\]
where the function $H$ was introduced in \eqref{hxi} and we choose the principle branch of the square root. Then $\zeta$ is conformal in some neighborhood of $1/\sqrt2$, $\zeta(1/\sqrt2)=0$, and $\zeta(x)>0$, $x>1/\sqrt2$. Further, define
\[
\tau(z) := -\zeta^{-1/3}(z)\int_{1/\sqrt2}^z\left(x-\frac1{\sqrt 2}\right)^{1/3}\frac{H^{-1/3}(x)}{\sqrt{1-x^2}}\,\mathrm{d}x.
\]
It readily follows that $\tau(z)$ is also conformal in some neighborhood of $1/\sqrt2$, $\tau(1/\sqrt2)=0$, and that
\begin{equation}
\label{tailDn}
\exp\left\{ -\frac32\zeta^{2/3}(z)\boldsymbol B^2-\tau(z)\zeta^{1/3}(z)\boldsymbol B\right\} = \boldsymbol D^n(z)
\end{equation}
by \eqref{hxi} and \eqref{Htoh}, and since $\Phi$ has value 1 at the point of $\RS$ whose natural projection is $1/\sqrt2$. Set
\[
\boldsymbol E_*(z) := \boldsymbol M(z)\diag\big(1,\rho_1^{-1}(z),\rho_2^{-1}(z)\big)\boldsymbol A^{-1}\big(n^{3/2}\zeta(z);n^{1/2}\tau(z)\big).
\]
It can easily be verified that
\[
\boldsymbol A_+ = \boldsymbol A_-\left\{
\begin{array}{lll}
\mathsf{T}_2\left(\begin{matrix} 0 & 1 \\ -1 & 0 \end{matrix}\right) & \text{on} & (0,\infty), \medskip \\
\mathsf{T}_3\left(\begin{matrix} 0 & 1 \\ -1 & 0 \end{matrix}\right) & \text{on} & (-\infty,0).
\end{array}
\right.
\]
As the entries of $\boldsymbol M$ as well as the entries of $\boldsymbol A^{-1}$ can have at most cubic root singularity at $1/\sqrt2$, the matrix $\boldsymbol E_*$ is holomorphic in $U_{1/\sqrt2}$. Then
\[
\boldsymbol P_*(z) := \boldsymbol E_*(z)\boldsymbol\Upsilon\big(n^{3/2}\zeta(z);n^{1/2}\tau(z)\big)\diag\big(1,\rho_1(z),\rho_2(z)\big)
\]
satisfies \hyperref[rhp2]{\rhp$_{1/\sqrt2}$}(a--c) (we always can adjust $\Delta_{i\pm}$ so that $\zeta$ maps them into $\{\re(z)=\pm\im(z)\}$). It also follows from \hyperref[rhup]{\rhup}(d) and \eqref{tailDn} that
\[
\boldsymbol P_*(z) = \boldsymbol M(z)\left(\boldsymbol I + \boldsymbol F_n(z) + \boldsymbol{\mathcal O}\big(n^{-2/3}\big) \right)\boldsymbol D^n(z)
\]
as $n\to\infty$ uniformly on $\partial U_{1/\sqrt2}$ (recall that $\tau(z)$ is conformal in $U_{1/\sqrt2}$ and vanishes at $1/\sqrt2$, which implies that $n^{1/2}\tau(z)$ remains bounded as $n\to\infty$ and therefore \hyperref[rhup]{\rhup}(d) is applicable), where
\[
\boldsymbol F_n(z) : = n^{-1/2}\zeta^{-1/3}(z) \diag\big(1,\rho_1^{-1}(z),\rho_2^{-1}(z)\big) \boldsymbol\Upsilon_1\big(n^{1/2}\tau(z)\big) \diag\big(1,\rho_1(z),\rho_2(z)\big).
\]
Since $\boldsymbol\Upsilon_1\big(n^{1/2}\tau(z)\big)\sim n^{1/2}\tau(z)$ as $z\to1/\sqrt2$ or $n\to\infty$, we can write
\[
(\boldsymbol M\boldsymbol F_n\boldsymbol M^{-1})(z) = \frac{\tau(z)\zeta^{-1/3}(z)}{(z-1/\sqrt2)^{2/3}}\left(\frac43\boldsymbol M_+\boldsymbol F_n^*\widehat{\boldsymbol M}_+^\mathsf{T} + \boldsymbol{\mathcal O}\left(\big|z-1/\sqrt2\big|^{1/3}\right)\right)
\]
by \eqref{Minverse} and \eqref{M-local}, where $\boldsymbol F_n^*$ is a matrix of constants (more precisely, it is one matrix of constants in $\{\im(z)>0\}$ and another one in $\{\im(z)<0\}$, see the definition of $\boldsymbol\Upsilon_1$) with entries whose moduli are uniformly bounded  as $n\to\infty$. As the fraction in the above equality is holomorphic around $1/\sqrt2$ and has value $-3/4$ at $1/\sqrt2$, we get that
\[
(\boldsymbol M\boldsymbol F_n\boldsymbol M^{-1})(z) = -\boldsymbol M_+\boldsymbol F_n^*\widehat{\boldsymbol M}_+^\mathsf{T} + \boldsymbol{\mathcal O}\left(\big|z-1/\sqrt2\big|^{1/3}\right).
\]
Since  $\big|z-1/\sqrt2\big|=n^{-1/2}$ on $\partial U_{1/\sqrt2}$, we have that
\[
\boldsymbol P_*(z) = \left(\boldsymbol I - \boldsymbol M_+\boldsymbol F_n^*\widehat{\boldsymbol M}_+^\mathsf{T} + \boldsymbol{\mathcal O}\big(n^{-1/6}\big) \right)\boldsymbol M(z)\boldsymbol D^n(z)
\]
uniformly on $\partial U_{1/\sqrt2}$. Thus, according to \eqref{zero-product}, \hyperref[rhp2]{\rhp$_{1/\sqrt2}$} is solved by
\[
\boldsymbol P_{1/\sqrt2}(z) := \left(\boldsymbol I + \boldsymbol M_+\boldsymbol F_n^*\widehat{\boldsymbol M}_+^\mathsf{T}\right)\boldsymbol P_*(z).
\]

\subsection{Final R-H Problem}

The final Riemann-Hilbert problem is \hyperref[rhr]{\rhr} from Section~\ref{ssec:FR-H} with ${\boldsymbol J}_{\boldsymbol Z}$ defined in \eqref{J-Z2}.
\begin{figure}[!ht]
\centering
\includegraphics[scale=.5]{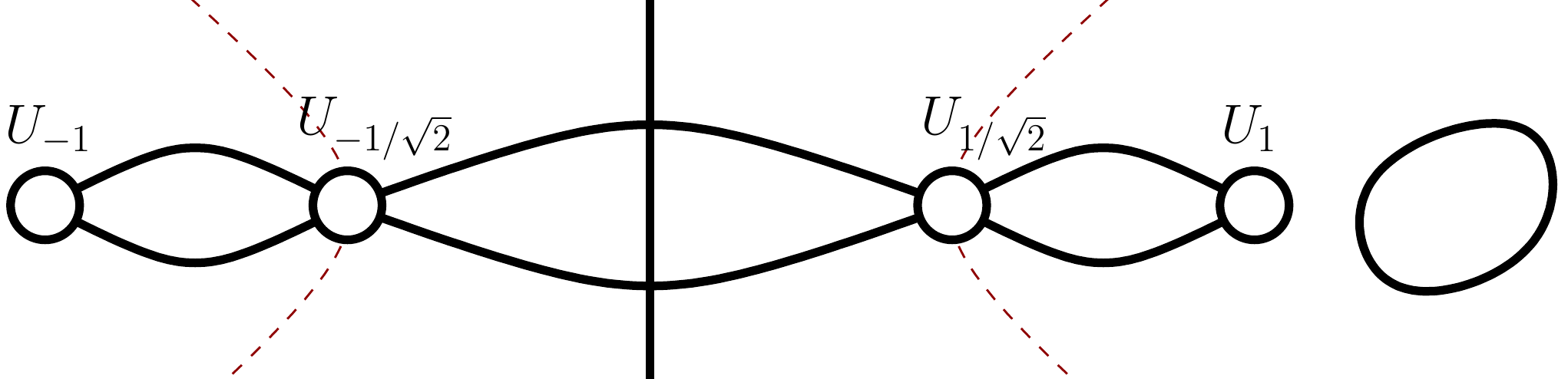}
\caption{\small The contour $\Sigma_{\boldsymbol R}$ without the circle $\{|z|=R\}$ (solid lines) and the borders of the domains $\Omega_{ijk}$ (dashed lines), see Figure~\ref{Oms}(b).}
\label{fig:Sigma_R-II}
\end{figure}
The same analysis shows that the jump matrices in \hyperref[rhr]{\rhr} are of order $\boldsymbol{\mathcal O}\big(n^{-1/6}\big)$ for $\N_*\ni n\to\infty$ and therefore \hyperref[rhr]{\rhr} is solvable for all $n\in\N_*$ large enough and satisfies \eqref{R-asymp} with $n^{-1}$ replaced by $n^{-1/6}$.

\subsection{Asymptotics of Hermite-Pad\'e Approximants}

It can readily be verified that formulae \eqref{Qn-Final}, \eqref{Rn1-Final}, and \eqref{Rn2-Final} remain valid. Let $\boldsymbol R$ be the solution of \hyperref[rhr]{\rhr}. Then the solution of  \hyperref[rhz3]{\rhz} is given by \eqref{Z-solution}, where $e\in\big\{\pm1,\pm1/\sqrt2,\infty\big\}$, for all $n\in\N_*$ large enough. The same argument as before shows that \eqref{bootstrap} holds in this case as well.

Since the first column of $\boldsymbol P_\infty$ is the same as the first column of $\boldsymbol{MD}^n$, \eqref{Z-Exact} remains valid when $k=0$. Thus, the proof of the first parts of \eqref{main-A} and \eqref{main-B} is exactly the same. The second parts of \eqref{main-A} and \eqref{main-B} are claimed to hold uniformly on each compact subset of the respective domains (because $\Gamma_i$ contains the point at infinity and $\mathcal N_\delta^{(i)}$ is bounded). Hence, given a compact subset, we always can enlarge $R$ in \hyperref[rhpi]{\rhp$_\infty$} so that this set is contained in $\{|z|<R\}$. This way \eqref{Z-Exact} is valid on this compact and the remaining part of the proof is the same as in Case I.

\section{Riemann-Hilbert Analysis: Case III}
\label{sec:CaseIII}

\subsection{Global Lenses}

\begin{figure}[!ht]
\centering
\includegraphics[scale=.5]{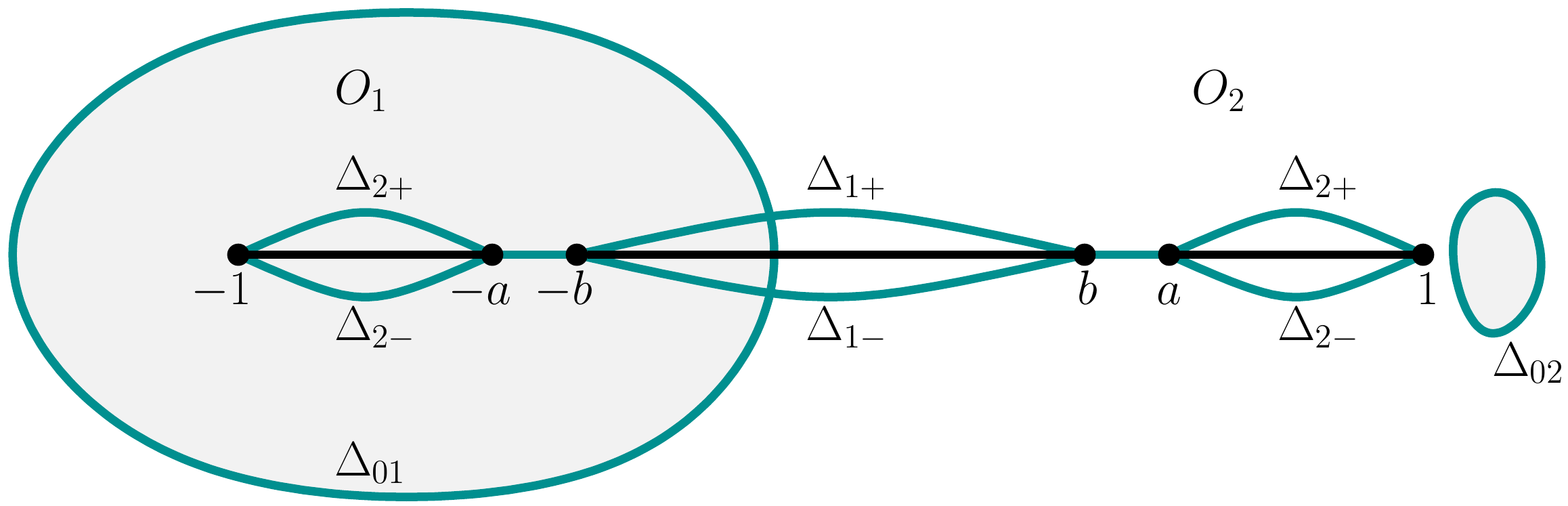}
\caption{\small The lens $\Sigma_{\boldsymbol Z}$, the domains $O_1$ (shaded region on the left) and $O_1$ (unshaded unbounded region), local lenses $\Delta_{1\pm}$ and $\Delta_{2\pm}$. The curve $\Delta_{02}$ is present only in Case~IIIa.}
\label{fig:Global-III}
\end{figure}

Let $\boldsymbol G_1(u)$ and $\boldsymbol G_2(v)$ be as in \eqref{global-matr}. Fix a domain, say $O_1$, that contains $[-1,-b]$ and whose boundary $\Delta_{01}:=\partial O_1$ is smooth, lies entirely in $\Omega_{012}$ (except for the point where it crosses $(-b,b)$), see Figure~\ref{Oms}(c), while crossing the real line at the origin, see Figure~\ref{fig:Global-III}. In Case IIIa, denote by $\Delta_{02}$ a smooth Jordan curve lying within the right component of $\Omega_{021}$ and containing within its interior all the singularities of $\rho_2/\rho_1$ (recall that $\rho_2/\rho_1\equiv\mathsf{const}$ in Case IIIb by Condition~\ref{cond}), see Figure~\ref{fig:Global-III}.
We also denote by $O_2$ the intersection of the exterior domains of $\Delta_{01}$ and $\Delta_{02}$.
We orient $\Delta_{01}$ and $\Delta_{02}$ counter-clockwise and set
\begin{equation}
\label{S3}
\boldsymbol S := \boldsymbol G_2\big((\rho_2/\rho_1)(\infty)\big)\boldsymbol Y\left\{
\begin{array}{lcl}
\boldsymbol G_2(-\rho_2/\rho_1)\boldsymbol G_1(\rho_1/\rho_2) & \text{in} & O_1, \medskip \\
\boldsymbol G_2(-\rho_2/\rho_1) & \text{in} & O_2, \medskip \\
\boldsymbol I & \text{in} & \overline\C\setminus\left(\overline{O_1\cup O_2}\right).
\end{array}
\right.
\end{equation}

Put $\Sigma_{\boldsymbol S}:=[-1,1]\cup\Delta_{01}\cup\Delta_{02}$. Then, according to \eqref{global-matr-modif}, $\boldsymbol S$ solves \rhs:
\begin{itemize}
\label{rhs3}
\item[(a)] $\boldsymbol S$ is analytic in $\C\setminus\Sigma_{\boldsymbol S}$ and $\lim_{z\to\infty} {\boldsymbol S}(z)\diag\left(z^{-2n},z^n,z^n\right) = \boldsymbol I$;
\item[(b)] ${\boldsymbol S}$ has continuous traces on $\Sigma_{\boldsymbol S}^\circ:=\Sigma_{\boldsymbol S}\setminus\{\pm1,\pm a,0\}$ that satisfy
\[
{\boldsymbol S}_+ = {\boldsymbol S}_-
\left\{
\begin{array}{lll}
{\boldsymbol J}(\rho_1,0) & \text{on} & \Delta_1^\circ\setminus\{0\}, \medskip \\
{\boldsymbol J}(0,\rho_2^*) & \text{on} & \Delta_2^\circ, \medskip \\
{\boldsymbol G}_1(\rho_1/\rho_2) & \text{on} & \Delta_{01}, \medskip \\
{\boldsymbol G}_2(\rho_2/\rho_1) & \text{on} & \Delta_{02},
\end{array}
\right.
\]
where $\rho_2^*$ is defined by \eqref{rhostar} and $\Delta_i$ in \eqref{chains};
\item[(c)] $\boldsymbol S$ satisfies \hyperref[rhy]{\rhy}(c) (see Section \ref{MOP}) with $[-1,1]$ replaced by $\Sigma_{\boldsymbol S}$.
\end{itemize}
If  \hyperref[rhs]{\rhs} is solvable, then so is  \hyperref[rhy]{\rhy}, and the solutions are connected via \eqref{S3}.

\subsection{Local Lenses}

As usual, we introduce additional arcs $\Delta_{1\pm}$ and systems of two arcs $\Delta_{2\pm}$ as in Figure~\ref{fig:Global-III}, all oriented from left to right. We further denote by $O_{i\pm}$ the open sets bounded by $\Delta_i$ and the arcs $\Delta_{i\pm}$, $i\in\{1,2\}$. Set
\begin{equation}
\label{Z3}
\boldsymbol Z:= \boldsymbol S\boldsymbol L_i^{\mp1} \quad \text{in} \quad O_{i\pm},
\end{equation}
where the matrices $\boldsymbol L_i$ are defined by \eqref{local-matr}. Put $\Sigma_{\boldsymbol Z} := \Sigma_{\boldsymbol S}\cup\Delta_{1+}\cup\Delta_{1-}\cup\Delta_{2+}\cup\Delta_{2-}$. Then it can readily be checked that ${\boldsymbol Z}$ solves \rhz:
\begin{itemize}
\label{rhz3}
\item[(a)] ${\boldsymbol Z}$ is analytic in $\C\setminus\Sigma_{\boldsymbol Z}$ and $\lim_{z\to\infty} {\boldsymbol Z}(z)\diag\left(z^{2n},z^{-n},z^{-n}\right) = \boldsymbol I$;
\item[(b)] ${\boldsymbol Z}$ has continuous traces on each side of $\Sigma_{\boldsymbol Z}^\circ:=\Sigma_{\boldsymbol Z}\setminus\{\pm1,\pm a,\pm b,0\}$ that satisfy
\[
{\boldsymbol Z}_+ = {\boldsymbol Z}_-
\left\{
\begin{array}{lll}
{\boldsymbol J}_i & \text{on} & \Delta_i^\circ\setminus\{0\}, \quad i\in\{1,2\}, \medskip \\
\boldsymbol J_{\boldsymbol Z} &\text{on} & \Sigma_{\boldsymbol Z}^\circ\setminus(\Delta_1\cup\Delta_2),
\end{array}
\right.
\]
where $\boldsymbol J_i$ are defined by \eqref{local-jumps} and
\begin{equation}
\label{J-Z3}
\boldsymbol J_{\boldsymbol Z} :=
\left\{
\begin{array}{lll}
{\boldsymbol J}(\rho_1,0) & \text{on} & (-a,-b)\cup(b,a), \medskip \\
{\boldsymbol G}_1(\rho_1/\rho_2) & \text{on} & \Delta_{01}, \medskip \\
{\boldsymbol G}_2(\rho_2/\rho_1) & \text{on} & \Delta_{02}, \medskip \\
{\boldsymbol L}_i & \text{on} & \Delta_{i\pm}, \quad i\in\{1,2\};
\end{array}
\right.
\end{equation}
\item[(c)] $\boldsymbol Z$ satisfies \hyperref[rhs3]{\rhs}(c) with $\Sigma_{\boldsymbol S}$ replaced by $\Sigma_{\boldsymbol Z}$.
\end{itemize}
As before, if \hyperref[rhz3]{\rhz} is solvable, then so is \hyperref[rhs3]{\rhs}, and the solutions are connected via \eqref{Z3}.

\subsection{Global Parametrix}

Let $\boldsymbol N=\boldsymbol C^{-n}\boldsymbol {MD}^n$ be given by \eqref{N}. Then $\det(\boldsymbol N)\equiv1$ and it is a solution of the following Riemann-Hilbert problem (\rhn):
\begin{itemize}
\label{rhn3}
\item[(a)] ${\boldsymbol N}$ is analytic in $\C\setminus(\Delta_1\cup\Delta_2)$ and $\lim_{z\to\infty} {\boldsymbol N}(z)\diag\left(z^{-2n},z^n,z^n\right) = \boldsymbol I$;
\item[(b)] ${\boldsymbol N}$ has continuous traces on each side of $\Delta_i^\circ$ that satisfy $\boldsymbol N_+ = \boldsymbol N_- {\boldsymbol J}_i$, $i\in\{1,2\}$.
\end{itemize}

\subsection{Local Parametrices}

Again, as in the previous cases, we need to solve \hyperref[rhp]{\rhp$_e$} for $e\in\{\pm1,\pm a,\pm b\}$. In fact, these local problems are exactly the same as in Case I. Thus, their solutions were constructed in Section~\ref{ssec:local-param}.

\subsection{Final R-H Problem}

Once more, the final Riemann-Hilbert problem is \hyperref[rhr]{\rhr} from Section~\ref{ssec:FR-H} with ${\boldsymbol J}_{\boldsymbol Z}$ defined in \eqref{J-Z3} and $e\in\{\pm a,\pm b,\pm 1\}$.
\begin{figure}[!ht]
\centering
\includegraphics[scale=.5]{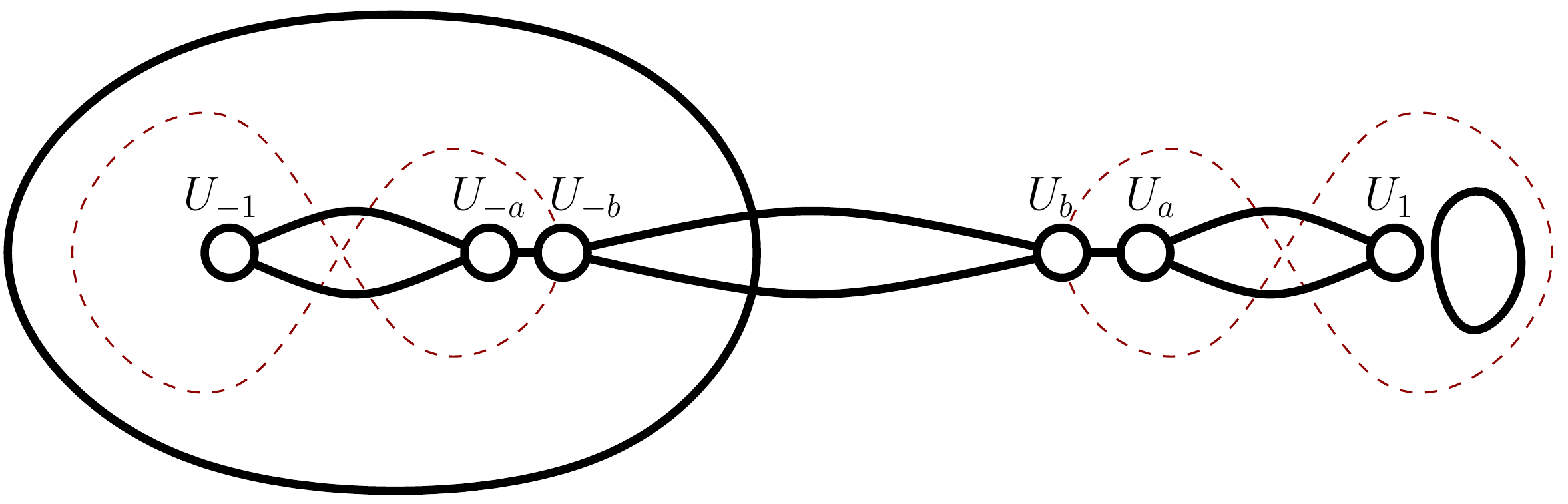}
\caption{\small The contour $\Sigma_{\boldsymbol R}$ (solid lines) and the relevant borders of the domains $\Omega_{ijk}$ (dashed lines), see Figure~\ref{Oms}(c,d).}
\label{fig:Sigma_R-III}
\end{figure}
Exactly the same analysis shows that the jump matrices in \hyperref[rhr]{\rhr} are of order $\boldsymbol{\mathcal{O}}(1/n)$ for $n\in\N_*$ and therefore \hyperref[rhr]{\rhr} is solvable for all $n\in\N_*$ large enough and satisfies \eqref{R-asymp}.

\subsection{Asymptotics of Hermite-Pad\'e Approximants}

As in Case I, one can verify that the formulae \eqref{Qn-Final}, \eqref{Rn1-Final}, \eqref{Rn2-Final}, \eqref{Z-Exact}, and \eqref{bootstrap} remain valid in this case as well. So the proof of \eqref{main-A} and \eqref{main-B} proceeds exactly along the same lines as in Case I.


\begin{thebibliography}{10}

\bibitem{Ang19}
A.~Angelesco.
\newblock Sur deux extensions des fractions continues alg\'ebriques.
\newblock {\em Comptes Rendus de l'Academie des Sciences, Paris}, 18:262--265,
  1919.

\bibitem{Ap88}
A.I. Aptekarev.
\newblock Asymptotics of simultaneously orthogonal polynomials in the
  {A}ngelesco case.
\newblock {\em Mat. Sb.}, 136(178)(1):56--84, 1988.
\newblock English transl. in {\it {M}ath. {USSR} {S}b.} 64, 1989.

\bibitem{Ap99}
A.I. Aptekarev.
\newblock Strong asymptotics of multiple orthogonal polynomials for {N}ikishin
  systems.
\newblock {\em Mat. Sb.}, 190(5):3--44, 1999.
\newblock English transl. in {\it Sb. Math.} 190(5):631--669, 1999.

\bibitem{Ap08b}
A.I. Aptekarev.
\newblock Asymptotics of {H}ermite-{P}ad\'e approximants for a pair of
  functions with branch points.
\newblock {\em Dokl. Akad. Nauk}, 422(4):443--445, 2008.
\newblock English trans. in {\it Dokl. Math.} 78(2):717--719, 2008.

\bibitem{AptBlK05}
A.I. Aptekarev, P.M. Bleher, and A.B.J. Kuijlaars.
\newblock Large $n$ limit of {G}aussian random matrices with external source,
  {Part II}.
\newblock {\em Comm. Math. Phys.}, 259:367--389, 2005.

\bibitem{ApKalLysToul09}
A.I. Aptekarev, V.A. Kalyagin, V.G. Lysov, and D.N. Toulyakov.
\newblock Equilibrium of vector potentials and uniformization of the algebraic
  curves of genus $0$.
\newblock {\em J. Comput. Appl. Math.}, 233(602--616), 2009.

\bibitem{ApKVA08}
A.I. Aptekarev, A.B.J.~Kuijlaars, and W.~Van Assche.
\newblock Asymptotics of {H}ermite-{Pa}d\'e rational approximants for two
  analytic functions with separated pairs of branch points (case of genus 0).
\newblock {\em Int. Math. Res. Papers}, 2008, rpm007, 128 pp.

\bibitem{ApLy10}
A.I. Aptekarev and V.G. Lysov.
\newblock Systems of Markov functions generated by graphs and the asymptotics of their Hermite-Padé approximants.
\newblock {\em Mat. Sb.}, 201(2):29--78, 2010.
\newblock English transl. in {\em Sb. Math.}, 201:183--234, 2010.

\bibitem{AptSt92}
A.I. Aptekarev and H.~Stahl.
\newblock Asymptotics of {H}ermite-{P}ad\'e polynomials.
\newblock In A.A. Gonchar and E.B. Saff, editors, {\em Progress in
  Approximation Theory}, pages 127--167, Berlin, 1992. Springer-Verlag.

\bibitem{uApToulVA}
A.I. Aptekarev, D.N. Toulyakov, and W.~{Van Assche}.
\newblock Hyperelliptic uniformization of algebraic curves of the third order.
\newblock {\em J. Comput. Appl. Math.}, 284:38--49, 2015.

\bibitem{uApToulYa}
A.I. Aptekarev, D.N. Toulyakov, and M.L. Yattselev.
\newblock On uniformization of a certain cubic algebraic curve of genus 2.
\newblock \emph{manuscript}.

\bibitem{Deift}
P.~Deift.
\newblock {\em Orthogonal Polynomials and Random Matrices: a Riemann-Hilbert
  Approach}, volume~3 of {\em Courant Lectures in Mathematics}.
\newblock Amer. Math. Soc., Providence, RI, 2000.

\bibitem{DKMLVZ99b}
P.~Deift, T.~Kriecherbauer, K.T.-R. McLaughlin, S.~Venakides, and X.~Zhou.
\newblock Strong asymptotics for polynomials orthogonal with respect to varying
  exponential weights.
\newblock {\em Comm. Pure Appl. Math.}, 52(12):1491--1552, 1999.

\bibitem{DZ93}
P.~Deift and X.~Zhou.
\newblock A steepest descent method for oscillatory {R}iemann-{H}ilbert
  problems. {A}symptotics for the m{K}d{V} equation.
\newblock {\em Ann. of Math.}, 137:295--370, 1993.

\bibitem{DeK11}
K.~Deschout and A.B.J. Kuijlaars.
\newblock Double scaling limit for modified {J}acobi-{A}ngelesco polynomials.
\newblock In P.~Br\"and\'en, M.~Passare, and M.~Putinar, editors, {\em Notions
  of Positivity and the Geometry of Polynomials}, Trends in Mathematics, pages
  115--161, Basel, 2011. Springer.

\bibitem{DLMF}
F.W.J.~Olver et~al. editors.
\newblock NIST digital library of mathematical functions.
\newblock \url{http://dlmf.nist.gov}.

\bibitem{FPrLL11}
U.~{Fidalgo Prieto} and G.~L\'opez Lagomasino.
\newblock Nikishin systems are perfect.
\newblock {\em Constr. Approx.}, 34(3):297--356, 2011.

\bibitem{FIK91}
A.S. Fokas, A.R. Its, and A.V. Kitaev.
\newblock Discrete {P}anlev\'e equations and their appearance in quantum
  gravity.
\newblock {\em Comm. Math. Phys.}, 142(2):313--344, 1991.

\bibitem{FIK92}
A.S. Fokas, A.R. Its, and A.V. Kitaev.
\newblock The isomonodromy approach to matrix models in {2D} quantum
  gravitation.
\newblock {\em Comm. Math. Phys.}, 147(2):395--430, 1992.

\bibitem{Gakhov}
F.D. Gakhov.
\newblock {\em Boundary Value Problems}.
\newblock Dover Publications, Inc., New York, 1990.

\bibitem{GRakh81}
A.A. Gonchar and E.A. Rakhmanov.
\newblock On convergence of simultaneous {P}ad\'e approximants for systems of
  functions of {M}arkov type.
\newblock {\em Trudy Mat. Inst. Steklov}, 157:31--48, 1981.
\newblock English transl. in {\it {P}roc. {S}teklov {I}nst. {M}ath.} 157, 1983.

\bibitem{GRakhS97}
A.A. Gonchar, E.A. Rakhmanov, and V.N. Sorokin.
\newblock Hermite-{P}ad\'e approximants for systems of {M}arkov-type functions.
\newblock {\em Mat. Sb.}, 188(5):33--58, 1997.
\newblock English transl. in {\it {M}ath. {USSR} Sbornik} 188(5):671--696,
  1997.

\bibitem{Herm73}
C.~Hermite.
\newblock Sur la fonction exponentielle.
\newblock {\em C. R. Acad. Sci. Paris}, 77:18--24, 74--79, 226--233, 285--293,
  1873.

\bibitem{Kal79}
V.A. Kalyagin.
\newblock On a class of polynomials defined by two orthogonality relations.
\newblock {\em Mat. Sb.}, 110(4):609--627, 1979.

\bibitem{KMcLVAV04}
A.B.J. Kuijlaars, K.T.-R. McLaughlin, W.~Van Assche, and M.~Vanlessen.
\newblock The {R}iemann-{H}ilbert approach to strong asymptotics for orthogonal
  polynomials on $[-1,1]$.
\newblock {\em Adv. Math.}, 188(2):337--398, 2004.

\bibitem{Mar95}
A.A. Markov.
\newblock Deux d\'emonstrations de la convergence de certaines fractions
  continues.
\newblock {\em Acta Math.}, 19:93--104, 1895.

\bibitem{Nik79}
E.M. Nikishin.
\newblock A system of {M}arkov functions.
\newblock {\em Vestnik Moskovskogo Universiteta Seriya 1, Matematika
  Mekhanika}, 34(4):60--63, 1979.
\newblock Translated in \emph{Moscow University Mathematics Bulletin} 34(4),
  63--66, 1979.

\bibitem{Nik80}
E.M. Nikishin.
\newblock Simultaneous {P}ad\'e approximants.
\newblock {\em Mat. Sb.}, 113(155)(4):499--519, 1980.

\bibitem{Nut84}
J.~Nuttall.
\newblock Asymptotics of diagonal {H}ermite-{P}ad\'e polynomials.
\newblock {\em J. Approx. Theory}, 42(4):299--386, 1984.

\bibitem{Pade92}
H.~Pad\'e.
\newblock Sur la repr\'esentation approch\'ee d'une fonction par des fractions
  rationnelles.
\newblock {\em Ann. Sci Ecole Norm. Sup.}, 9(3):3--93, 1892.

\bibitem{Privalov}
I.I. Privalov.
\newblock {\em Boundary Properties of Analytic Functions}.
\newblock GITTL, Moscow, 1950.
\newblock German transl., VEB Deutscher Verlag Wiss., Berlin, 1956.

\bibitem{Rakh11}
E.A. Rakhmanov.
\newblock On the asymptotics of {H}ermite-{P}ad\'e polynomials for two
  {M}arkov-type functions.
\newblock {\em Mat. Sb.},  202(1):133--140, 2011.
\newblock English transl. In {\em Sb. Math.}, 202(1):127--134, 2011.

\bibitem{uSt2}
H.~Stahl.
\newblock Asymptotics of {H}ermite-{P}ad\'e polynomials and related
  approximants. {A} summary of results.
\newblock Manuscript, 79 pages.

\bibitem{St87}
H.~Stahl.
\newblock Asymptotics of {H}ermite-{P}ad\'e polynomials and related convergence
  results -- a summary of results.
\newblock In {\em Nonlinear numerical methods and rational approximation},
  volume~43 of {\em Math. Appl.}, pages 23--53. Reidel, Dordrecht, 1987.

\bibitem{Strebel}
K.~Strebel.
\newblock {\em Quadratic Differentials}, volume~5 of {\em Ergebnisse der
  Mathematik und ihrer Grenzgebiete (3)}.
\newblock Springer-Verlag, Berlin, 1984.

\bibitem{Toul09}
D.N. Tulyakov.
\newblock Difference schemes with power-growth bases perturbed by the spectral parameter.
\newblock {\em Mat. Sb.}, 200(5):129--158, 2009.
\newblock English transl. in {\em Sb. Math.}, 200(5):753--781, 2009.

\bibitem{GerKVA01}
W. Van Assche, J.S. Geronimo, and A.B.J. Kuijlaars.
\newblock Riemann-{Hilbert} problems for multiple orthogonal polynomials.
\newblock In {\em Special functions 2000: current perspective and future
  directions}, number~30 in NATO Sci. Ser. II Math. Phys. Chem., pages 23--59,
  Dordrecht, 2001. Kluwer Acad. Publ.

\bibitem{Zver71}
E.I. Zverovich.
\newblock Boundary value problems in the theory of analytic functions in
  {H}\"older classes on {R}iemann surfaces.
\newblock {\em Russian Math. Surveys}, 26(1):117--192, 1971.

\end{thebibliography}
\end{document}